\font\tenmib=cmmib10
\font\sevenmib=cmmib10 scaled 800
\font\titolo=cmbx12
\font\titolone=cmbx10 scaled\magstep 2

\font\cs=cmcsc10
\font\sc=cmcsc10
\font\css=cmcsc8

\font\ninerm=cmr9
\font\ottorm=cmr8
\textfont5=\tenmib
\scriptfont5=\sevenmib
\scriptscriptfont5=\fivei

\font\euftw=eufm9 scaled\magstep1
\font\euftww=eufm7 scaled\magstep1
\font\euftwww=eufm5 scaled\magstep1
\font\msytw=msbm9 scaled\magstep1

\font\msytwww=msbm5 scaled\magstep1

\font\indbf=cmbx10 scaled\magstep2

\def\st{\scriptstyle}

\font\ottorm=cmr8\font\ottoi=cmmi8\font\ottosy=cmsy8%
\font\ottobf=cmbx8\font\ottott=cmtt8%
\font\ottocss=cmcsc8%
\font\ottosl=cmsl8\font\ottoit=cmti8%
\font\sixrm=cmr6\font\sixbf=cmbx6\font\sixi=cmmi6\font\sixsy=cmsy6%
\font\fiverm=cmr5\font\fivesy=cmsy5
\font\fivei=cmmi5
\font\fivebf=cmbx5%

\def\ottopunti{\def\rm{\fam0\ottorm}%
\textfont0=\ottorm\scriptfont0=\sixrm\scriptscriptfont0=\fiverm%
\textfont1=\ottoi\scriptfont1=\sixi\scriptscriptfont1=\fivei%
\textfont2=\ottosy\scriptfont2=\sixsy\scriptscriptfont2=\fivesy%
\textfont3=\tenex\scriptfont3=\tenex\scriptscriptfont3=\tenex%
\textfont4=\ottocss\scriptfont4=\sc\scriptscriptfont4=\sc%
\textfont5=\tenmib\scriptfont5=\sevenmib\scriptscriptfont5=\fivei
\textfont\itfam=\ottoit\def\it{\fam\itfam\ottoit}%
\textfont\slfam=\ottosl\def\sl{\fam\slfam\ottosl}%
\textfont\ttfam=\ottott\def\tt{\fam\ttfam\ottott}%
\textfont\bffam=\ottobf\scriptfont\bffam=\sixbf%
\scriptscriptfont\bffam=\fivebf\def\bf{\fam\bffam\ottobf}%
\setbox\strutbox=\hbox{\vrule height7pt depth2pt width0pt}%
\normalbaselineskip=9pt\let\sc=\sixrm\normalbaselines\rm}
\let\nota=\ottopunti%

%
%
%
%
%
%
%

\global\newcount\numsec\global\newcount\numapp
\global\newcount\numfor\global\newcount\numfig
\global\newcount\numsub
\numsec=0\numapp=0\numfig=0
\def\veroparagrafo{\number\numsec}\def\veraformula{\number\numfor}
\def\veraappendice{\number\numapp}\def\verasub{\number\numsub}
\def\verafigura{\number\numfig}

\def\section(#1,#2){\advance\numsec by 1\numfor=1\numsub=1\numfig=1%
\SIA p,#1,{\veroparagrafo} %
\write15{\string\Fp (#1){\secc(#1)}}%
\write16{ sec. #1 ==> \secc(#1)  }%
\hbox to \hsize{\titolo\hfill \number\numsec. #2\hfill%
\expandafter{\alato(sec. #1)}}\*}

\def\appendix(#1,#2){\advance\numapp by 1\numfor=1\numsub=1\numfig=1%
\SIA p,#1,{A\veraappendice} %
\write15{\string\Fp (#1){\secc(#1)}}%
\write16{ app. #1 ==> \secc(#1)  }%
\hbox to \hsize{\titolo\hfill Appendix A\number\numapp. #2\hfill%
\expandafter{\alato(app. #1)}}\*}

\def\senondefinito#1{\expandafter\ifx\csname#1\endcsname\relax}

\def\SIA #1,#2,#3 {\senondefinito{#1#2}%
\expandafter\xdef\csname #1#2\endcsname{#3}\else
\write16{???? ma #1#2 e' gia' stato definito !!!!} \fi}

\def \Fe(#1)#2{\SIA fe,#1,#2 }
\def \Fp(#1)#2{\SIA fp,#1,#2 }
\def \Fg(#1)#2{\SIA fg,#1,#2 }

\def\etichetta(#1){(\veroparagrafo.\veraformula)%
\SIA e,#1,(\veroparagrafo.\veraformula) %
\global\advance\numfor by 1%
\write15{\string\Fe (#1){\equ(#1)}}%
\write16{ EQ #1 ==> \equ(#1)  }}

\def\etichettaa(#1){(A\veraappendice.\veraformula)%
\SIA e,#1,(A\veraappendice.\veraformula) %
\global\advance\numfor by 1%
\write15{\string\Fe (#1){\equ(#1)}}%
\write16{ EQ #1 ==> \equ(#1) }}

\def\getichetta(#1){\veroparagrafo.\verafigura%
\SIA g,#1,{\veroparagrafo.\verafigura} %
\global\advance\numfig by 1%
\write15{\string\Fg (#1){\graf(#1)}}%
\write16{ Fig. #1 ==> \graf(#1) }}

\def\etichettap(#1){\veroparagrafo.\verasub%
\SIA p,#1,{\veroparagrafo.\verasub} %
\global\advance\numsub by 1%
\write15{\string\Fp (#1){\secc(#1)}}%
\write16{ par #1 ==> \secc(#1)  }}

\def\etichettapa(#1){A\veraappendice.\verasub%
\SIA p,#1,{A\veraappendice.\verasub} %
\global\advance\numsub by 1%
\write15{\string\Fp (#1){\secc(#1)}}%
\write16{ par #1 ==> \secc(#1)  }}

\def\Eq(#1){\eqno{\etichetta(#1)\alato(#1)}}
\def\eq(#1){\etichetta(#1)\alato(#1)}
\def\Eqa(#1){\eqno{\etichettaa(#1)\alato(#1)}}
\def\eqa(#1){\etichettaa(#1)\alato(#1)}
\def\eqg(#1){\getichetta(#1)\alato(fig. #1)}
\def\sub(#1){\0\palato(p. #1){\bf \etichettap(#1).}}
\def\asub(#1){\0\palato(p. #1){\bf \etichettapa(#1).}}

\def\equv(#1){\senondefinito{fe#1}$\clubsuit$#1%
\write16{eq. #1 non e' (ancora) definita}%
\else\csname fe#1\endcsname\fi}
\def\grafv(#1){\senondefinito{fg#1}$\clubsuit$#1%
\write16{fig. #1 non e' (ancora) definito}%
\else\csname fg#1\endcsname\fi}
\def\secv(#1){\senondefinito{fp#1}$\clubsuit$#1%
\write16{par. #1 non e' (ancora) definito}%
\else\csname fp#1\endcsname\fi}

\def\equ(#1){\senondefinito{e#1}\equv(#1)\else\csname e#1\endcsname\fi}
\def\graf(#1){\senondefinito{g#1}\grafv(#1)\else\csname g#1\endcsname\fi}
\def\figura(#1){{\css Figura} \getichetta(#1)}
\def\secc(#1){\senondefinito{p#1}\secv(#1)\else\csname p#1\endcsname\fi}
\def\sec(#1){{\S\secc(#1)}}
\def\refe(#1){{[\secc(#1)]}}

\def\BOZZA{\bz=1
\def\alato(##1){\rlap{\kern-\hsize\kern-1.2truecm{$\scriptstyle##1$}}}
\def\palato(##1){\rlap{\kern-1.2truecm{$\scriptstyle##1$}}}
}

\def\alato(#1){}
\def\galato(#1){}
\def\palato(#1){}


{\count255=\time\divide\count255 by 60 \xdef\hourmin{\number\count255}
        \multiply\count255 by-60\advance\count255 by\time
   \xdef\hourmin{\hourmin:\ifnum\count255<10 0\fi\the\count255}}

\def\oramin{\hourmin }

\def\data{\number\day/\ifcase\month\or gennaio \or febbraio \or marzo \or
aprile \or maggio \or giugno \or luglio \or agosto \or settembre
\or ottobre \or novembre \or dicembre \fi/\number\year;\ \oramin}
\setbox200\hbox{$\scriptscriptstyle \data $}

\newdimen\xshift \newdimen\xwidth \newdimen\yshift \newdimen\ywidth

\def\ins#1#2#3{\vbox to0pt{\kern-#2\hbox{\kern#1 #3}\vss}\nointerlineskip}

\def\eqfig#1#2#3#4#5{
\par\xwidth=#1 \xshift=\hsize \advance\xshift
by-\xwidth \divide\xshift by 2
\yshift=#2 \divide\yshift by 2
\line{\hglue\xshift \vbox to #2{\vfil
#3 \includegraphics{#4.ps}
}\hfill\raise\yshift\hbox{#5}}}

\def\8{\write12}  


\let\a=\alpha \let\b=\beta  \let\g=\gamma  \let\d=\delta \let\e=\varepsilon
  \let\h=\eta   \let\th=\theta \let\k=\kappa \let\l=\lambda
    \let\n=\nu    \let\x=\xi     \let\p=\pi    \let\r=\rho
\let\s=\sigma \let\t=\tau   \let\f=\varphi 
\let\ps=\psi   \let\o=\omega
\let\G=\Gamma \let\D=\Delta  \let\Th=\Theta\let\L=\Lambda 
\let\P=\Pi         
\let\O=\Omega 

\def\\{\hfill\break} \let\==\equiv

\let\io=\infty 
\def\ap{{\it a priori\ }}
\let\0=\noindent

\def\ie{{i.e. }}\def\eg{{e.g. }}
\let\dpr=\partial

\def\tende#1{\,\vtop{\ialign{##\crcr\rightarrowfill\crcr
 \noalign{\kern-1pt\nointerlineskip} \hskip3.pt${\scriptstyle
 #1}$\hskip3.pt\crcr}}\,}
\def\circage{\lower2pt\hbox{$\,\buildrel > \over {\scriptstyle \sim}\,$}}
\def\otto{\,{\kern-1.truept\leftarrow\kern-5.truept\to\kern-1.truept}\,}
\def\fra#1#2{{#1\over#2}}

\def\EE{{\cal E}}\def\MM{{\cal M}} \def\VV{{\cal V}}
\def\FF{{\cal F}} \def\HHH{{\cal H}}
\def\NN{{\cal N}} 
\def\RR{{\cal R}}  
\def\AAA{{\cal A}} \def\SS{{\cal S}}

\def\T#1{{#1_{\kern-3pt\lower7pt\hbox{$\widetilde{}$}}\kern3pt}}
\def\VVV#1{{\underline #1}_{\kern-3pt
\lower7pt\hbox{$\widetilde{}$}}\kern3pt\,}
\def\W#1{#1_{\kern-3pt\lower7.5pt\hbox{$\widetilde{}$}}\kern2pt\,}

\def\lis{\overline}\def\tto{\Rightarrow}

\def\indica{\leaders \hbox to 0.5cm{\hss.\hss}\hfill}
\def\guida{\leaders\hbox to 1em{\hss.\hss}\hfill}

\def\hh{{\bf h}}  \def\AA{{\bf A}} \def\qq{{\bf q}}
\def\BB{{\bf B}}   \def\pp{{\bf p}}
   
\def\aaa{{\bf a}}\def\bbb{{\bf b}}\def\hhh{{\bf h}}\def\III{{\bf I}}

\def\ul{\underline}


\mathchardef\aa   = "050B
\mathchardef\bb   = "050C
\mathchardef\ggg  = "050D
\mathchardef\xxx  = "0518
\mathchardef\zzzzz= "0510
\mathchardef\oo   = "0521
\mathchardef\lll  = "0515
\mathchardef\mm   = "0516
\mathchardef\Dp   = "0540
\mathchardef\H    = "0548
\mathchardef\FFF  = "0546
\mathchardef\ppp  = "0570
\mathchardef\nn   = "0517
\mathchardef\ff   = "0527
\mathchardef\pps  = "0520
\mathchardef\FFF  = "0508
\mathchardef\nnnnn= "056E

\def\to{\rightarrow}

\def\qed{\hfill\raise1pt\hbox{\vrule height5pt width5pt depth0pt}}

\def\Val{{\rm Val}}
\def\indic{\hbox{\raise-2pt \hbox{\indbf 1}}}

\def\RRR{\hbox{\msytw R}}

\def\NNN{\hbox{\msytw N}} 
 \def\ZZZ{\hbox{\msytw Z}}
 \def\zzz{\hbox{\msytwww Z}}
\def\TTT{\hbox{\msytw T}}

\def\vvv{\hbox{\euftw v}}    \def\vvvv{\hbox{\euftww v}}
\def\vvvvv{\hbox{\euftwww v}}\def\www{\hbox{\euftw w}}
\def\wwww{\hbox{\euftww w}}  \def\wwwww{\hbox{\euftwww w}}
\def\vvr{\hbox{\euftw r}}    

\def\ul#1{{\underline#1}}
\def\Sqrt#1{{\sqrt{#1}}}
\def\V0{{\bf 0}}
\def\defi{\,{\buildrel def\over=}\,}
\def\lhs{{\it l.h.s.}\ }


\newcount\mgnf  
\mgnf=0 

\ifnum\mgnf=0
\def\openone{\leavevmode\hbox{\ninerm 1\kern-3.3pt\tenrm1}}%
\def\*{\vglue0.3truecm}\fi
\ifnum\mgnf=1
\def\openone{\leavevmode\hbox{\ninerm 1\kern-3.63pt\tenrm1}}%
\def\*{\vglue0.5truecm}\fi


\newcount\tipobib\newcount\bz\bz=0\newcount\aux\aux=1
\newdimen\bibskip\newdimen\maxit\maxit=0pt


\tipobib=2
\def\9#1{\ifnum\aux=1#1\else\relax\fi}

\newwrite\bib
\immediate\openout\bib=\jobname.bib
\global\newcount\bibex
\bibex=0
\def\verabib{\number\bibex}

\ifnum\tipobib=0
\def\cita#1{\expandafter\ifx\csname c#1\endcsname\relax
\hbox{$\clubsuit$}#1\write16{Manca #1 !}%
\else\csname c#1\endcsname\fi}
\def\rife#1#2#3{\immediate\write\bib{\string\raf{#2}{#3}{#1}}
\immediate\write15{\string\C(#1){[#2]}}
\setbox199=\hbox{#2}\ifnum\maxit < \wd199 \maxit=\wd199\fi}
\fi
\ifnum\tipobib=1
\def\cita#1{%
\expandafter\ifx\csname d#1\endcsname\relax%
\expandafter\ifx\csname c#1\endcsname\relax%
\hbox{$\clubsuit$}#1\write16{Manca #1 !}%
\else\probib(ref. numero )(#1)%
\csname c#1\endcsname%
\fi\else\csname d#1\endcsname\fi}%
\def\rife#1#2#3{\immediate\write15{\string\Cp(#1){%
\string\immediate\string\write\string\bib{\string\string\string\raf%
{\string\verabib}{#3}{#1}}%
\string\Cn(#1){[\string\verabib]}%
\string\CCc(#1)%
}}}%
\fi
\ifnum\tipobib=2%
\def\cita#1{\expandafter\ifx\csname c#1\endcsname\relax
\hbox{$\clubsuit$}#1\write16{Manca #1 !}%
\else\csname c#1\endcsname\fi}
\def\rife#1#2#3{\immediate\write\bib{\string\raf{#1}{#3}{#2}}
\immediate\write15{\string\C(#1){[#1]}}
\setbox199=\hbox{#2}\ifnum\maxit < \wd199 \maxit=\wd199\fi}
\fi

\def\Cn(#1)#2{\expandafter\xdef\csname d#1\endcsname{#2}}
\def\CCc(#1){\csname d#1\endcsname}
\def\probib(#1)(#2){\global\advance\bibex+1%
\9{\immediate\write16{#1\verabib => #2}}%
}

\def\C(#1)#2{\SIA c,#1,{#2}}
\def\Cp(#1)#2{\SIAnx c,#1,{#2}}

\def\SIAnx #1,#2,#3 {\senondefinito{#1#2}%
\expandafter\def\csname#1#2\endcsname{#3}\else%
\write16{???? ma #1,#2 e' gia' stato definito !!!!}\fi}

\bibskip=10truept
\def\hboxto{\hbox to}

\catcode`\{=12\catcode`\}=12
\catcode`\<=1\catcode`\>=2
\immediate\write\bib<
        \string\halign{\string\hboxto \string\maxit%
        {\string #\string\hfill}&%
        \string\vtop{\string\parindent=0pt\string\advance\string\hsize%
        by -.5truecm%
        \string#\string\vskip \bibskip
        }\string\cr%
>
\catcode`\{=1\catcode`\}=2
\catcode`\<=12\catcode`\>=12

\def\raf#1#2#3{\ifnum \bz=0 [#1]&#2 \cr\else
\llap{${}_{\rm #3}$}[#1]&#2\cr\fi}

\newread\bibin

\catcode`\{=12\catcode`\}=12
\catcode`\<=1\catcode`\>=2
\def\chiudibib<
\catcode`\{=12\catcode`\}=12
\catcode`\<=1\catcode`\>=2
\immediate\write\bib<}>
\catcode`\{=1\catcode`\}=2
\catcode`\<=12\catcode`\>=12
>
\catcode`\{=1\catcode`\}=2
\catcode`\<=12\catcode`\>=12

\def\makebiblio{
\ifnum\tipobib=0
\advance \maxit by 10pt
\else
\maxit=1.truecm
\fi
\chiudibib
\immediate \closeout\bib
\openin\bibin=\jobname.bib
\ifeof\bibin\relax\else
\raggedbottom
\input \jobname.bib
\fi
}

\openin13=#1.aux \ifeof13 \relax \else
\input #1.aux \closein13\fi
\openin14=\jobname.aux \ifeof14 \relax \else
\input \jobname.aux \closein14 \fi
\immediate\openout15=\jobname.aux

\def\biblio{\*\*\centerline{\titolo References}\*\nobreak\makebiblio}


\ifnum\mgnf=0
   \magnification=\magstep0
   \hsize=15truecm\vsize=20.4truecm\voffset2.truecm\hoffset.5truecm
   \parindent=0.3cm\baselineskip=0.45cm\fi
\ifnum\mgnf=1
   \magnification=\magstep1\hoffset=0.truecm
   \hsize=15truecm\vsize=20.4truecm
   \baselineskip=18truept plus0.1pt minus0.1pt \parindent=0.9truecm
   \lineskip=0.5truecm\lineskiplimit=0.1pt      \parskip=0.1pt plus1pt\fi


\mgnf=0   
\openin14=\jobname.aux \ifeof14 \relax \else
\input \jobname.aux \closein14 \fi
\openout15=\jobname.aux

\footline={\rlap{\hbox{\copy200}}\tenrm\hss \number\pageno\hss}
\def\fiat{}

%
\fiat


\overfullrule=0pt
\headline{{\tenrm \number\pageno:\ {\it Degenerate elliptic resonances}\hss
}}


\centerline{\titolone
Degenerate elliptic resonances}

\vskip1.truecm \centerline{{\titolo
Guido Gentile}$^{\dagger}$ {\titolo and Giovanni Gallavotti}${}^{*}$}
\vskip.2truecm
\centerline{{}$^{\dagger}$ Dipartimento di Matematica,
Universit\`a di Roma Tre, Roma, I-00146}
\vskip.1truecm
\centerline{{}$^{*}$ Dipartimento di Fisica,
Universit\`a di Roma ``La Sapienza'', I-00185}
\vskip1.truecm

\line{\vtop{
\line{\hskip1truecm\vbox{\advance \hsize by -2.1 truecm
\0{\cs Abstract.}
{\it Quasi-periodic motions on invariant tori of an integrable system
of dimension smaller than half the phase space dimension may continue
to exists after small perturbations. The parametric equations of the
invariant tori can often be computed as formal power series in the
perturbation parameter and can be given a meaning via resummations.
Here we prove that, for a class of elliptic tori, a resummation
algorithm can be devised and proved to be convergent, thus extending
to such lower-dimensional invariant tori the methods employed to prove
convergence of the Lindstedt series either for the maximal (\ie KAM)
tori or for the hyperbolic lower-dimensional invariant tori.}} \hfill}
}}

\*\*\*\*
\section(1,Introduction)

\0Quasi-integrable analytic Hamiltonian systems are described by
Hamiltonians of the form $\HHH=\HHH_{0}(\III)+\e\HHH_{1}(\ff,\III)$,
where $(\ff,\III)\in \TTT^{d}\times \AAA$, with $\AAA$ an open
domain in $\RRR^{d}$, are conjugate coordinates
(called angle-action variables), the functions $\HHH_{0}$ and $\HHH_{1}$
are analytic in their arguments, and $\e$ is a small real parameter.
We shall consider for simplicity only Hamiltonians of the form
$$ \HHH = {1\over 2} \III\cdot\III + \e f(\ff) ,
\Eq(1.1) $$
where $\cdot$ denotes the inner product in $\RRR^{d}$.

Kolmogorov's theorem (KAM theorem) yields, for $\e$ small enough, the
existence of many invariant tori for Hamiltonian systems of the form
\equ(1.1): such tori can be parameterized by the corresponding
rotation vectors, at least if the latter satisfy some Diophantine
conditions.  On the other hand Poincar\'e's theorem states the
existence of periodic orbits, which can be parameterized by rotation
vectors satisfying $d-1$ resonance conditions (so that after
a simple linear canonical map one can assume that the
rotation vector is $(\o_{1},0,0,\ldots,0)$).

A natural question is what happens of the invariant tori
corresponding, in absence of perturbations, to rotation vectors
satisfying $s$ resonance conditions, with $1 \le s \le d-2$.  If we
fix the rotation vector as $(\oo,\V0)\=(\o_1,\ldots,\o_r,0,\ldots,0)$
and parameterize the invariant torus for $\e=0$ with the action value
$\III=\V0$ then, after translating the origin in $\RRR^d$ by
$(\oo,\V0)$ and setting $\III=(\AA,\BB)\in
\RRR^r\times\RRR^s,\ff=(\aa,\bb)\in\TTT^r\times\TTT^s$, the
Hamiltonian \equ(1.1) becomes
$$ \HHH = \oo\cdot\AA + {1\over 2} \AA \cdot \AA +
{1\over 2} \BB\cdot\BB + \e f(\aa,\bb) ,
\Eq(1.2) $$
where $(\aa,\AA)\in \TTT^{r}\times\RRR^{r}$ and $(\bb,\BB)\in\TTT^{s}
\times\RRR^{s}$ are conjugate variables, with $r+s=d$, and $\cdot$
denotes the inner product both in $\RRR^{r}$ and in $\RRR^{s}$.
Here we impose that $\oo$ is a vector in $\RRR^{r}$
satisfying
$$ \left| \oo\cdot\nn \right| \ge C_{0} |\nn|^{-\t_{0}} \qquad
\forall \,\nn\in \ZZZ^{r}\setminus \{\V0\} ,
\Eq(1.3) $$
with $C_{0}>0$ and $\t_0\ge r-1$, which is called the {\it Diophantine
condition}; we shall define by $D_{\t_0}(C_{0})$
the set of rotation vectors in $\RRR^{r}$ satisfying \equ(1.3).
We also write
$$ f(\aa,\bb) = \sum_{\nn\in\zzz^{r}} e^{i\nn\cdot\aa} f_{\nn}(\bb) .
\Eq(1.4) $$
We shall suppose that $f$ is analytic in a strip
around the real axis of the variables $\aa,\bb$,
so that there exist constants $F_{0},F_{1},\k_{0}$ such that
$|\dpr_{\bb}^{q}f_{\nn}(\bb)|\le q!F_{0}F_{1}^{q}e^{-\k_{0}|\nn|}$
for all $\nn\in\ZZZ^{r}$ and all $\bb\in\TTT^{s}$.

There are quite a few results on the above problem,
essentially solved, under the assumptions of Theorem 1 below,
in Ref. \cita{JLZ}, and on closely related problems.
We summarize our understanding of the existing results in Appendix \secc(A1).

The equations of motion for the system \equ(1.2),
written in terms of the angle variables alone, are
$$ \ddot \aa = - \e\dpr_{\aa} f(\aa,\bb) , \qquad
\ddot \bb = - \e\dpr_{\bb} f(\aa,\bb) ,
\Eq(1.5) $$
so that, once a solution of \equ(1.5) is found, the action variables
are immediately obtained by a simple differentiation:
$\AA=\dot\aa-\oo$, $\BB=\dot\bb$.

We look for solutions of \equ(1.5), for $\e\neq 0$,
conjugated to the free solution $(\aa_{0}+\o t,\bb_{0},0,0)$, i.e.
we look for solutions of the form
$$ \aa(t) = \pps + \aaa (\pps,\bb_{0}; \e) , \qquad
\bb(t) = \bb_{0} + \bbb(\pps,\bb_{0}; \e) ,
\Eq(1.6) $$
for some functions $\aaa$ and $\bbb$, real analytic
and $2\p$-periodic in $\pps\in\TTT^{r}$, such that the motion
in the variable $\pps$ is governed by the equation $\dot\pps=\oo$.
We shall prove the following result.

\*

\0{\bf Theorem 1.}
{\it Consider the Hamiltonian \equ(1.2), with $\oo\in D_{\t_{0}}(C_{0})$
and $f$ analytic and periodic in both variables.
Suppose $\bb_{0}$ to be such that
$$
\dpr_{\bb}f_{\V0}(\bb_{0}) = \V0 ,\Eq(1.7) $$
and assume that the eigenvalues $a_{1},\ldots,a_{s}$
of the matrix $\dpr_{\bb}^{2} f_{\V0}(\bb_{0})$ are pairwise distinct
and positive, \ie for some constant $a>0$ one has
$a_{i}, \, a_j-a_i \, > \, a\,>\,0$ for all $j>i=1,\ldots,s$.
\\
Then there exist a constant $\lis\e>0$ and a set
$\EE\subset (0,\lis\e)$ such that the following holds.
\\
(i) For all $\e\in\EE$ there are solutions of \equ(1.5)
of the form \equ(1.6), where the two functions
$\aaa(\pps,\bb_{0};\e)$ and $\bbb(\pps,\bb_{0};\e)$
are real analytic and $2\p$-periodic in the variables $\pps\in\TTT^{r}$.
\\
(ii) The relative Lebesgue measure of $\EE\cap(0,\e)$ with respect to
$(0,\e)$ tends to $1$ for $\e\to 0$.
\\
(iii) The functions $\aaa,\bbb$ can be extended to Lipschitz
functions of $\e,\pps$ in $[0,\lis\e]\times\TTT^{r}$.}

\*

\0{\it Remarks.}
(1) From the literature one might expect that the
non-resonance condition on the eigenvalues of $\dpr_{\bb}f_{\V0}(\bb_{0})$
could be avoided; see Appendix \secc(A1).
\\
(2) The case of negative $\e$ was dealt with in Ref. \cita{GG}, with
techniques close to the ones introduced here, and it corresponds to
the case of hyperbolic tori.
\\
(3) The case of mixed stationarity, \ie $\det
\dpr^2_{\bb}f_{\V0}(\bb_{0}) \ne0$ and eigenvalues of
$\dpr^2_{\bb}f_{\V0}(\bb_{0})$ of mixed signs (with non-degeneracy
of the positive ones), can be treated in
exactly the same way discussed in this paper and the above result
extends to this case; cf. Theorem 2 in Section \secc(7).
\\
(4) For $\e\not\in\EE$ the smooth extension in (iii) does not
represent parametric equations of invariant tori: it just says that
their values in the physically interesting set $\EE$ (which turns out
to have dense complement in $[0,\lis\e]$) can be smoothly
interpolated in $\e$. Such (non-unique) extensions are commonly used
for interpolation purposes and are called Whitney extensions.

\*

{\it The novelty and the purpose of our work is the development of a
method of proof based on the existence of a formal power series expression
for the functions $(\aaa,\bbb)$ and its multiscale analysis producing
a rearrangement of its terms, involving summing
many divergent series, which turns it into an absolutely
convergent series.}

\*

The paper is organized as follows.
In Section \secc(2) we recall the basic formalism, following
Ref. \cita{GG}, and in Section \secc(3) we give a simple example of
resummation.

In Section \secc(4) we set up terminology and discuss {\it
heuristically} the ideas governing our resummations, by explaining why
they have to be performed by a multiscale analysis of the series
(which we call ``Lindstedt series'') representing a formal
expansion of the quasi-periodic motions in powers of $\e$.

The singularities are first ``probed'' down to a scale in which
possible resonances between the proper frequencies, \ie the components
of $\oo$, and the normal frequencies, \ie the square roots of
the eigenvalues of $\e\dpr^{2}_{\bb} f_{\V0}(\bb_0)$, are still irrelevant.
The analysis of such singularities leads to what we call {\it
non-resonant} or {\it high frequency resummations}, which
can be treated by the method of Ref. \cita{GG}, \ie of
the hyperbolic case, in which no resonances at all were possible
between proper frequencies and normal frequencies (simply because,
for the Hamiltonian \equ(1.2) the latter did not exist).
Further probing of the singularities leads to what we call
the {\it resonant} (or {\it infrared}) {\it resummations}: the
analysis is more elaborated and it requires new ideas, obtained by
combining the ideas in Ref. \cita{GG} with the ones introduced in
Ref. \cita{Ge}.

In Section \secc(5) we discuss the non-resonant resummations while the
new infrared resummations are studied in Section \secc(6) where a
``fully renormalized series'' is obtained, \ie a resummation of the
series defining the formal expansion of the quasi-periodic solution
\equ(1.6) of the equations of motion \equ(1.5),
{\it which we prove to be absolutely convergent}.

The paper is a self-contained discussion of the construction and of
the convergence of the resummed series. This includes a self-contained
{\it description} of the well-known formal series,
\cita{JLZ}, \cita{GG}. Once this is achieved one has to check that the
defined functions do actually represent parametric equations of
invariant tori: for this we follow, in Appendix \secc(A5),
the analysis of Refs. \cita{GG} and \cita{Ge}.

\*\*
\section(2,Tree formalism)

\0We look for a formal power series  expansion (in $\e$) of the
parametric equations $\hhh=(\aaa,\bbb)$ of the invariant torus
close to the torus $\aa=\pps,\bb=\bb_0$
$$
\hhh(\pps;\e) =
\sum_{k=1}^{\io} \e^{k} \hhh^{(k)}(\pps) =
\sum_{\nn\in\zzz^{r}} e^{i\nn\cdot\pps} \hhh_{\nn}(\e) =
\sum_{k=1}^{\io} \e^k \sum_{\nn\in \zzz^{r}}
e^{i\nn\cdot\pps} \, \hhh^{(k)}_{\nn} , \Eq(2.1) $$
where we have not explicitly written the dependence on $\bb_{0}$.
The power series is easy to derive, see for instance Ref. \cita{GG}:
however its convergence turns out to be substantially harder than the
convergence proof of the Lindstedt series for the maximal KAM tori.
The series constructed below for our problem, which we still
call ``Lindstedt series'', is naturally described in terms of trees.
The coefficients $\hhh^{(k)}_{\nn}$ can be computed as sums of
``values'' that we attribute to trees whose nodes and lines carry a
few labels, which we call ``decorated trees''.

The formalism to define trees, decorations and values
has been described many times and used in the proof of several
stability results in Hamiltonian mechanics. Usage of graphical tools
based on trees in the context of KAM theory has been advocated
recently in the literature as an interpretation of Ref. \cita{E1};
see for instance Refs. \cita{Ga}, \cita{GG}, \cita{BGGM} and \cita{BaG}.

\*
\eqfig{200pt}{141.6pt}{
\ins{-29pt}{75pt}{\rm root}
\ins{14pt}{71pt}{$\nn\kern-2pt=\kern-2pt\nn_{\ell_{0}}$}
\ins{24pt}{91pt}{$\ell_{0}$}
\ins{50pt}{71pt}{$\vvv_0$}
\ins{46pt}{91pt}{$\nn_{\vvvv_0}$}
\ins{86pt}{105pt}{$\h$}
\ins{127pt}{100.pt}{$\vvv_1$}
\ins{191.7pt}{133.3pt}{$\vvv_5$}
\ins{156.7pt}{138pt}{$\h'$}
\ins{121.pt}{125.pt}{$\nn_{\vvvv_1}$}
\ins{92.pt}{41.6pt}{$\vvv_2$}
\ins{158.pt}{83.pt}{$\vvv_3$}
\ins{191.7pt}{100.pt}{$\vvv_6$}
\ins{191.7pt}{71.pt}{$\vvv_7$}
\ins{191.7pt}{-8.3pt}{$\vvv_{11}$}
\ins{191.7pt}{16.6pt}{$\vvv_{10}$}
\ins{166.6pt}{42.pt}{$\vvv_4$}
\ins{191.7pt}{54.2pt}{$\vvv_8$}
\ins{191.7pt}{37.5pt}{$\vvv_9$}
\ins{32pt}{139pt}{$\h$}
\ins{53pt}{130pt}{$\g$}
\ins{12pt}{134pt}{$\g'$}
\ins{63pt}{115pt}{$\vvv$}
\ins{0pt}{117pt}{${\vvv}'$}
}{Fig1}{}
\*\*
\line{\vtop{\line{\hskip1.5truecm\vbox{\advance\hsize by -4.1 truecm
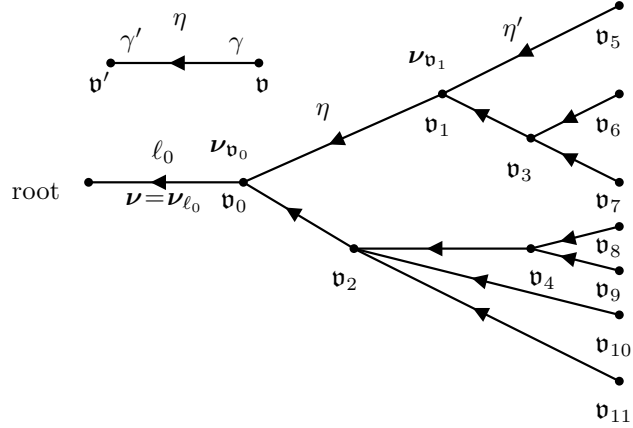
\0{{\css Figure 1.} {\nota A tree $\th$ with $12$ nodes;
one has $p_{\vvvvv_{0}}=2,p_{\vvvvv_{1}}=2,p_{\vvvvv_{2}}=3,
p_{\vvvvv_{3}}=2, p_{\vvvvv_{4}}=2$.  The length of the lines should
be the same but it is drawn of arbitrary size. The endnodes
$\vvvv_{i}$, $i=5,\ldots,11$ can be either nodes or leaves of the
tree. The separated line illustrates the way to think of the label
$\h=(\g',\g)$.\vfil}} } \hfill} }}
\*
\*

\kern-4mm
A tree $\th$ (see Fig. 1) is defined as a partially ordered set of points,
connected by oriented {\it lines}. The lines are
consistently oriented toward the {\it root}, which is the leftmost
point $\vvr$; the line entering the root is called the {\it root line}.
If a line $\ell$ connects two points $\vvv_{1},\vvv_{2}$ and is
oriented from $\vvv_2$ to $\vvv_1$, we say that $\vvv_{2} \prec \vvv_{1}$
and we shall write $\vvv_{2}'\defi\vvv_{1}$ and $\ell_{\vvvv_{2}}\defi\ell$;
we shall say also that $\ell$ exits $\vvv_{2}$ and enters $\vvv_{1}$.
More generally we write $\vvv_{2} \prec \vvv_{1}$ if $\vvv_{2}$
is on the path of lines connecting $\vvv_{1}$ to the root. The points
different from the root will be called the {\it nodes} of the tree.

Each line from $\vvv$ to $\vvv'$ carries a pair $\h$ of labels
$\h=(\g,\g')$ ranging in $\{1,\ldots,d\}$ (marked in Fig. 1 only on
some of the lines for clarity of the drawing). The labels $\g$ and $\g'$
should be regarded as associated with $\vvv$ and $\vvv'$, respectively;
hence with each node $\vvv$ with $p_{\vvvv}$ entering lines
$\ell_{1},\ldots,\ell_{p_{\vvvvv}}$ one can associate $p_{\vvvv}+1$
labels $\g_{0},\g_{1},\ldots,\g_{p_{\vvvvv}}$, with
$\g_{0}=\g_{\ell_{\vvvvv}}$ and $\g_{j}=\g_{\ell_{i}}'$.  Also the
root line (from $\vvv_0$ to the root) carries two such labels and the
one associated with the final extreme of the root line will be called
the {\it root label}.

Fixed any $\ell_{\vvvv} \in \th$, we shall say that the subset of $\th$
containing $\ell_{\vvvv}$ as well as all nodes $\www\preceq \vvv$
and all lines connecting them is a {\it subtree} of $\th$
with root $\vvv'$; of course a subtree is a tree.

Given a tree, with each node $\vvv$ we associate a {\it harmonic} or
{\it mode}, as called in Ref. \cita{GG}, \ie a label
$\nn_{\vvvv}\in\ZZZ^{r}$.  We shall denote by $V(\th)$ the set of
nodes and by $\L(\th)$ the set of lines.  The number $k=|V(\th)|$ of
nodes in the tree $\th$, equal to the number $|\L(\th)|$ of lines,
will be called the {\it order} of $\th$.

We call a node with one entering line and $\V0$ harmonic
label a {\it trivial node}.

With any line $\ell=\ell_{\vvvv}$ we associate (besides the
above mentioned pair $\h_{\ell}=(\g_{\ell},\g_{\ell}')$ of labels
assuming values in $\{1,\ldots,d\}$) a {\it momentum label}
$\nn_{\ell}\in\ZZZ^{r}$ defined as
$$ \nn_{\ell} \= \nn_{\ell_{\vvvv}} = \sum_{\wwww\in V(\th) \atop
\wwww \preceq \vvvv} \nn_{\wwww} .
\Eq(2.2) $$
{\it We shall assume that no tree contains trivial nodes with
the entering line with $\V0$ momentum}: this is an important
restriction, as we shall see.  We call {\it degree} $P(\th)$ of a tree
the order of the tree minus the number of $\V0$ momentum lines, so
that $|V(\th)|-P(\th)$ is their number.

We call $\Th_{\nn,k,\g}$ the set of trees $\th$ of order $k$, \ie with
$|V(\th)|=k$ nodes, and $\Th^o_{\nn,k,\g}$ the set of trees
of degree $k$, \ie with $P(\th)=k$. One has
$\Th_{\nn,k,\g}\ne \Th^o_{\nn,k,\g}$.

Each tree $\th$ ``decorated'' by labels in the described way will have a
{\it value} which is defined in terms of a product of several factors.
\\
$\bullet$ With each node $\vvv$ we associate a {\it node factor}
$$ F_{\vvvv} = \prod_j \dpr_{\g_j} f_{\nn_{\vvvv}}(\bb_0),\Eq(2.3)$$
where the labels $\g_j$ are the $p_{\vvvv}+1$ labels associated with the
extreme $\vvv$ of the $p_{\vvvv}$ lines entering the node $\vvv$ and
of the line exiting it, and the derivatives $\dpr_\g$, with
$\g=1,2,\ldots,r$, have to be interpreted as factors
$(i\nn_{\vvvv})_\g$. Hence $F_{\vvvv}$ is a tensor of rank $p_{\vvvv}+1$.
\\
$\bullet$ With each line $\ell$ carrying labels $\h_\ell=(\g_{\ell},
\g_{\ell'})$ and momentum $\nn_\ell$ we associate a matrix,
called {\it propagator},
$$ \eqalign{
&G_{\ell} \= \d_{\g_{\ell},\g_{\ell}'}
{ 1 \over (\oo\cdot\nn_{\ell})^2 } , \kern3.1cm{\rm if}\ \nn_\ell\ne\V0,\cr
&G_{\ell} \= -\e^{-1} \,(\dpr^2_{\bb}
f_{\V0}(\bb_0))^{-1}_{\g_{\ell},\g_{\ell}'}\,
\chi({\g_{\ell},\g_{\ell}'}>r), \qquad {\rm if}\ \nn_\ell=\V0,\cr}
\Eq(2.4) $$
where $\chi({\g_{\ell},\g_{\ell}'}>r)$ is $1$ if both
$\g_{\ell}$ and $\g_{\ell}'$ are strictly greater than $r$,
and $0$ otherwise.

Given the definitions \equ(2.3) and \equ(2.4) define a {\it value function}
$\Val$, which with each tree $\th$ of order $k$ associates a {\it tree value}
$$ \Val(\th) =\fra{\e^k}{k!}
\Big( \prod_{\vvvv\in V(\th)} F_{v} \Big)
\Big( \prod_{\ell \in \L(\th)} G_{\ell} \Big) ,
\Eq(2.5) $$
where, by the definitions, all labels $\g_{i}$ associated with the
nodes appear twice because they appear also in the propagators: we
make in \equ(2.5) the {\it summation convention} that repeated $\g$
labels associated with nodes and lines are summed over,
with the exception of the label $\g$ associated with the root
(because we do not consider it a node and the corresponding
label $\g$ appears only once in \equ(2.5)).
Therefore \equ(2.5) is a number labeled by $\g=1,\ldots,d$, \ie
$\Val(\th)$ is a vector.

\*

\0{\it Remarks.}
(1) The trees can be drawn in various ways: we can limit the
arbitrariness by demanding that the length of the segments
representing the lines is $1$ (unlike the drawings in the above
figures) and that the angles between the lines are irrelevant. {\it
The combinatorics being very important}, because it matters in the
check of cancellations essential for the analysis, we adopt the
convention that trees are drawn on a plane, have lines of unit length
{\it and carry an identifier label}, that we call {\it number label}
(not shown in the above figures) which distinguishes the lines from
each other even if we ignore the other labels attached to
them. Furthermore two trees that can be superposed by pivoting the
lines merging into the same node $\vvv$, around $\vvv$ itself, are considered
{\it identical}. This is a convention which is useful for checking
cancellations: however  it is by no means the only possible one. Others
are possible and often very convenient in other respects
\cita{G3}, \cita{GM1}, but in a given work a choice has to be made
once and for all.
\\
(2) A line $\ell$ carrying $\V0$ momentum is somewhat special. We
could visualize the part of the tree preceding such lines by
encircling it into a dotted circle: such a representation has been
used in earlier papers, \eg in Ref. \cita{GG}, calling the subtree
$\th_\ell$ with $\ell$ as root line a {\it leaf}. Here, however, we
shall avoid using a special word for the $\V0$ momentum lines and the
subtrees preceding them.
\\
(3) We can think of the propagators as matrices of the form
$$ G_{\ell} = \Big( \matrix{
G_{\ell,\a\a} & G_{\ell,\a\b} \cr
G_{\ell,\b\a} & G_{\ell,\b\b} \cr} \Big),
\Eq(2.6) $$
where $G_{\ell,\a\a}$, $G_{\ell,\a\b}$, $G_{\ell,\b\a}$ and
$G_{\ell,\b\b}$ are $r\times r$, $r\times s$, $s\times r$ and $s\times
s$ matrices.
\\
(4) The value of a tree $\th$ defined above has no pole at
$\e=0$ if $\Val(\th)\ne0$ because every line
with $\V0$ momentum is preceded by at least two nodes, so that the
total power of $\e$ to which the value is proportional is always
non-negative and, in fact, it is necessarily positive: we need to take
into account that $\dpr_{\bb} f_{\V0}(\bb_0)\=0$ and that our trees
contain no trivial nodes with one entering line with $\V0$ momentum.
Note that ${\rm Val}(\th)$ is a monomial in $\e$ of degree $P(\th)$.
\\
(5) In the case of maximal tori and if $\Val(\th)\ne0$
there are {\it no lines with $\V0$
momentum} for systems described by the Hamiltonians \equ(1.1):
indeed $s=0$, see also \cita{Ga}. In this case the number
of nodes, \ie the tree order, coincides with the power of $\e$
associated with the monomial in $\e$ defined by the tree value, \ie
with the tree degree. In general, however, the order $|V(\th)|$ of a
tree can be larger than its degree $P(\th)$: $|V(\th)|\ge P(\th)\ge
\fra12|V(\th)|$.

\*

The above definitions uniquely attribute a value to each tree.  The
following result states the existence of formal solutions to \equ(1.5)
which are conjugated to the unperturbed motion \equ(1.4), provided the
value $\bb_{0}$ is suitably fixed.  The proof is an algebraic check
which does not distinguish the possible signs of $\e$ and can be taken
from Ref. \cita{GG} where it is done in the case $\e<0$.

\*

\0{\cs Lemma 1.}
{\it The Fourier transform of the power series solution
$\hh=(\aaa,\bbb)$ of \equ(1.5) of the form \equ(2.1) is obtained by
writing (the definition of $\Th^o_{k,\nn,\g}$ follows \equ(2.2))
$$ \e^{k} h^{(k)}_{\nn,\g} =
\sum_{\th\in \Th^o_{k,\nn,\g} } \Val(\th)
\Eq(2.7) $$
for all $\nn\in\ZZZ^{r}$, all $k\in\NNN$ and $\g=1,\ldots,d$.}

\*

The expression \equ(2.7) is well defined at fixed $k$ and the
sum over $k$ gives what we call
the {\it formal power series solution} for the equations
for the parametric representation \equ(2.1), \equ(1.6) of the invariant tori.

\*\*
\section(3, The simplest resummation)

\0The power series in $\e$ in
\equ(2.1) and its Fourier transform defined by the sum
over $k$ of \equ(2.7) may be not convergent {\it as a power series}
(as far as we know). The problem is difficult
because if in \equ(2.7) we replace $\Val(\th)$ with $|\Val(\th)|$  the
series certainly diverges.

Our aim, as stated in the introduction, is to show that nevertheless a
meaning to the series can be given. We shall show that the tree values
can be further decomposed into sums of several other quantities and
that the various contributions to the series can be rearranged by
suitably collecting them into families: the sums of the contributions
from each family leave us with a new series (no longer a power series
in $\e$) which is in fact convergent and its sum solves the problem
of constructing the parametric representations $\hh=(\aaa,\bbb)$,
\equ(2.1), of the invariant tori at least for all $\e\in\EE$, with
$\EE$ a set with $0$ as a density point (\ie a Lebesgue point).

For this purpose we need to define and consider more involved trees
and more involved definitions of their values.  We begin by remarking
that trees may contain {\it trivial nodes}, \ie nodes with $\V0$
harmonic separating two lines with equal momentum $\nn\ne\V0$.

One can suppose that {\it no tree contains trivial nodes} provided we
use for all lines, with momentum $\nn\ne\V0$ and labels $\g,\g'$
associated with the extremes, the {\it new} propagators
$$ \lis g(x;\e)\defi(x^2-M_0)^{-1},\qquad x\defi \oo\cdot\nn,\qquad
M_0\defi\e\left(
\matrix{ 0 & 0 \cr 0 &
\dpr_{\bb}^{2}f_{\V0}(\bb_{0}) \cr} \right) .
\Eq(3.1)
$$
This is a {\it resummation of many divergent series} obtained by
adding the values of trees obtained from a tree without trivial nodes
by ``insertion'' of an arbitrary number of trivial nodes on the
branches with momentum $\nn\ne\V0$: this requires summing series, one
per branch of a tree {\it without trivial nodes}, which are geometric
series with ratio given by the $d\times d$ matrix
$z=\fra{M_{0}}{(\oo\cdot\nn)^{2}}$; $|z|$ can be larger than $1$
because the $s$ non-zero eigenvalues $\e a_j$, $j=1,\ldots,s$, of
$M_0$ are unrelated to $x=\oo\cdot\nn$.\footnote{${}^1$}{
\nota Note that since the tree lines are numbered (\ie they are regarded as
distinct) adding $p$ nodes on a line $\ell$ changes the combinatorial
factor $k!^{-1}$ in \equ(2.5) into $(k+p)!^{-1}$: however the new $p$
lines thus produced can be chosen in ${k+p} \choose p$ ways and ordered
in $p!$ ways so that we can ignore the extra number labels on $\ell$ and
use as combinatorial factor $(k+p)!^{-1}{k+p\choose
p}p!=k!^{-1}$.\vfil}

Therefore replacing $\sum_{p=0}^\io z^p$ by $(1-z)^{-1}$ {\it is not
rigorous and needs to be eventually justified}.  Certainly we must at
least suppose that $x^2-M_0$ can be inverted: otherwise the values of
the trees representing the new series might even be meaningless! (\ie if
some lines will have momentum $\nn$ such that $\det (x^2-M_0)=0$).
This happens for a dense set of $\e$'s and we have to exclude such
$\e$'s by imposing conditions on the eigenvalues
$\l^{[0]}_{r+j}\=\e a_j,\, j=1,\ldots,s$, \ie on $\e$.

For uniformity of notations it is convenient to assume that $\e$ is in
an interval $(\e_{\rm min},4\e_{\rm min}]$ related to
the largest eigenvalue $\l^{[0]}_d\=a_s\e$ of $M_0$ by

$$\l^{[0]}_d\=\e \,a_s\,\in\, I_C\,\defi
\Big(\fra14 C^2,{C}^2 \Big], \qquad
C\,\defi\, C_{0} 2^{-n_0},\qquad n_0 \ge 0 ,
\Eq(3.2)$$
where $C_{0}$ is the Diophantine constant in \equ(1.3) (fixed
throughout the analysis); thus $I_C$ is an interval of size $O(C^2)$
(\ie $\fra3{4}C^2$).  In other words we find it convenient to measure
$\e$ in units of $C_{0}^2 a_s^{-1}$ via an integer $n_0$. {\it We \ap assume,
for simplicity, the restrictions $a_{s}\e\le C_{0}^{2}$ and $\e\le1$}.

To give a meaning to $(x^2-M_0)^{-1}$ it would suffice to require
$|x^2-\e a_j|\ne0$ for all $j$ thereby excluding ``only'' a
denumerable ({\it dense}) set of values of $\e$, of $0$ length; however
stronger conditions will be needed in order to analyze the convergence
problems and we begin by imposing them in a form which will be useful
later.  Setting for later use $\ul\l_j^{[0]}(\e)\defi
\l^{[0]}_j$, the conditions that we impose on $\l^{[0]}_{j}$, \ie on
$\e$, are that for all $x=\oo\cdot\nn\neq 0$ and for all independent
choices of the signs $+$ or $-$
$$ \G(x)\defi \min_{j\ge i} \Big\{
\Big| |x|-\Sqrt{\ul\l^{[0]}_{j}(\e)}\Big| \,,\,
 \Big| x \pm
\Sqrt{\ul\l^{[0]}_{j}(\e)} \pm \Sqrt{\ul\l^{[0]}_{i}(\e)} \Big|
\Big\} \ge 2^{-(\lis n_0-1)/2} \fra{C_{0}}{|\nn|^{\t_{1}}}
\Eq(3.3)$$
for $\t_{1}$ suitably large and $\lis n_0$ suitably larger than $n_0$,
see \equ(4.2). This excludes a closed set of values of $\e$
in the considered interval $I_C$, \equ(3.2): its measure
can be estimated without difficulties.
Let
$$ \t_{1} = \t_{0} + r + 1 ,
\Eq(3.4) $$
be a convenient, although somewhat arbitrary, choice; then
the total measure of the excluded set is
$$ \le\,2^{-(\lis n_0-1)/2} C^2 \, K,\Eq(3.5)$$
where $K$ is a suitable constant; see Appendix \secc(A2).
Hence the measure of the complement of the set $\EE_{\lis n_0-1}$ where
\equ(3.3) is verified is a {\it small fraction} of order $C^{1/2}$
of the measure of the interval $I_C$, whose size is proportional to $C^2$,
in which we let $\e$ vary, at least if $\lis n_0$ is large.

\*\*
\section(4,Resummations: semantic and heuristic considerations)

\0Replacing the propagators $x^{-2}$ of the lines by
$(x^2-M_0)^{-1}$ we obtain a representation of the parametric
equations $\hh$ involving simpler trees (\ie trees with
no trivial nodes).  The new representation is
a series in which each term is well defined if $\e$ is in the large
set $\EE_{\lis n_0-1}\subset I_C$ in which
\equ(3.3) holds. This is quite different from the original Lindstedt
series in \equ(2.7) whose terms are well defined for all $\e$.

We should also stress that the resummed series is in a sense more
natural: the $\V0$ momentum lines now appear as less anomalous
because their propagator is much more closely related to
$(x^2-M_0)^{-1}$. One can say that it is just the latter evaluated at
$x=0$ with the meaningless entries (\ie the first $r$ diagonal
entries) replaced by $0$. Another way of saying the latter property is
that lines $\ell$ with $\V0$ momentum and labels $\g_{\ell},
\g_{\ell}'\le r$ are forbidden. One should not be surprised by this fact:
it is the generalization of the corresponding property in the case
of maximal tori ($r=d$) in which this means that lines with $\V0$ momentum
are forbidden. The latter property goes back to Poincar\'e's
theory of the Lindstedt series and is the key to the proof of
the KAM theorem and of cancellations which make the formal Lindstedt series
for maximal tori absolutely convergent; see Refs. \cita{E1} and \cita{Ga}.
However the new series is still only a formal representation because
it is by no means clear that it is absolutely converges.

The next natural idea is to try to establish convergence by further
modifying the propagators, changing at the same time the trees structure,
until one achieves a formal representation whose convergence
will be ``easy'' to check. Once achieved a formal representation
which is convergent we shall have to check that it really solves the
equations for $\hh$.

The modification of the trees structure will be performed by steps:
at each step, labeled by an integer $n=0,1,\ldots$, the propagators of
the lines with non-zero momentum
will have been modified acquiring labels $[0],[1],\ldots [n-1]$, or
the label $[\ge n]$, indicating that they are given no longer by
$(x^2-M_0)^{-1}$ but by a matrix proportional to
$(x^2-\MM^{[\le p]})^{-1}$, if their label is
$[p]$, with $p<n$, or (with a different proportionality factor) 
to $(x^2-\MM^{[\le n]})^{-1}$,
if their label is $[\ge n]$; here $\MM^{[\le p]}$ are suitable
matrices. {\it Here and in the following the symbols $[\le n]$
and $[\ge n]$ are consistently used. Hence $[\ge n]$ {\it does not}
denote the set of scales $[p]$ with $p\ge n$, and in fact it is
just a different scale; likewise $[\le n]$ does not ``include'' $[p]$
even if $p\le n$. In other words one has to regard the
symbols $[\le n]$, $[n]$ and $[\ge n]$ as unrelated symbols.
This might appear unusual but it turns
out to be a good notation for our purposes}.

The proportionality factor depends on $x$ and contains cut-off functions
which vanish unless $x^2-\MM^{[\le p]}$ has smallest eigenvalue
of order $O(2^{-2p}C_0^{2})$; the cut-offs are so devised that if the
propagator does not vanish its denominator has a minimum size
proportional to $2^{-2p}$ and the ratio between its minimum and
maximum values will be bounded above and below by a {\it
$p$-independent constant}. No modification will be made of the
propagators of the $\V0$ momentum lines: for uniformity of notation
we shall attach a label $[\io]$ to such lines.

Considering trees with no trivial nodes in which each line carries also
an extra {\it scale} label $[\io],[0],[1],\ldots [n-1],[\ge n]$ a new formal
representation of $\hh$ will be obtained by assigning, to the trees,
values defined by the same formula in \equ(2.5),
with the propagators $G_\ell$ replaced by the new
propagators, that we denote $g^{[p]}_{\ell}$
if the line $\ell$ carries the label $[p]$, with $p=\io,0,\ldots,n-1$,
and $g^{[\ge n]}_{\ell}$ if the line carries the label $[\ge n]$.
When the line $\ell$ is on scale $[p]$, with $p=0,\ldots,n-1$,
or $[\ge n]$ or $[\infty]$, 
then the corresponding propagator will be proportional to
$(x^2-\MM^{[\le p]})^{-1}$ or $(x^2-\MM^{[\le n]})^{-1}$ or 
$(\e \dpr^2_{\bb} f_\V0(\bb_0))^{-1}$.

The construction will be performed in such a way that
the matrices $(x^2-\MM^{[\le p]})$  will be defined by series which
{\it will be proved to be convergent}; furthermore if we only considered
the contributions to the formal representation of $\hh$ coming from
trees {\it in which no propagator carries the ``last label'' $[\ge n]$
then the corresponding series would be convergent}.

We express the latter property by saying that {\it the performed
resummations regularize the formal representation of $\hh$ down to
scale $[n-1]$, or that the propagators singularities are probed down
to scale $[n-1]$}. The problem of course remains to understand the
contributions from the trees containing lines with label $[\ge n]$.
The construction will be such that their propagators are also properly
defined because the matrices $\MM^{[\le n]}$ will always be well
defined by convergent series (as we shall see).
However for the lines whose label is $[\ge n]$ no useful positive
lower bound, not even $n$-dependent, can be given on the smallest
eigenvalue of the denominators in the corresponding propagators.

We shall say that the lines with scale $[\ge n]$ probe the singularity
all the way down to the smallest frequencies or all the way down in
the {\it infrared} scales. Thus in spite of the convergence of the
contributions to $\hh$ coming from trees with labels
$[\io],[0],[1],\ldots,[n-1]$ the representation of $\hh$ remains formal.

Therefore we shall proceed by increasing the value of $n$ trying to
take the limit $n\to\io$. This is the procedure followed in the case
of the theory of hyperbolic tori in Ref. \cita{GG}. In that case,
however, the propagators denominators $(x^2-\MM^{[\le n]})$ had
eigenvalues {\it always bounded below proportionally to $x^2$}. Indeed
the last $s$ eigenvalues of $\MM^{[\le n]}$ were negative whereas
the first $r$ remained close to zero within $O(\e x^2)$ 
(a non-trivial property, however, due to remarkable cancellations well
known in the KAM theory, \cita{Ga}).

Here the matrices $x^2-\MM^{[\le n]}$ will be shown to have
the first $r$ eigenvalues differing by a factor $(1+O(\e^2))$
and the last $s$ differing by $O(\e^{2})$ with respect to those of
$x^2-M_0$ (which has by construction $r$ eigenvalues $x^2$ and $s$
eigenvalues $x^2-\e a_j$ $j=1,\ldots,s$).  Thus the denominators can become
small because {\it either} $x^2$ gets close to $0$ {\it or} because it
gets close to $\e a_1,\ldots \e a_s$.  Therefore the regularization
will have to be split in two parts. The first part will concern
regularizing the scales $[p]$ with $p$ such that the eigenvalues of
$x^2-\MM^{[\le p]}$ remains bounded below proportionally to $x^2$;
we shall call this part of the analysis the {\it high frequencies
resummation}. The other part, which we shall call the {\it infrared
resummation}, will concern the regularization of the scales $[p]$, in
which $x$ can be so close to some $\e a_j$ that the denominators
cannot be bounded below proportionally to $x^2$.

We associate with each momentum $\nn$ the frequency $x=\oo\cdot\nn$
and we measure the strength of this resonance by the integer $p$
if $D(x;\e)\simeq C_{0}^2 2^{-2p}$, with
$$ D(x;\e)=\min_j\Big|x^2-\ul\l^{[0]}_{j}(\e) \Big|
\defi \Big|x^2-\ul\l^{[0]}_{j_{\e}(x)}(\e) \Big| .
\Eq(4.1)$$
Therefore the condition that the resonance strength of the frequency
$x$ be bounded
below proportionally to $|x|$ is that $p$ is not too large compared to
$n_0$ defined in \equ(3.2), so that $x^2$ stays away from the corresponding
eigenvalue ${\ul\l^{[0]}_j(\e)}$ by more than a small fraction of
the minimum separation $\d$ between the distinct eigenvalues.
For instance we can require $D(x;\e)\ge 2^{-2(\lis n_{0}+1)}C_{0}^2
\ge \d/4$. This gives $p\le \lis n_0$, with
$$ \lis n_0=n_0+\lis n, \qquad
\lis n \defi - 1 + \fra12\log_2\fra{1}{\r}, \qquad
\r= \fra{1}{4} a_{s}^{-1}\min\{ {a_{1}}, \min_{j} \{ {a_{j+1}}-{a_{j}}\} \} .
\Eq(4.2) $$
In fact the requirement could be fulfilled with $\lis n$ one unit
larger: the interest of using the above value of $\lis n$ will emerge
later (if $s=1$ one interprets $\r=\fra14$).

\*

We then perform the analysis by defining recursively the matrices
$\MM^{[\le p]}(x;\e)$ for $p=0,\ldots,\lis n_0$ with eigenvalues
$\l^{[p]}_j(x,\e)$ verifying for a suitable constant $\g>0$
$$ |\l^{[p]}_j(x,\e)-\l^{[0]}_j(\e)|<\g \e^2,\qquad p\le \lis n_0,
\Eq(4.3)$$
so that if the label $p$ of the line with frequency $x$ is
$p\le \lis n_{0}$ then one has, if $\fra12 a_s 2^{-2(\lis n+1)}-\g\e\ge0$,
$$ |x^2-\l^{[p]}_j(x,\e)|\ge \fra12 D(x,\e)+\fra12 D(x,\e)-
\g \e^2\ge \fra12 D(x,\e)\ge 2^{-2(\lis n+2)} |x|^2
\Eq(4.4)$$
where the last step is obvious if $|x|^2\ge 2 \ul{\l}^{[0]}_{d}(\e)$,
otherwise it follows from the inequality
$$ {D(x;\e)} \ge 2^{-2(\lis n +1)}
2^{-2n_{0}}C_{0}^2 \ge 2^{-2(\lis n+1)} \ul{\l}^{[0]}_{d}(\e)\ge
2^{-2(\lis n+1)-1} |x|^2 .
\Eq(4.5) $$
We can say that for $p\le \lis n_0$ the strength of the
singularity is dominated by the distance $|x|$ to the origin,
\ie by the ``classical'' small divisors $x^{-2}$ provided, of course,
the matrices $x^2-\MM^{[\le p]}(x;\e)$ remain close enough to $x^2 - M_0$
(which we shall check). Furthermore the convergence
of the sum of all values of trees with no line
label $[\ge \lis n_0]$ will be performed exactly
along the lines of Ref. \cita{GG} because the bound \equ(4.4) guarantees
that in evaluating such trees one does not probe the singularities
close to the eigenvalues of $M_0$.

The departure from the method in Ref. \cita{GG} occurs when we consider
trees in which lines bear the label $[\ge \lis n_0]$.
The problem will again be studied by a
multiscale analysis which will have to be suitably modified to allow
probing the new singularities arising from the resonances between the
frequencies $x$ and the $\Sqrt{\ul\l^{[0]}_j(\e)},\, j>r$. The difficulty is that
the propagator $g^{[\ge \lis n_0]}_{\ell}$ {\it will not be singular
exactly at the frequencies $\Sqrt{\ul\l^{[0]}_j(\e)}\ne0$}
but at the frequencies fixed by the roots of the eigenvalues
$\l^{[\le \lis n_0]}_j(x;\e)$ of the matrices $\MM^{[\le \lis n_0]}(x;\e)$.
The latter not only are slightly different from those of $M_0$ but
will turn out to depend also on $x$.

This means that $D(x;\e)$ and even $\D(x;\e)=\min_{j} \big| \,
x^2- \l^{[\le \lis n_0]}_j(x;\e) \big|$ no longer provide a
good estimate of the strength of the singularity,
because $D,\D$ vanish at the ``wrong places''.
In fact we shall have to perform a multiscale analysis to
resolve the infrared singularities, and it will happen that at each of the
new scales with labels $[p]$, with $p\ge \lis n_0$, the singularities
will keep moving. 

Suppose to have regularized the series up to scale
$[n-1]$, with $n > \lis n_0$, introducing suitably matrices
$\MM^{[\le p]}(x;\e)$, with $p=\lis n_0,\ldots, n$, thus pushing the
probe of the singularities down to scales $C_{0} 2^{-n}$; then to
avoid meaningless expressions
we shall have to impose on the eigenvalues of the last propagator,
proportional to $(x^2-\MM^{[\le n]}(x;\e))^{-1}$, a condition
like \equ(3.3). Since the eigenvalues depend on $n$ and $x$ this
risks to imply that we have to discard too many $\e$'s; in the limit
$n\to\io$: when, finally, the singularities will have been probed on all
scales, or even for large enough scales,
we might be left with an empty set of $\e$'s rather than with a
set of almost full measure.

Physically the difficulty shows up because of the possibility of
resonances between the proper frequencies of the quasi-periodic motion
on the tori and the normal frequencies. It will be
studied and solved in Section \secc(6) below, while in Section
\secc(5) we shall discuss the simpler regularization of the series for $\hh$
on the high frequency scales.

\*

The spirit informing the analysis is very close to the techniques used
in harmonic analysis, in quantum field theory and in statistical
mechanics, known as ``renormalization group methods'' (see
Refs. \cita{Fe}, \cita{GM1}, \cita{Ga01} and \cita{Ga02}). The
methods are also based on a ``multiscale decomposition'' of the
propagators singularities. We introduced and adopt the above
terminology because we feel that it is suggestive and provides useful
intuition at least to the readers who have some acquaintance with the
renormalization group approach and multiscale analysis.

\*\*
\section(5,Non-resonant resummations)

\0The resummations  will be defined via trees with no trivial nodes and
with lines bearing further labels. Moreover the definition of propagator
will be changed, hence the values of the trees will be different
from the ones in Section \secc(3): they are constructed recursively.


\eqfig{340pt}{50pt}
{\ins{9pt}{38pt}{$\psi_0$}
\ins {30pt}{-4pt}{$\st C_0^2/4$}
\ins {60pt}{-4pt}{$\st C_0^2$}
\ins {90pt}{-4pt}{$\st D$}
\ins {150pt}{-4pt}{$\st C^2_0/4$}
\ins {180pt}{-4pt}{$\st C^2_0$}
\ins{129pt}{38pt}{$\chi_0$}
\ins {210pt}{-4pt}{$\st D$}
\ins{249pt}{38pt}{$\lis\chi_0$}
\ins{245pt}{-4pt}{$\st C^2_0/4^2$}
\ins {270pt}{-4pt}{$\st C^2_0/4$}
\ins {300pt}{-4pt}{$\st C^2_0$}
\ins {330pt}{-4pt}{$\st D$}
}
{Fig2}{}
\*
\line{\vtop{\line{\hskip2.5truecm\vbox{\advance\hsize by -5.1 truecm
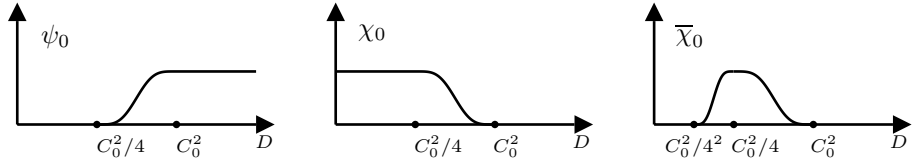
\0{{\css Figure 2.} {\nota The first graph is $\psi_0$,
the second is $\chi_0$ and the third
is $\lis\chi_0\=\psi_1\chi_0$.\vfill}}}\hfill}}}
\*

Instead of the sharp multiscale decomposition considered in
Ref. \cita{GG} here it will be convenient to work with a smooth one as
in Ref. \cita{Ge}.  Let $\psi(D)$ be a $C^{\io}$ non-decreasing compact
support function defined for $D\ge0$, see Fig. 2, such that
$$ \psi(D) =
1 , \quad {\rm for} \quad D \ge C_{0}^2,\qquad \psi(D) =
0 , \quad {\rm for} \quad D \le C_{0}^2/4, \Eq(5.1) $$
where $C_{0}$ is the Diophantine constant in \equ(1.3), and 
let $\chi(D)=1-\ps(D)$.
Define also $\psi_{n}(D)=\psi(2^{2n}D)$ and
$\chi_n(D)=\chi(2^{2n}D)$ for all $n\ge 0$.
Hence $\psi_{0}=\psi,\,\chi_{0}=\chi$ and
$$ 1\=\ps_n(\D(x;\e))+\chi_n(\D(x;\e)),\qquad\hbox{for all}\ n\ge0 ,
\Eq(5.2)$$
for all choices of the function $\D(x;\e)\ge0$: in particular for
$\D(x,\e)=D(x)$ with $D(x)$
defined in \equ(5.3) below. We set the following notations.
\*

\0{\cs Definition 1.} {\it
Let $\lis n_0,\lis n$ be as in \equ(4.2) and $D(x;\e)$ as in \equ(4.1).
\\
(i) Divide the interval $I_C\=[\e_{\rm min},4\e_{\rm min}]$, where $\e$
varies, see \equ(3.2), into a finite number of small intervals $I$
of size smaller than $\fra{1}{2}\e_{\rm min}\r$, see \equ(4.2), \ie
smaller than a fraction of the minimum separation between the
eigenvalues $0,a_1,\ldots,a_s$. Define
$$ D(x;I)=\min_{\e\in I}D(x;\e)=
\min_{\e \in I} \min_j\Big|x^2-\ul\l^{[0]}_{j}(\e)
\Big| = \min_{\e\in I}\Big|x^2-\ul\l^{[0]}_{j(x)}(\e)\Big| .
\Eq(5.3) $$
where $j(x)$ is the smallest value of $j$ for which the last equality
holds: exceptionally there might be $2$ such labels. The $j(x)$
is $\e$-independent, by construction, for $\e\in I$.}

\*

\0{\it Remarks.} (1) Note that, as a consequence of the definition
of the intervals $I$ and of $D(x;I)$ as given by \equ(5.3),
one has, for all $\e\in I $,

$$ \min_{j} \Big| x^2-\ul\l^{[0]}_{j}(\e) \Big|
\ge \fra{1}{2} \Big| x^2-\ul\l^{[0]}_{j(x)}(\e) \Big| ,
\Eq(5.4) $$
\\
(2) If $\e$ is in one of the intervals $I$ and  $x$ verifies
$D(x;I)\le C^2_0 2^{-2 \lis n_{0}}$ then there is only one value
of $j$ for which last equality in \equ(5.3) holds.
\\
(3) {\it We shall fix, from now on, $\e$ in one of the
intervals $I\subseteq I_C$}. Remark that $D(x;I)$ is piecewise
linear in $x^2$ with slope equal to $1$ in absolute value 
for $x$ in the regions where it will be considered (see below)
and we simplify the notation by setting
$$ D(x)\defi D(x;I) .
\Eq(5.5) $$
\\
(4) The number of intervals $I\subset I_{C}$ can and will be taken
independent of $\e_{\rm min}$, \ie of the interval $I_C$ where
$\e$ varies, and equal to a fixed integer $\le 6\r^{-1}$.
\\
(5) From now on we only consider trees with no trivial nodes.

\*

A simple way to represent the value of a tree as sum of many terms
is to make use of the identity in \equ(5.2).
The propagator $ g(x;\e)\= g^{[\ge 0]}(x;\e) \defi
(x^2-M_0)^{-1}$ of each line with non-zero momentum (hence
with $x\neq0$) is written as
$$ g^{[\ge 0]}(x;\e) = \psi_{0}(D(x)) \, g^{[\ge 0]}(x;\e) +
\chi_{0}(D(x))\,g^{[\ge 0]}(x;\e) \defi
g^{[0]}(x;\e) + g^{\big\{\ge 1\big\}}(x;\e) ,
\Eq(5.6)$$
and we note that $g^{[0]}(x;\e)$ vanishes if $D(x)$
is smaller than $(C_{0}/2)^2$, see Fig. 2.

If we replace each $g^{[\ge0]}(x;\e)$ with the sum in \equ(5.6) then
the value of each tree of order $k$ is split as a sum of up
to $2^k$ terms\footnote{${}^2$}{\nota Not necessarily $2^k$ because
there might be
lines on scale $[\infty]$ whose propagator is not decomposed.} 
which can be identified by affixing
on each line with momentum $\nn\ne\V0$ a label $[0]$ or $\big\{\ge 1\big\}$.
Further splittings of the tree values can be achieved as follows.

\*

\0{\cs Definition 2.}
{\it For $p=1,\ldots,\lis n_0$, let
$\MM^{[p]}(x;\e)$ be matrices with eigenvalues
$\l_j^{[p]}(x;\e)$, $p=1,\ldots,n$; we set $\MM^{[0]}(x;\e)\=M_0$ and
$\MM^{[\le n]}(x;\e)=\sum_{p=0}^{n} \MM^{[p]}(x;\e)$.
Define for $0<n\le \lis n_0-1$
$$\eqalign{
& g^{[n]}(x;\e)\defi
\fra{\psi_{n}(D(x))\prod_{m=0}^{n-1}
\chi_m(D(x))}{
x^2-\MM^{[\le n]}(x;\e)}, \cr
& g^{\big\{\ge n\big\}}(x;\e)\defi
\fra{\prod_{m=0}^{n-1}\chi_m(D(x))}
{x^2-\MM^{[\le n-1]}(x;\e)} ,
\cr
& g^{[\ge n]}(x;\e)\defi \fra{\prod_{m=0}^{n-1}
\chi_m(D(x))}{
x^2-\MM^{[\le n]}(x;\e)} , \cr}
\Eq(5.7) $$
and $g^{[0]}(x;\e)=\psi_{0}(D(x))\,(x^2-M_0)^{-1}$. 
We call the labels $[n],\{\ge n\},[\ge n]$ scale labels.}

\*

\0{\it Remarks.}
(1) The products $\prod_{m=0}^{n-1}\chi_m(D(x))$ can be simplified to
involve only the last factor: we keep the notation above as it is a
notation that naturally reflects the construction. The
propagators  $g^{\big\{\ge n\big\}}$ play a subsidiary
role and are here for later reference.
\\
(2) The matrices $\MM^{[p]}(x;\e)$ will be defined recursively under
the requirement that the functions $\hh$ defining the parametric
equations of the invariant torus will be expressed in terms of trees
whose lines carry scale labels indicating that their values are
computed with the propagators in \equ(5.7).
\\
(3) Note that if we defined $\MM^{[\le p]}(x;\e)\=M_0$, \ie
$\MM^{[p]}\=0$ for $p>0$, then (recall that we consider only trees
without trivial nodes) we would naturally decompose
(see below for details) the tree values into
sums of many terms keeping obviously each total sum constant
by repeatedly using \equ(5.2), thus  meeting the requirement
in Remark (2) above. {\it This would be of no interest}.
Therefore we shall try to define the matrices $\MM^{[p]}(x;\e)$
so that the sum of the values of {\it new trees} (with no trivial
nodes and whose nodes and lines $\ell$ still carry harmonic and
momentum labels as well as scale labels $[\io],[0],\ldots,[n-1],[\ge n]$)
remain the same provided their values are evaluated by using the
propagators in \equ(5.7) and {\it we shall try to define
$\MM^{[\le p]}(x;\e)$, so that there is also control of the convergence}.
\\
(4) {\it In other words we try to obtain a graphical representation of
$\hh$, involving values of trees which are easier to study at the
price of needing more involved propagators}.  This is a typical method
employed in KAM theory \cita{GBG}, and in other fields.

\*

To define recursively the matrices we introduce the notions of
clusters and of self-energy clusters of a tree whose lines and
nodes carry the same labels  introduced so far and {\it in addition}
each line carries a scale label which can be either $[\io]$,
if the momentum of the line is zero, or $[p]$,
with $p=0,\ldots,\lis n_0-1$, or $[\ge\lis n_0]$.
Given a tree $\th$ decorated in this way we give
the following definition, for $n\le \lis n_0$.

\*

\0{\cs Definition 3.} (Clusters)
\\
{\it (i) A {\rm cluster} $T$ on scale $[n]$, with $0 \le n$,
is a maximal set of nodes and lines connecting them
with propagators of scales $[p]$, $p \le n$, one of which, at least,
of scale exactly $[n]$. We denote with $V(T)$ and $\L(T)$ the
set of nodes and the set of lines, respectively, contained in $T$.
The number of nodes in $T$ will define the order of $T$, and it
will be denoted with $k_{T}$.
\\
(ii) The $m_{T}\ge 0$ lines entering the cluster $T$ and the possible
line coming out of it (unique if existing at all) are called the
{\rm external lines} of the cluster $T$.
\\
(iii) Given a cluster $T$ on scale $[n]$, we shall
call $n_{T}=n$ its scale.}

\*

\0{\it Remarks.}
(1) For instance if $n=0$ the scale of the lines in the cluster
can only be $[0]$.
\\
(2) Here $n \le \lis n_{0} - 1$. However the definition above is given
in such a way that it will extend unchanged when also
scales larger than $\lis n_{0}$ will be introduced.
\\
(3) The clusters of a tree can be regarded as sets of lines
hierarchically ordered by inclusion and have hierarchically
ordered scales.
\\
(4) A cluster $T$ is not a tree (in our sense);
however we can uniquely associate a
tree with it by adding the entering and the exiting lines and by
imagining that the lower extreme of the exiting line is the root and
that the highest extremes of the entering lines are nodes carrying a
harmonic label equal to the momentum flowing into them.

\*

\0{\cs Definition 4.} (Self-energy clusters)
\\
{\it (i) We call {\rm self-energy cluster} of a tree $\th$ any cluster $T$
of scale $[n]$ such that $T$ has only one entering line
$\ell_{T}^{2}$ and one exiting line $\ell_{T}^{1}$,
and furthermore $ \sum_{\vvvv\in V(T)} \nn_{\vvvv} = \V0$.
\\
(ii) The degree of a self-energy cluster is the number of nodes.}

\*

\0{\it Remark.}
The essential property of a self-energy cluster is that it has
necessarily just one entering line and one exiting line, and both have
{\it equal momentum} (because $\sum_{\vvvv\in V(T)}\nn_{\vvvv}=\V0$).
Note that both scales of the external lines of a self-energy cluster
$T$ are strictly larger than the scale of $T$ as a cluster,
but they can be different from each other by at most one unit.
Furthermore the degree of a self-energy cluster is $\ge2$. 
Of course no self-energy cluster can contain any line on scale $[\io]$.

\*
\0{\cs Definition 5.}  (Self-energy matrices)
\\
{\it (i) Let $\Th^{\RR}_{k,\nn,\g}$ be the set of trees of order $k$
with root line momentum $\nn$ and root label $\g$ which contain
neither self-energy clusters nor trivial nodes.  Such trees will be
called {\rm renormalized trees}.
\\
(ii) We denote with $\SS^{\RR}_{k,n}$ the set of self-energy clusters
of order $k$ and scale $[n]$ which do not contain any other
self-energy cluster nor any trivial node; we call them {\rm
renormalized self-energy clusters} on scale $n$.
\\
(iii) Given a self-energy cluster $T\in\SS^{\RR}_{k,n}$
we shall define the {\it self-energy value} of $T$ as the
matrix\footnote{${}^3$}{\nota This is a matrix because the self-energy
cluster inherits the labels $\g,\g'$ attached to the endnode of the
entering line and to the initial node of the exiting line.\vfil}
$$ \VV_{T}(\oo\cdot\nn;\e) = \fra{\e^{k}}{(k-1)!}
\Big( \prod_{\ell\in \L(T)} g^{[n_{\ell}]}_{\ell} \Big)
\Big( \prod_{\vvvv\in V(T)} F_{\vvvv} \Big) ,
\Eq(5.8) $$
where $g^{[n_{\ell}]}_{\ell}=g^{[n_{\ell}]}(\oo\cdot\nn_{\ell};\e)$.
Note that, necessarily, $n_\ell\le n$. The $k_T-1$ lines of the
self-energy cluster $T$ will be imagined as distinct and to carry
a number label ranging in $\{1,\ldots, k_T-1\}$.}

\*

The recursive definition of the matrices
$\MM^{[n]}(x;\e)$, $n\ge 1$, will be (if the series converges)
$$
\MM^{[n]}(x;\e)
=  \Big(\prod_{p=0}^{n-1}
\chi_{p}(D(x))\Big) \sum_{k=2}^{\io}
\sum_{T \in \SS^{\RR}_{k,n-1} } \VV_{T}(x;\e)
\defi \Big(\prod_{p=0}^{n-1} \chi_{p}(D(x))\Big) M^{[n]}(x;\e) ,
\Eq(5.9)$$
where the self-energy values are evaluated by means of the
propagators on scales $[p]$, with $p=0,\ldots,n$, which makes sense
because we have already defined the propagators on scale $[0]$ and
the matrices $\MM^{[0]}(x;\e)\=M_{0}$ (cf. Definition 2).

\*

With the above new definitions we have the formal identities
$$ h_{\nn,\g} = \sum_{k=1}^{\io} \sum_{\th \in
\Th^\RR_{k,\nn,\g}} \Val(\th),
\Eq(5.10) $$
where we have redefined the {\it value} of a tree $\th \in
\Th^\RR_{k,\nn,\g}$ as
$$ \Val(\th) = \fra{\e^k}{k!}\Big( \prod_{\ell\in \L(\th) }
g^{[\h_{\ell}]} (\oo\cdot\nn_{\ell};\e) \Big)
\Big( \prod_{\vvvv\in V(\th)} F_{\vvvv} \Big) ,
\Eq(5.11) $$
with $[\eta_\ell]=[\io],[0],\ldots,[\lis n_{0}-1],[\ge \lis n_{0}]$.
Note that  \equ(5.10) is not a power series in $\e$.

The statement in \equ(5.10) requires some thought,
but it turns out to be a tautology, see also Ref. \cita{GG},
and Ch. VIII in Ref. \cita{GBG}, {\it if one neglects convergence problems}
which, however, will occupy us in the rest of this paper.
A sketch of the argument is as follows.

Imagine that we have only scales $[\io],[0],\ldots,[n-1],[\ge n]$,
\ie we have performed the scale decomposition of the
propagators up to scale $[n-1]$ and we have not decomposed the
propagators on scale $[\ge n]$ and that we have checked the statement
\equ(5.9) and \equ(5.10) (trivially true for $n=0$).

Given a tree $\th\in\Th^\RR_{k,\nn,\g}$ with lines carrying labels
$[p]$ with $p=0,\ldots,n-1$ or $[\ge n]$ or $[\io]$,
we can split the propagators $g^{[\ge n]}(x;\e)$
as $g^{[n]}(x;\e)+g^{\big\{\ge n+1\big\}}(x;\e)$
as in \equ(5.6) with $g^{[n]}(x;\e)=\psi_n(D(x))g^{[\ge n]}(x;\e)$
and $g^{\big\{\ge n+1\big\}}(x;\e)=\chi_n(D(x))g^{[\ge n]}(x;\e)$.
In this way we get new trees {\it which in
general contain self-energy clusters of scale $[n]$}.
We can in fact construct infinitely many trees with
self-energy clusters of scale $[n]$ simply by {\it inserting} an
arbitrary number of them on any line $\ell$ with scale $\{\ge n+1\}$.

The values of the trees obtained by $k \ge 0$
such self-energy {\it insertions} on a given line of a tree in
$\Th^{\RR}_{k,\nn,\g}$ can be arranged into a geometric progression:
in fact they differ only by a factor equal to the value of the integer
power $k$ in $g^{\big\{\ge n+1\big\}}(x;\e) \big(M^{[n+1]}(x;\e)
g^{\big\{\ge n+1\big\}}(x;\e)\big)^{k+1} $ if $M^{[n+1]}(x;\e)$ is
defined as in \equ(5.9), where the $\VV_{T}(x;\e)$ are evaluated by using
as propagators $g^{[p]}(x;\e)$, with $0\le p\le n$ or $p=\io$,
for the lines carrying a scale label $[p]$.
{\it Summation over $k$ will simply change
$g^{\big\{\ge n+1\big\}}(x;\e)$ into $g^{[\ge n+1]}(x;\e)$ and at the
same time one shall have to consider only trees with no self-energy
cluster of scale $[n]$ nor of scale $[p]$ with $p<n$
and with lines carrying scale labels
$[\io],\ldots,[n]$ or $[\ge n+1]$.} In this way we prove \equ(5.10) for
all $n\le \lis n_0$: we could continue, but for the reasons outlined
in Section \secc(4), we decide to stop the resummations at this scale.

In other words the above is a generalization of the simple resummation
considered in Section \secc(3). The result is still {\it as formal as the
Lindstedt series we started with} even assuming convergence of the
series in \equ(5.9). In fact
the consequent expression for $\hh$
cannot even be, if taken literally, correct because as in
Section \secc(3) {\it the denominators in the new expressions
could even vanish because no lower cut-off operates on
the lines with scale $[\ge \lis n_0]$ as the third of \equ(5.7) shows}.

\*

To proceed we must first check that the series \equ(5.9) defining
$M^{[n]}(x;\e)$ are really convergent. In spite of the last comment
this will be true because in the evaluation of $M^{[n]}(x;\e)$ {\it
the only propagators needed have scales $[p]$ with $p\le n-1$} so
that, see the factors $\psi_{n}(D(x)),\chi_{n}(D(x))$ in \equ(5.7),
their denominators not only do not vanish but have controlled sizes
that can be bounded below proportionally to $x^{2}$ by \equ(4.4), \ie simply by a
constant times $C_{0}^2|\nn|^{-2\t_{0}}$, see \equ(1.3).

In Ref. \cita{GG} it has been shown {\it by a purely algebraic symmetry
argument} that, as long as one can prove convergence of the series
in \equ(5.9), the matrices $M^{[n]}(x;\e)$ are Hermitian and
$(M^{[n]}(x;\e))^{T}=M^{[n]}(-x;\e)$. Furthermore we should expect
that the eigenvalues of the matrix $\MM^{[\le n]} (x;\e)$
should be approximately located either near $0$ or
near $\e a_1,\ldots,\e a_s$ at least within $O(\e^2)$.


The expectation relies on Ref. \cita{GG} (see Eq. (3.25)) where the
following ``{\it cancellations result}\/'' is derived for $n_0$ large
enough (hence for $\e$ small because $2^{-2n_0-2}<\e a_s \le
2^{-2n_0}C^2_0$): we reproduce the proof in Appendix \secc(A3) below,
adapting it to the present notations.

\*

\0{\cs Lemma 2.}
{\it If $n_0$ is large enough and $n\le \lis n_0=n_0+\lis n$
(see \equ(4.2)) then the following properties hold.
\\
(i) The matrices $\MM^{[\le n]}(x;\e)$, $x=\oo\cdot\nn$,
are Hermitian and can be written as
$$\MM^{[\le n]}(x;\e) = \left( \matrix{
\MM^{[\le n]}_{\a\a}(x;\e) & \MM^{[\le n]}_{\a\b}(x;\e) \cr
\MM^{[\le n]}_{\b\a}(x;\e) & \MM^{[\le n]}_{\b\b}(x;\e) \cr} \right)
\Eq(5.12) $$
where the labels $\a$ run over $\{1,\ldots,r\}$ and $\b$ over
$\{r+1,\ldots,s\}$.
\\
(ii) One has $\MM^{\le [n]}(x;\e) = (\MM^{[\le n]}(-x;\e))^T$, so that
the eigenvalues of $\MM^{[\le n]}(x;\e)$ verify the symmetry property
$\l^{[n]}_{j}(x;\e)=\l^{[n]}_{j}(-x,\e)$, \ie they are
functions of $x^2$. \footnote{${}^4$}{\nota
For instance if $r=s=2$ and $f({\scriptstyle\aa},{\scriptstyle\bb})=
f_{\V0}({\scriptstyle\bb})+f_1({\scriptstyle\bb})
\cos\a_1 + f_2({\scriptstyle\bb})\cos\a_2$, {\it to lowest order in $x,\e$},
one has $M^{[\le n]}_{\a\a}(x;\e) = 3\e^{2}x^{2}(2\o_{u}^{4})^{-1}
[f_{u}^{2}({\scriptstyle\bb})+|\dpr_{\scriptscriptstyle\bb}
f_{u}({\scriptstyle\bb})|^{2}]\d_{u,v}$,
$M^{[\le n]}_{\a\b} = i\e^{2}x(2\o_{v}^{3})^{-1} \dpr_{\b_{v}}
[(f_{u}^{2}({\scriptstyle\bb}) + |\dpr_{\scriptscriptstyle\bb}
\f_{u}({\scriptstyle\bb})|^{2})]$,
and $M^{[\le n]}_{\b\b} = \e \dpr_{\scriptscriptstyle\bb}^{2}
f_{0}(\scriptstyle\bb)$, $u,v=1,2$.}
\\
(iii) Let $\dpr^\pm_x$ be right and left $x$-derivatives, then
$$\eqalignno{
&\|\MM^{[n]}(x,\e)\|\le B\, \e^2\,e^{-\k_1 2^{n/\t}},\quad
\|\dpr^\pm_x \MM^{[\le n]}(x,\e)\|\le B \e^2 a_s^{-1/2},\quad
\|\dpr^\pm_\e \MM^{[\le n]}(x,\e)\|\le B \, \e ,\cr
&\|\MM^{[n]}_{\a\a}(x;\e)\|\le B \,e^{-\k_1 2^{n/\t}}\,
\min\{\e^2,\e \, x^{2} a_{s}^{-1} \} , \cr
&\|\MM^{[n]}_{\a\b}(x;\e)\|\le B \,e^{-\k_1 2^{n/\t}}\,
\min\{\e^2,\e^{\fra32} |x|\,a_{s}^{-1/2} \} , & 
\eq(5.13) \cr
&\|\MM^{[n]}_{\b\b}(x;\e)\|\le B \,e^{-\k_1 2^{n/\t}}\,
\e^2 , \cr} $$
for $n\le \lis n_0$ and for suitable $\lis n_0$-independent constants
$B,\k_1,\t>0$; one can take $\t=\t_{0}$.}

\*

While $\k_1$ is dimensionless the constants $A',A,B$ have same
dimension (of a frequency square): this is the purpose of introducing
appropriate powers of $a_s$.

General properties of matrices and \equ(5.13) imply, see Appendix
\secc(A4),
\kern-3mm
$$\eqalign{
&A'< |\dpr_\e\l^{[n]}_{j}(x;\e)| <\, A,\qquad
a_s^{\fra12}
|\dpr_{x}^\pm\l^{[n]}_{j}(x;\e)|\, <\,   A \,\e^{2}, \qquad j>r, \cr
&A'< |\dpr_{\e} ( \l^{[n]}_{j}(x;\e)-\l^{[n]}_{i} (x;\e) ) |,
\qquad \kern3.cm i\ne j>r  , \cr
&|\l_{j}^{[n]}(x;\e)-\l_{j}^{[n-1]}(x;\e)|\le\e^2
B\,e^{-\k_1 2^{n/\t}},\kern2.4cm j>r,\cr
&|\l^{[n]}_{j}(x;\e)| < A\,\min\{\e^2,\e\,x^2 a_{s}^{-1} \} ,\qquad
\kern3.2cm j\le r, \cr}
\Eq(5.14)$$
where $A',A>0$ are $n,n_0$-independent constants, and $\t=\t_{0}$.

\*

\0{\it Remarks.}
(1) The first three bounds on the eigenvalues in \equ(5.14), follow
from the first line of \equ(5.13) by using the self-adjointness of the
matrices $\MM^{[\le n]}(x;\e)$; see Appendix \secc(A4).
The other bounds in \equ(5.13) imply the last bound in \equ(5.14);
see Appendix \secc(A4).
\\
(2) The natural domain of definition in $x$ of $\MM^{[n]}(x,\e)$, $n>0$, will
turn out to be $D(x)\le 2^{-2(n-1)}C^2_0$, but we imagine that it is
defined for all $x$ by continuing it as a constant from its limit
value. In fact this is not important because, as we shall see, only
the values of $\MM^{[n]}(x,\e)$ with $D(x)\le 2^{-2(n-1)}C^2_0$ enter
into the analysis.  Smoothness means differentiability in $\e\in I_C$
and a right and left differentiability in $x$. The lack of
differentiability in $x$, but the existence of right and left $x$
derivatives, is due to the fact that the function $D(x)$ admits right
and left derivatives: hence lack of differentiability in $x$ appears
as an artifact of the method.  This lack of smoothness (unpleasant but
inessential for our purposes) can be eliminated by changing $D(x)$
into a new $\widetilde D(x)$ which is smooth for $x^2$ between
successive ${\ul\l_{j}(\e)}$'s and, at the same time, it is bounded
above and below proportionally to $D(x)$. But this would make the
discussion needlessly notationally involved and we avoid it.
\\
(3) One should also  remark that, although we excluded  some values of
$\e$ (\ie we required $\e\in \EE_{\lis n_0-1}$, see \equ(3.3)), here all
$\e\in I_C$  are allowed. The restriction on $\e$
plays  no role in  the high  frequency resummations:  so far  its only
purpose is  to avoid  divisions by  $0$ and to assign  a finite  value to
contributions to $\hh$ from trees with with propagators on scale $[\ge
\lis  n_{0}]$ (which  could  be infinite  because  of the  lack of  an
infrared cut-off in their expressions; see the third of \equ(5.7)).
\\
(4) The bounds on the entries of $\MM^{[n]}(x;\e)$ in the second and
third lines of \equ(5.13) arise from cancellations that are checked
in Ref. \cita{GG} via a sequence of algebraic identities on the
Lindstedt series coefficients and {\it the real difficulty lies in the
proof of convergence}. The algebraic mechanism for the cancellations
is briefly recalled in Appendix \secc(A3), for completeness.
\\
(5) Loosely speaking (as mentioned in Section \secc(4)) the reason why
the above result holds with $\lis n_0$-independent constants, and why
its proof can be taken from Ref. \cita{GG}, is that if the
scales of the propagators are constrained to be $[p]$ with $p< \lis n_0$
the propagators denominators can be estimated by $2^{-2(\lis n+1)-2} x^2$
by \equ(4.4) and by the Remark (1) after Definition 1.
This means that one can proceed as in the hyperbolic tori cases in
which boundedness, from below, proportionally to $x^2$ of the
propagators denominators was the main feature exploited and {\it no
restriction} on $\e$ had to be required, other than suitable smallness.

\*

The lemma can be proved by imitating the convergence proof of the
KAM theorem, see for instance Ref. \cita{GG}; however in the following
Appendix \secc(A3) the part of the proof which is not reducible to a
purely algebraic check is repeated, for completeness.

We have therefore constructed a new representation of the formal series
for the function $\hh$ of the parametric equations for the
invariant torus: in it only trees with lines carrying a scale label
$[\io],[0],\ldots,[\lis n_0-1]$ or $[\ge \lis n_0]$ and
{\it no self-energy clusters} are present.
The above lemma will be the starting block of the construction that follows.

\*\*
\section(6,Renormalization: the infrared resummation)

\0Convergence problems still arise from the propagators
$g^{[\ge \lis n_0]}(x;\e)$, which become uncontrollably large for
$x=\oo\cdot\nn$ close to the eigenvalues of $M_{0}$
because the bound \equ(4.4) which allowed control
of the divisors in terms of the classical small divisors
(\ie in terms of $|x|$) does not hold any more.
Hence we must change strategy.
\*

\0{\cs Definition 6.}
{\it Given $d\times d$ Hermitian matrices $\MM^{[\le n]}(x;\e)$,
$n=\lis n_0,\lis n_0+1,\ldots$, with eigenvalues
$\l^{[n]}_{j}(x;\e)$, we introduce the following notations.
\\
(i) The sequence of {\rm self-energies}
$\ul\l^{[n]}_{j}(\e)$ is defined for $n\ge \lis n_0$ by
$$ \ul\l^{[n]}_{j}(\e)\defi
\l^{[n]}_{j} \Big(\Sqrt{\ul\l^{[n-1]}_{j}(\e)},\e\Big), \qquad
\ul\l_{j}^{[\lis n_0-1]}(\e)\defi \l^{[0]}_{j} ,
\Eq(6.1)$$
provided $\ul\l^{[n]}_{j}(\e)\ge0$, $n=\lis n_0,\lis n_0+1,\ldots$.
\\
(ii) The {\rm propagator divisors} are defined for $n\ge \lis n_0$ by
$$ \D^{[n]}(x;\e) \defi \Big| x^2 -
{\ul\l^{[n]}_{j(x)}(\e)} \Big| ,
\Eq(6.2) $$
where $j(x)$ is the label where the minimum of
$\Big|x^2-{\ul\l^{[n]}_{j}(\e)}\Big|$ is reached.
}

\*

\0{\it Remarks.}
(1) The self-energies are defined
recursively starting from those of the matrix $M_0$ whose first $r$
eigenvalues are $0$. Hence, {\it as long as one can extend the last of}
\equ(5.14) and as long as {\it the self-energies $\ul\l^{[n]}_{j}(\e)$
remain close to the original value $\l_{j}^{[0]}$},
as we shall check for $\e$ small enough, one has
$\ul\l^{[n]}_{j}(\e)=0$ for $j=1,\ldots,r$
and $\ul\l^{[n]}_{j}(\e)>0$ for $j>r$.
\\
(2) Under the same conditions and if $\D^{n]}(x;\e)\simeq 2^{-2n} C_0^2$
the label $j(x)$ depends only on $M_0$, hence it is $n$-independent, and
furthermore it is constant at $x$ fixed, as $\e$ varies in the
intervals $I$ introduced in Definition 1 
(because for large $n$ the frequency $x$ is constrained to be close to
one of the $\ul \l_j^{[n]}(\e)$).
\\
(3) The name of {\it propagator divisor} assigned to $\D^{[n]}(x,\e)$
in \equ(6.2) reflects its later use as a lower bound on the
denominator of a propagator, see Remark (7) to the inductive
assumption below.

\*

By repeating the analysis of Section \secc(4) we can represent the
function $\hh$ via sums of values of trees in which lines can carry
scale labels $[\io],[0],\ldots,[\lis n_0-1],[\lis n_0],[\lis
n_0+1],\ldots$ and which contain no self-energy clusters and no
trivial nodes (\ie are renormalized trees, see Definition 5 in Section
\secc(5)). The new propagators will be defined by the same procedure
used to eliminate the self-energy clusters of scales $[n]$ with $n\le
\lis n_0-1$. However we shall now determine the scale of a line in
terms of the recursively defined $\D^{[n]}(x;\e)$ rather than in terms
of $D(x)$: the latter becomes too rough to resolve the separation between
the eigenvalues and their variations.

Let $X_{\lis n_0-1}(x)\defi
\prod_{m=0}^{\lis n_0-1}\chi_m(D(x))$, $Y_{n}(x;\e)\defi
\prod_{m=\lis n_0}^{n}\chi_{m}(\D^{[m]}(x;\e))$ for $n\ge \lis n_0$ and
$Y_{\lis n_0-1}\=1$: the definition of the new propagators will be
$$\eqalign{
g^{[\lis n_0]} \defi & \kern3mm
X_{\lis n_0-1}(x)\, \psi_{\lis n_0}(\D^{[\lis n_0]}(x;\e))\,
(x^2-\MM^{[\le \lis n_0]}(x;\e))^{-1} , \cr
g^{[\lis n_0+1]} \defi & \kern3mm
X_{\lis n_0-1}(x)\,\chi_{\lis n_0} (\D^{[\lis n_0]}(x;\e))\,
\psi_{\lis n_0+1}(\D^{[\lis n_0+1]}(x;\e)) \,
(x^2-\MM^{[\le \lis n_0+1]}(x;\e))^{-1} , \cr
&\ldots\cr
g^{[n]} \defi & \kern3mm
X_{\lis n_0-1}(x)\,Y_{n-1}(x;\e)\,\psi_{n}(\D^{[n]}(x;\e))
\,(x^2-\MM^{[\le n]}(x;\e))^{-1} , \cr}
\Eq(6.3) $$
and so on, using indefinitely the identity $1\=\psi_n(\D^{[n]}(x;\e))+
\chi_n(\D^{[n]}(x;\e))$ to generate the new propagators.

\*

In this way we obtain a formal representation of $\hh$ as a sum of
tree values in which only renormalized trees (\ie with neither trivial nodes
nor self-energy lines, see Definition 4 in Section \secc(4))
and in which each line $\ell$ carries a {\it scale label}
$[n_\ell]$. This means that we can formally write $\hh$ as in
\equ(5.10), with $\Val(\th)$ defined according to \equ(5.11), but now
the scale label $[n_{\ell}]$ is such that $n_{\ell}$ can assume all
integer values $\ge 0$ or $\io$, and no line carries a scale label like
$[\ge n]$: {\it only scale labels like $[n]$ are possible}.
The corresponding propagators $g^{[n_{\ell}]}(\oo\cdot\nn_{\ell};\e)$
will be defined as follows. 

\*

\0{\cs Definition 7.}
{\it Given a sequence $\MM^{[\le m]}(x;\e)$ as in Definition 6,
$m\ge1$, let $\MM^{[n]}(x;\e)=\MM^{[\le n]}(x;\e)-\MM^{[\le
n-1]}(x;\e)$ with  $\MM^{[\le0]}\=\MM^{[0]}\=M_0$
so that $\MM^{[\le n]}(x;\e)=\sum_{m=0}^{n}\MM^{[m]}$ $(x;\e)$. Setting
$\D^{[n]}(x;\e)\=D(x)$ if $n\le \lis n_0$, define for all $n\ge 0$
$$ g^{[n]}(x;\e)=
\fra{\psi_{n}(\D^{[n]}(x;\e))\prod_{m\ge 0}^{n-1}
\chi_m(\D^{[m]}(x;\e))}{x^2-\MM^{[\le n]}(x;\e)} .
\Eq(6.4)$$
(for $n=0$ this means $\psi_0(D(x))\,(x^2-M_0)^{-1}$).  We say that
$g^{[n]}_{\ell}= g^{[n]}(\oo\cdot\nn_{\ell};\e)$ is a propagator with
scale $[n]$.  The matrices $\MM^{[m]}(x;\e)$ will be defined as in
Section \secc(5) for
$n\le\lis n_0$ and will be defined recursively also for $n>
\lis n_0$ in terms of the self-energy clusters $\SS^\RR_{k,n-1}$
introduced in Definition 4, Section \secc(5), setting for $n>\lis
n_0$ (see \equ(5.9))
$$ \MM^{[n]}(x;\e)= \Big(\prod_{m=0}^{
n}\chi_m(\D^{[m]}(x;\e))\Big)\sum_{k=2}^{\io}
\sum_{T \in \SS^{\RR}_{k,n-1} } \VV_{T}(x;\e) ,
\Eq(6.5)$$
where the self-energy values $\VV_{T}(x;\e)$ are evaluated by means
of propagators on scales less than $[n]$. Note that 
we have already defined (consistently with \equ(6.5)))
the matrices $\MM^{[\le n]}$ with $n\le \lis n_0$ and
the propagators on scale $[\io],[0],\ldots,$ $[\lis n_0-1]$
(so that \equ(6.4) defines also $g^{[\lis n_0]}(x;\e)$).}
 
\*

\0{\it Remark.} (1) Some propagators may vanish being proportional to
a product of cut-off functions. If a propagator of a line
has scale $[n]$ and does not vanish then, see \equ(6.4),

$$2^{-2(n+1)}C^2_{0}\le\D^{[n]}(x;\e)\le 2^{-2(n-1)}C^2_{0}\Eq(6.6)
$$
\\
(2) Our definitions of the matrices $\MM^{[\le n]}(x;\e)$ for $n> \lis
n_0$ will be such that given the node harmonics of a tree 
hence by \equ(6.6) the scale $[n]$ that is
attributed to a line can only assume up to two consecutive values
unless the propagator (hence the tree value) vanishes.
\\
(3) We may and shall imagine that scale labels are assigned
arbitrarily to each line of a given tree {\it with the constraint that
no self energy clusters are generated}; however the tree will have a
non-zero value only if the scale labels are such that all propagators
do not vanish. This means that only up to two consecutive scale labels
can be assigned to each line if the tree value is not zero. The
``ambiguity'' on the scale labels for a line is due to the
use of the non-sharp $\chi$ and $\psi$ functions of Fig. 2.

\*

We make an inductive assumption on the propagators on the scales
$[m]$, $0 \le m < n$.

\*

\0{\cs Inductive assumption.}
{\it Let $\lis n_0\=n_0+\lis n$ (see \equ(4.2)) and suppose
$n_0$ large enough; then
\\
(i) For $0\le m\le n-1$ the matrices $\MM^{[m]}(x;\e)$
are defined by convergent series for all $\e\in I_C$ and, for all $x$,
they are Hermitian, and $\MM^{[m]}(x;\e)=(\MM^{[m]}(-x,\e))^T$.
Furthermore they satisfy the same relations as \equ(5.13),
hence \equ(5.14), with $n$ replaced by $m$, for all $0<m<n-1$,
with suitably chosen (new, possibly different) constants $\k_1,A,A',B,\t$.
One can take $\t=2\t_{1}$.
\\
(ii) There exist $K>0$ and open sets $\EE^o_m$, $m=0,\ldots, n$, with
$\EE^o_{m}\subset I_C$, such that, defining recursively
$\ul\l^{[m]}_{j}(\e)$ in terms of $\ul\l^{[m-1]}_{j}(\e)$ for $m=\lis
n_0,\ldots,n-1$ by (i) in Definition 6 above, while setting
$\ul\l^{[m]}_{j}(\e)\=\l^{[0]}_j$ for $m=0,\ldots,\lis n_0-1$,
and defining $\t_{1} \defi \t_{0} + r + 1$, see
\equ(3.4), one has for $\e\not\in \EE^o_m$
$$ \eqalign{
\G^{[m]}(x;\e)
 &=\min \Big\{\min_{j}
\Big| x\pm\Sqrt{\ul\l^{[m]}_{j}(\e)} \Big| \,,\,
\, \min_{j\ge i}\Big| x \pm
\Sqrt{\ul\l^{[m]}_{j}(\e)} \pm
\Sqrt{\ul\l^{[m]}_i(\e)} \Big|\,\Big\} \ge
2^{-\fra12 m} \fra{C_{0}}{|\nn|^{\t_{1}}} , \cr
|\EE_m^o|&\le K 2^{-\fra12 m} C^2 , \cr}
\Eq(6.7) $$
for all $m\le n-1$ and all $x$.}

\*

\0{\it Remarks.} Assuming validity of the hypothesis for $m< n$ we
note a few of its implications.
\\
(1) So far we have only checked the hypothesis for scales $[m]$ with
$m\le \lis n_0$, as expressed by Lemma 2 in Section \secc(5), \ie for
the high frequency propagators. If (i) is proved also for
$m=n$ then we can impose \equ(6.7)
immediately by excluding a set $\EE_m^o$ of $\e$'s of
measure estimated by $2^{-\fra12 m} C^2 K$ with $K$ a constant that
can be bounded in terms of $A',A$ by introducing the constants
$\r_{m}$ and $\r'_m$ as in \equ(A2.1), with
$\ul{\l}^{[0]}_{j}(\e)$ replaced by $\ul{\l}^{[m]}_{j}(\e)$, and
proceeding as done in Appendix \secc(A2) for the case $n\le \lis
n_0$. Note that since the self-energies
$\ul\l_j^{[m]}(\e)$ are $\=\ul\l_j^{[0]}(\e)$ for all
$m=0,\ldots,\lis n_0-1$ one will have, for such $m$'s,
$\EE^o_m\=I_C/\EE_{\lis n_0-1}$, see \equ(3.3). It is very important
to keep in mind, in the above argument, that the self-energies
either are $0$ (for $j\le r$) or are close within $O(\e^2)$ to the positive
eigenvalues of $M_0$, and they are {\it differentiable} in $\e$ and
to the right and left of each  $x$ by (i); see \equ(5.14).
\\
(2) If a line with a scale $[n]$ has
vanishing propagator (\ie $g^{[n]}(x;\e)=0$ because of the $\chi,\psi$
cut-off functions in the definition \equ(6.4))
but $n$ differs at most by one unit from the integer $n'$
such that $g^{[n']}(x;\e) \neq 0$. Thus if we consider
$\D^{[n]}(x,\e)$ we can bound it by changing the
inequalities \equ(6.6) into $C^2_{0} \,2^{-2(n+2)} <
\D^{[n]}(x;\e) \le C^2_{0}\, 2^{-2(n-2)}$. The remark will be
useful later when we shall exploit it in the discussion of the
cancellations which we shall study to check the inductive hypothesis.
\\
(3) By \equ(5.13) and \equ(5.14), and (I) in Appendix \secc(A4),
we deduce that ${\l_{j}^{[m]}(x;\e)}$, hence ${\ul\l_{j}^{[m]}(\e)}$,
do not change by more than $C\,B\, \e^2
\sum_{n\ge \lis n_0} e^{-\k_1 2^{n/(2\t_{1})}}$,
with respect to ${\l^{[0]}_j(\e)}$, if $\e < \lis \e_{1}$
(and $m\ge \lis n_0$).
\\
(4) Hence if $\e$ is small enough the self-energies, \ie
${\ul\l^{[m]}_{j}(\e)}$, have distance bounded above by $2{a_s\e}$ and
below by $\fra12{\e} \min \big\{ {a_1},\min_{j}\{ a_{j+1}- a_1 \} \big\}
= 2 \r \,{\e a_{s}}$ with $\r$ defined in \equ(4.2),
if $\e$ is small enough, say $\e< \lis\e_{2}$.
\\
(5) Therefore by Remark (4) we see that the distance of $|x|^2$ from the
closest value ${\ul\l^{[m]}_{j}(\e)}$ is smaller than one fourth,
up to corrections $O(\e^{2})$, the distance between
the distinct values of ${\ul\l^{[m]}_{j}(\e)}$,
if $m$ is large enough compared to $n_0$, \ie if $2
C_{0}^2 2^{-2m} < \r {\e a_{s}}$ (or $m-n_0\ge \lis n$ as implied by
the definition \equ(4.2) of $\lis n$).  This
means that $j(x)$ is $m,\e$-independent and it coincides
with the label minimizing $\big| x^2-{|\l^{[m]}_{j}(x;\e)|}\big|$
for all $m\ge\lis n_0$ and all $\e\in I$.
\\
(6) $\ul\l^{[\lis n_0-1]}_{j}(\e)\=\l^{[0]}_{j}$ are $x$-independent
and, by their definition,
the same remains true for {\it all } $\ul\l^{[m]}_{j}(\e)$. The
self-energy
$\ul\l^{[m]}_{j}(\e)$ will be thought of as a reference position for
the $j$-th eigenvalue on scale $[m]$, $m\le n-1$.
\\
(7) As noted in Remark (5) the quantity
$\big| x^2-{\l^{[n]}_{j(x)}(x;\e)}\big|$
is the smallest denominator appearing in the value of the
propagator of a line with momentum $\nn$
if $g^{[n]}(x;\e)\ne0$ (here $x=\oo\cdot\nn$).
{\it The key to the analysis is the check
that the quantities $\D^{[n]}(x;\e)$ can be used to
bound below the denominators of the non-vanishing propagators of scale
$[n]$}. If $\l^{[n]}_{j(x)}(x;\e)<0$ one has
$x^{2}-\l^{[n]}_{j(x)}(x;\e)\ge x^{2}$, so that the assertion
is trivially satisfied: therefore the really interesting case is
when $\l^{[n]}_{j(x)}(x;\e)\ge 0$ (which includes the cases $j(x)>r$).
If $x$ has scale $[n]$ with $n\ge\lis n_0$ one has
$$\eqalign{
& \Big||x|- \Sqrt{|\l^{[n]}_{j(x)}(x;\e)|} \Big| \ge
\Big| |x|- \Sqrt{\ul\l^{[n]}_{j(x)}(\e)} \Big| -
\Big| \Sqrt{\ul\l^{[n]}_{j(x)}(\e)} -
\Sqrt{\l^{[n]}_{j(x)}(x;\e)} \Big| \cr
& \qquad \ge \fra12\Big| |x|- \Sqrt{\ul\l^{[n]}_{j(x)}(\e)} \Big| +
2^{-(n+3)}C_{0}-\Big|\Sqrt{\l^{[n]}_{j(x)}
\big( \Sqrt{\ul\l^{[n-1]}_{j(x)}(\e)},\e \big) } -
\Sqrt{\l^{[n]}_{j(x)}(x;\e)} \Big| \cr
& \qquad \ge \fra12 \Big||x|- \Sqrt{\ul\l^{[n]}_{j(x)}(\e)}\Big| ,\qquad
\tto \qquad \|x^2-\MM^{[n]}(x,\e)\|\ge \fra1{2^{3}}\sqrt{\fra {a_1}{a_s}}
\,\big|x^2-\ul\l^{[n]}_{j(x)}(\e)\big| , \cr}
\Eq(6.8) $$
having used the lower cut-off $\psi_{n}(\D^{[n]}(x;\e))$ in the propagator
(see \equ(6.3)) to obtain the first two terms in the second line
(and added a further factor $2^{-1}$ in order to extend the
result also to the propagators considered in Remark (3)),
while the upper cut-off $\chi_{n-1}(\D^{[n-1]}(x;\e))$
has been used to obtain positivity of the difference between
the second and third terms in the second line,
after applying \equ(5.14), for $n\ge \lis n_{0}$, to get
$$ {\max_x |\dpr^{\pm}_x\l_{j(x)}^{[n]}(x;\e)|}
\le \widehat B\,\e^{2} ,\ j(x)>r,\qquad
{|\l_{j(x)}^{[n]}(x;\e)|}\le \widehat B\,\e\,|x|^{2} \,\le\, \e\,
C_{0}^{2}\, 2^{-2n},\ j(x)\le r ,
\Eq(6.9)$$
for some $\widehat B$, so that the last term in the second line of
\equ(6.8) can be bounded above proportionally to $\e 2^{-n} C_{0}$.
Hence the first inequality in the
last line of \equ(6.8) follows if $\e$ small
enough, say $\e\le\lis \e_{3}$ for some $\lis \e_{3}$, {\it fixed
{\it independently of $n$}}. The latter constraint can be achieved
simply by taking $n_0$ large enough, see \equ(3.2).
The last implication follows from \equ(6.9)
if $j(x)\le r$ because $\ul\l^{[n]}_{j(x)}(\e)=0$.
Otherwise if $j(x)>r$ and $|x|,\sqrt{\ul\l^{[n]}_{j(x)}(\e)},
\sqrt{\ul\l^{[n]}_{j(x)}(x,\e)}\in[\fra12\sqrt{\e\,a_1},2\sqrt{a_s\,\e}]$
one has $(|x|+\sqrt{\ul\l^{[n]}_{j(x)}(x,\e)})/
(|x|+\sqrt{\ul\l^{[n]}_{j(x)}(\e)}) \ge 2^{-2} \sqrt{a_1/a_s}$, 
as long as $\e<\lis \e_{2}$ (see Remark (4) above):
implying again the \equ(6.8).
Hence $\D^{[n]}(x;\e)$ can be effectively used to estimate the size of the
non-vanishing propagators which is, therefore, closely related to the
scale of the corresponding lines.
\\
(8) The Diophantine condition \equ(3.3) and \equ(6.7) will play from
now on a key role. We begin by remarking that if the inductive
hypothesis is proved {\it all lines will eventually acquire a well
defined scale label}: in fact fixed $x$ one cannot have $\D^{[n]}(x,\e) \le
2^{-2n} C_0^2$ for all $n$ because this implies\footnote{${}^5$}{\nota
As $|a^2-b^2|< c^2$ implies $|a-b|<c$ for $a,b,c>0$.} $||x|-
\Sqrt{\ul\l^{[n]}_{j(x)}(\e)}|< 2^{-n} C_0$ , which soon or later becomes
incompatible with the first of \equ(6.7). This explains why there
is no trace left of the propagators $g^{[\ge n]}(x,\e)$.

\*

To estimate the corrections to the self-energy as $n$ increases it is clear
that we must estimate the size of $\MM^{[n]}(x;\e)$. For this purpose
we need the following result.

\*

\0{\cs Lemma 3.} {\it There is $\lis\e$ small and
constants $\k_1,A,A',B$ such that if $\e<\lis\e$ and the inductive
hypothesis is assumed for $\,0\le m\le n-1$
then the matrix $\MM^{[n]}(x;\e)$ can be bounded by \equ(5.13)
and the inductive hypothesis holds for $m=n$.}

\*
Hence the hypothesis holds for all
$n$ since we have already checked it for $n=0,\ldots,\lis n_0$ (Lemma
2): the new constants $\k_1,A,A',B$ will be different from the ones
determined in Lemma 2.
\*

\0{\it Proof.} For $n\le \lis n_{0}$ the bound \equ(5.13)
is covered by Lemma 2. So we can assume $n \ge \lis n_{0}+1$.
Suppose first $\e\in \cap_{m=\lis n_0-1}^{n-1}\EE_m$,
with $\EE_{m}=I_{C}\setminus \EE_{m}^{o}$, so that the
Diophantine property \equ(6.7) holds for all $m\le n-1$.
Consider a self-energy cluster $T$ in
$\cup_{k=2}^\io\SS_{k,n-1}^\RR$. If the entering and exiting lines
(with propagators of scale $[\ge n]$) have momenta $\nn$
we begin by showing that
$$ \sum_{\vvvv\in V(T)} |\nn_{\vvvv}| > 2^{(n-6)/(2\t_{1})} .
\Eq(6.10) $$
Indeed the cluster contains at least one line $\ell=\ell_{\vvvv}$
with propagator which we can suppose to be not
vanishing and which has scale $[n-1]$.
We can write $\nn_{\ell}=\nn_{\ell}^{0}+\s_{\ell} \nn$, where
$\s_\ell=0,1$ and we set
$\oo\cdot\nn=x$,
$ \nn_{\ell}^{0} = \sum_{\wwwww \in V(T) \atop \wwwww \preceq \vvvvv}
\nn_{\wwww} $, and finally $x_{\ell}=\oo\cdot\nn_{\ell}$.

Since the line $\ell$ is not on scale $[n-2]$ (as it is on scale
$[n-1]$) it follows from
\equ(6.3) that
$$ \big| |x_{\ell}|- \Sqrt{\ul\l^{[n-2]}_{j(x_{\ell})}(\e)} \big|
\le 2^{-(n-2)} C_{0},
\Eq(6.11) $$
Therefore if \equ(6.10) does not hold and if $\s_\ell=0$,
by the first part of the Diophantine conditions \equ(6.7), one finds
$\big||x_\ell|-\Sqrt{\ul\l^{[m]}_{i}(\e)}\big|> C_{0}2^{-m/2}
2^{-(n-6)/2}$ for all $m\le n-1$ and for all $1\le i \le d$,
which would be in contradiction with \equ(6.11).

If instead $\s_\ell=1$ we shall use the second part of
the Diophantine conditions \equ(6.7) and get a contradiction.
Remark that $x$ can be assumed to be on scale $[q]$ with $q\ge n$
because of the cut-off functions in \equ(6.5) so that one has
$\big||x|-\Sqrt{\ul\l^{[p]}_{j(x)}(\e)}\big|\le C_0 2^{-p}$ for
$p\le n-1$. Hence if $x_\ell$ satisfies \equ(6.11) we get,
by assuming that \equ(6.10) does not hold,
$$ \eqalignno{
2^{3-n}C_{0} & \ge
\Big| |x_{\ell}|-\Sqrt{\ul{\l}^{[n-2]}_{j(x_{\ell})}(\e)} \Big| +
\Big| |x|-\Sqrt{\ul{\l}^{[n-2]}_{j(x)}(\e)} \Big| &\eq(6.12)\cr
& \ge \Big| x_{\ell}-x + \h_{\ell}
\Sqrt{\ul{\l}^{[n-2]}_{j(x_\ell)}}(\e) +
\h \Sqrt{\ul{\l}^{[n-2]}_{j(x)}(\e)} \Big| \cr
& \ge { C_{0} \over 2^{(n-2)/2} |\nn_{\ell}-\nn|^{\t_{1}} }
= { C_{0} \over 2^{(n-2)/2} |\nn_{\ell}^{0}|^{\t_{1}} }
\ge 2^{4-n} C_{0} , \cr}
 $$
for some $\h,\h_\ell=\pm1$, which again leads to a contradiction,
so that \equ(6.10) holds also in such a case.

Every node factor contributes to $\MM^{[n]}$
a factor $f_{\nn_{\vvvvv}}$ bounded by $F_0 e^{-\k_{0} |\nn_{\vvvv}|}$;
there are $\le (4d^{2})^k k!$ self-energy clusters,
$4^{k}$ scales (for each line there are only two scales
for which the propagator is not zero, and one has to allow
also a scale different by one unit from that which
corresponds to have a nonvanishing propagator,
see Remark (3) after the inductive assumption),
and $\NN_{m}(T)$ lines of scale $m=0,1,\ldots,n$ in each self-energy
cluster $T$ contributing to $M^{[n]}(x;\e)$
and not to the $M^{[m]}(x;\e)$, with $m<n$.
Thus the bound on $M^{[n]}(x;\e)$ is
$$ G_0\sum_{k=2}^\io \e^{k} G^{k}_{1}
e^{-\fra12\k_0\sum_{\vvvvv\in V(T)}|\nn_{\vvvvv}|}
e^{- G_{2} 2^{n/(2\t_{1})}} \prod_{m=0}^n 2^{2m\NN_{m}(T)} ,
\Eq(6.13) $$
for suitable constants $G_0,G_{1}, G_{2}$, explicitly computable by the
above remarks; for $k<4$ the exponent of $\e$ can be replaced by $2$,
see Remark after Definition 4. The estimate of the number $\NN_{m}(T)$
is given in Appendix \secc(A3) (cf. in particular
Section \secc(A3.4)), and gives $\NN_{m}(T) \le E_{m}
\sum_{\vvvv\in V(T)}|\nn_{v}|$, with $E_{m}=2^{(6-m)/(2\t_{1})}$ and
$\t_1$ in \equ(3.4), which shows convergence of
the series in \equ(6.13) if $\e$ is small enough, say $\e < \lis \e$.

We can and shall assume that $\lis\e$ does not exceed
$\min\{\lis\e_{1},\lis\e_{2},\lis\e_{3}\}$, with $\lis\e_{1}$,
$\lis\e_{2}$ and $\lis\e_{3}$ introduced earlier (see Remarks 3, 4, 7
after the inductive hypothesis). {\it The rest of the argument repeats
the analysis in Appendix \secc(A3) with minor notational changes}:
we only hint at the details in Appendix \secc(A3.4).

Under the considered hypotheses the matrices $\MM^{[n]}(x;\e)$ are
well defined, by the above discussion on convergence of the defining
series on the set $\cap_{m=\lis n_0-1}^{n-1}\EE_m$.
The symmetry in item (i) is due to algebraic identities valid
for the Lindstedt series.
They are detailed  in Ref. \cita{GG}, Appendix A5,
for $\e<0$: being of algebraic nature the argument does not
depend on the sign of $\e$ and it holds unchanged in the present case.

The second and third lines of inequalities in \equ(5.13) embody the
cancellations. We need to check the cancellations to make sure
for instance that the structure of the matrix $\MM^{[n]}(x;\e)$
preserves the eigenvalues and the Whitney smoothness:
a danger being that the first $r$ eigenvalues become ``detached''
from $0$, \ie no longer can be bounded by $\e x^2$, see \equ(5.14).
For instance a bound like $O(\e^2)$ would
not be enough as it would imply that the self-energies
$\ul\l^{[n]}_j(\e)$ may become different from zero for $j\le r$.

Since the function $\MM^{[n]}(x;\e)$ is defined on the complement of a
dense open set differentiability in the sense of Whitney can be proved
(as usual) by computing a formal derivative and then showing that it
is continuous and that it can also be used as a bound in
interpolations.\footnote{${}^6$}{\nota More precisely in its simplest
form Whitney's theorem states that if $F(x)$ is a function defined on
a closed set $C$ of the interval $[0,1]$ and if there is a continuous
function $F'(x)$ defined on $C$ and if for some $\g>0$ and all $x,y\in C$
one has $|F(y)-F'(x)(y-x)|< \g |x-y|$ (we call this an {\it
interpolation bound}) then there is a continously differentiable
function $\lis F(x)$ extending $F$ to $[0,1]$ and with derivative
$\lis F'(x)$, with $\max|\lis F'(x)|<\g$, extending $F'(x)$.\vfil}

The computation of the formal derivatives proceeds as the
computation of the actual derivatives done in the proof of Lemma 2
(in Appendix \secc(A3)). One proves formal right and left continuous
differentiability of the matrices $\MM^{[n]}(x;\e)$
on the closed set $\cap_{m=\lis n_0-1}^{n-1}\EE_m$
simply by differentiating term by term the value of
each cluster contributing to $\MM^{[n]}(x;\e)$. This involves
differentiating matrices like $(x^2-\MM^{[\le p]}(x;\e))^{-1}$, \ie the
matrices $\MM^{[p]}(x;\e)$ with $p<n$, which are differentiable
by the inductive assumption, or it involves differentiating the
cut-off functions $\psi_p,\chi_p$ and the locations
$\ul\l^{[p]}_j(\e)$ with $j>r$ (because $\ul\l^{[p]}_j(\e)\=0$
for $j\le r$) which appear in the form $\D^{[p]}(x,\e)$ in the
arguments of the cut-off functions. All such quantities are
differentiable in $\e$ and right and left
differentiable in $x$ by the inductive assumption;
furthermore all terms arising from differentiation either of
$\MM^{[p]}(x;\e)$ or $\ul\l^{[p]}_j(\e)$, with $p<n$,
appear multiplied by some power of $\e$, so that
the inductive assumption is found to hold also for $p=n$
(for a similar discussion see Ref. \cita{Ge}).

Note that $\D^{[n]}(x;\e)$ depend on $j(x)$ but as $\e$ varies within
the interval $I$, see (ii) in definition 1, $j(x)$ is not only
$\e$-independent but it is also constant in $x$ for $x$ varying in
small intervals near the eigenvalues of $M_0$ and, therefore,
in intervals widely spaced because $n\ge\lis n_0$: this is due
to the cut-off functions which force $x$ to be close to a single eigenvalue if
the propagator of the corresponding line is different
from $0$. Hence for $n\ge \lis n_0$ we do not have to differentiate
the function $j(x)$ (neither with respect to $x$ nor with respect to
$\e$ from which it does not depend);
for $n<\lis n_0$ the function $j(x)$ is constant to the
right and to the left of every point.

The $n$-independence of the constants $A',A,B$ appearing in the
inductive hypothesis is proved word by word as the corresponding
argument in Appendix \secc(A3); the constant $\k_1$ has been estimated
above (see $G_2$ in \equ(6.13)) and is $n$-independent.

The interpolation bound, see footnote ${}^6$, necessary for
defining the Withney derivatives, holds because in
comparing two contributions to $\MM^{[n]}(x;\e)$ with different $x$ or
different $\e$ the difficulty might only come from the comparison of
$(x_\ell^2-\MM^{[\le p]}(x_\ell,\e))^{-1}$ evaluated
at two different points and for one line $\ell$ at a time:
this can be done algebraically by using the resolvent identity
$$ \eqalignno{
& \left( x_{\ell}^2-\MM^{[\le p]}(x_\ell,\e) \right)^{-1} -
\left( {x'}_{\ell}^2- \MM^{[\le p]}(x'_\ell,\e') \right)^{-1}
= \left( x_{\ell}^2-\MM^{[\le p]}(x_{\ell},\e) \right)^{-1}\cdot
& \eq(6.14) \cr
& \cdot\left( {{x'}}_{\ell}^2 - x_\ell^2+ \MM^{[\le p]}(x'_{\ell},\e') -
\MM^{[\le p]}(x_{\ell},\e)\right)
\left( {x'}_\ell^2-\MM^{[\le p]}(x'_{\ell},\e') \right)^{-1} , \cr} $$
which involves only denominators evaluated at $x,\e$'s which are in
the set where they are controlled by the \equ(6.7) and therefore can
be estimated in the same way as the formal derivatives. The Withney
extension is therefore possible keeping control of the
bounds for all $\e$'s (small as above) and $x$.
The dependence on $x$ may involve the functions $D(x)$ (for
$p\le \lis n_0-1$) so that the differentiability in $x$ will be possible only
to the right and to the left of each point (this involves a natural
generalization of Whitney's theorem).

The cancellations analysis (\ie the proof of the second and third
inequalities in \equ(5.13)) is inductive and has been
performed several times in the literature, see Refs. \cita{Ga} and
\cita{GG}: in Appendix \secc(A3) we have repeated it following the
version in Ref. \cita{GM1} with some minor modifications.  The same
proof applies to the present case (being a purely algebraic check).

The inequalities \equ(5.13) imply the \equ(5.14) and therefore we get
differentiability of the matrices $\MM^{[\le n]}(x;\e)$
and of the self-energies.
This allows us to impose validity of \equ(6.7) by excluding
a few more values of $\e$ by Remark (1) to the inductive hypothesis.

Therefore we conclude that $\MM^{[n]}(x;\e)$ is defined
and verifies \equ(5.13) (with suitably chosen constants $\k_1,A',A,B$)
in the same domain $\e<\lis \e$, where the matrix $\MM^{[\le p]}(x;\e)$
is already defined for $p\le n-1$. {\it Of course $\MM^{[n]}$ 
will be relevant for
our analysis only on the set $\cap_{m=\lis n_0-1}^n \EE_m$} and the
extension outside such set is only useful to simplify the analysis
as it allows us to use freely interpolations formulae, mainly to check
\equ(6.7).  The matrix
$\MM^{[\le n-1]}(x;\e)$ verifies the inductive assumption although it
has physical meaning only for $\e\in\cap_{m=\lis n_0-1}^\io\EE_n$,
where $\EE_{n}$ is the domain in which \equ(6.7) holds for $m\le n$. \qed

\*

Having checked that the series defining the
$\MM^{[\le n]}(x;\e)$, hence the self-energies, converge and verify the
bounds in the inductive hypothesis we still have
to check that the fully renormalized series for $\hh$,
which has thus been shown to make sense term by term, converges and
that its sum is actually a function $\hh$ satisfying the
equations for the parametric representation of invariant tori.

To study convergence we can take again advantage of the method,
already used in the proof
of Lemmas 2,3 above to estimate the number of lines on scale $n$ in a
self-energy cluster containing no self-energy clusters.  Indeed also
for renormalized trees one can prove a bound like $\NN_{m}(\th)\le
E_{m}\sum_{\vvvv\in V(\th)}|\nn_{v}|$ for the number $\NN_{m}(\th)$ of
lines on scales $m$ contained in $\L(\th)$ with $E_m$ fast decreasing with
$m$: $E_m\defi 2^{(6-m)/(2\t_{1})}$ (see Appendix \secc(A3)).
Hence convergence in the region $\EE\in\cap_{n=\lis n_0-1 }^\io \EE_n$
follows because if we only sum values of trees without
self-energy clusters then we can use the above bound on $\NN_{m}(\th)$.

The renormalized trees may contain lines of scales $[\infty]$ which so
far played no role. Their propagators are bounded below by a constant
times $\e$; however their number cannot be larger than $\fra12k$ in
trees of order $k$ (see Remark (5) in Section \secc(2)).
Therefore they may reduce the factor $\e^k$
normally present in the value of a graph with $k$ nodes to $\e^{\fra12
k}$; hence this will not affect the convergence of the series other
than by putting a more severe constant on the maximum value of $\e$.

The set $\EE^0_n$, complement of $\EE_n$ in $I_C$, has measure
estimated by $C^2 2^{-n/2} K$ for $\e\in[(\fra12 C)^2, C^2]$ $=$ $I_C$.
Since $C=2^{-n_0}C_0$ and $n\ge \lis n_0-1>n_0$ 
this is a very small fraction of the interval $I_C$ and the smaller
the closer is $I_C$ to $0$. This means 
that the set of $\e$'s for which the whole construction can be
performed has $0$ as a density point. Note that the resummation just
defined  is a real resummation of our series only for
$\e\in\cap_{n=\lis n_{0}-1}^{\io}\EE_n$,
and there it gives a well defined function.

\*

The check that the functions $\hh(\pps)$ defined by the convergent
renormalized series evaluated at $\pps=\oo t$ do actually solve the
equations of motion can be performed by repeating the corresponding
analysis in Ref. \cita{Ge}.  The equation that
$\hhh=(\aaa,\bbb)$ has to solve is $\hhh =
\e g (\dpr_{\aa} f(\pps + \aaa, \bb_{0} + \bbb),
\dpr_{\bb} f(\pps + \aaa, \bb_{0} + \bbb))$ where $g$
is the pseudo-differential operator $(\oo\cdot\nn)^{-2}$.  The proof is of
algebraic nature and ultimately follows from the fact that the series
we are considering is a resummation of Lindstedt's series which is
a formal solution of the problem.  This explains why the various
algebraic identities necessary for the check actually hold and the
proof proceeds exactly as in Section 8 of Ref. \cita{Ge}: we reproduce
the argument and the chain of identities in Appendix \secc(A5).
Therefore the proof of Theorem 1 in Section \secc(1) is complete.

\*\*
\section(7,Concluding remarks)

\0The analysis can be immediately extended to the case in which the
matrix $\dpr^2_{\bb}f_{\V0}(\bb_{\V0})$ has some non-degenerate
positive eigenvalues and some additional negative ones. The negative
eigenvalues give no problems and they can be treated as in the case of
Ref. \cita{GG} in
which all eigenvalues are negative. The negative
eigenvalues do not give rise to new small divisors, unlike the
positive ones; in more physical language the proper time scales (\ie
real proper frequencies) of the tori cannot resonate with the time scales
of hyperbolic type (\ie imaginary) introduced by the perturbation. Hence
the following generalization of Theorem 1 holds.

\*

\0{\bf Theorem 2.} {\it If the matrix $\dpr^2_{\bb}f_{\V0}(\bb_0)$ is not
singular and has pairwise distinct eigenvalues the
conclusions (i), (ii) and (iii) of Theorem 1 in Section \secc(1) follow
also in this case.}

\*

The present work has developed a combinatorial approach to the proof
that the frequencies of elliptic type possibly introduced by
the perturbation do not resonate with the proper frequencies of the tori
at least if $\e$ is not too special in a small interval $[0,\lis\e]$.
\ie if it is in a set $\EE\subset [0,\lis\e]$ of large measure near $0$:
Nevertheless the complement of $\EE$ is an open dense set in
$[0,\lis\e]$. The results hold for the Hamiltonian \equ(1.2) and the
special resonances $(\oo,\V0)$ considered: they can be extended to the
most general resonances of Hamiltonians like \equ(1.2) with a general
quadratic form for the kinetic part (\ie with $\fra12 {\bf I}\cdot{\bf
I}$ replaced by $\fra12 {\bf I}\cdot Q{\bf I}$ with $Q$ a non-degenerate
$d\times d$ matrix).

\*

The case of $\dpr^2_{\bb}f_{\V0}(\bb_{\V0})$ with degenerate eigenvalues 
seems quite different from the one treated here.  Degeneracy will be
removed to order $O(\e^2)$ under generic conditions. However $O(\e^2)$
is also the order of variation of the self-energies and one has to
find a way to perform the resummations even between scale $\lis n_{0}$
and scale $2 \lis n_0$, which is the scale at which the singularities
of the propagator are split apart and one shall be able to proceed in
the same way as we did in the case of non-degenerate eigenvalues.
\*

The Lipschitz regularity in $\e$ in Theorems 1 and 2 can be replaced by
$C^{k}$ regularity for any $k$ by exploiting the comments in Remark
(2) to Lemma 2 and Remark (2) in Appendix \secc(A3.2).

Unfortunately there seems to be no example known in which one can
check that the power series studied here are divergent {\it as power series}
in $\e$. Note that the ({\it infinitely many}) 
divergent series that have arisen in this paper
are obtained by first splitting the coefficient of order $k$ in the
Lindstedt power series and then collecting contributions from the
different orders in $\e$: the latter form divergent series for which we have
assigned a summation rule. Therefore we have not proved divergence of
the Lindstedt series as power series in $\e$: in this sense
({\it unlikely}) convergence of the Lindstedt series has not been ruled out
({\it yet}). Nor there is any uniqueness result on the value of the
renormalized series. The latter depends on quite a few arbitrary
choices (even in the hyperbolic cases): for instance the cut-off
shapes in Fig. 2 are quite arbitrary and in principle the allowed
$\e$'s will change with the choice.

Furthermore, although we have not really checked all necessary
details, it seems to us that our method also shows that, given a
value $\e_0$ for which the renormalized series converges, one can find
a complex domain of $\e$ which is open, reaches the real axis
with a vertical cusp at $\e_0$ and extends to an
open region including a segment $(-\h,0)$ on the negative real
axis. In this domain the renormalized series should converge 
taking on the real axis  real values parameterizing an {\it hyperbolic}
torus with the same rotation vector. However since there are no
uniqueness proofs we cannot guarantee that each such extension does not
correspond to a {\it different} torus (close within any power of $\e$
to any other torus of the same type). This would signal a ``giant
bifurcation'' that one would like to exclude: in Ref. \cita{GG} an attempt
was made to show uniqueness by estimating the size of the Lindstedt
series coefficients aiming at applying the theory of Borel transforms.
However we could not prove good enough bounds. We obtained $k!^\a$ growth with
a too large $\a$ (given our estimated size of the domain of analyticity
in $\e$) to apply uniqueness results from the theory of Borel
summations .

\*\*
\appendix(A1, A brief review of earlier results)

\0The system which is usually studied in literature when
the problem of persistence of lower-dimensional elliptic tori
is studied, is of the form
$$ \HHH = \oo(\x) \cdot\AA + \sum_{k=1}^{s} \O_{k}(\x)
\left( q_{k}^{2} + p_{k}^{2} \right) + P(\aa,\AA,\qq,\pp) ,
\Eqa(A1.1) $$
where $(\aa,\AA,\pp,\qq)\in\TTT^{r}\times\RRR^{r}\times\RRR^{s}\times
\RRR^{s}$. The function $P$ is analytic in its arguments, and $\x$ is a
parameter in $\RRR^{r}$; the function $P$ is a {\it perturbation}:
this means that a rescaling of the actions could allow us to introduce
a small parameter $\e$ in front of the function $P$. The frequencies
of the harmonic oscillators are called {\it normal frequencies}; the
case $\O_{k}(\x)=\O_{k}=$ constant (that is with the normal frequencies
independent of $\x$) is a particular case, and it is usually
referred to as the ``constant frequency case''. Existence of
invariant tori for the system \equ(A1.1) was originally proved by
Mel'nikov \cita{Me1}, \cita{Me2}, new proofs were produced
by Eliasson \cita{E1}, Kuksin \cita{Ku}, and P\"oschel \cita{P1}.
The case $s=1$ is easier, and it was earlier solved by Moser \cita{Mo}.
Later proofs were given by R\"ussmann, see for instance Ref. \cita{R}.
See also the very recent Ref. \cita{LW}.

For $P=0$ the dimension of the tori is $r<d$ and the variables
$(\qq,\pp)$ move around stable equilibrium points, hence such tori are
called {\it elliptic lower-dimensional tori}.

The conditions under which the quoted results are proved are, besides
the usual Diophantine condition \equ(1.3)
on $\oo$, two non-resonance conditions involving one and two normal
frequencies (the so called first and second Mel'nikov conditions,
originally introduced in Ref. \cita{Me1}); in particular one has to
impose that the normal frequencies are non-degenerate (i.e. they have
to be all different from each other).

Recently proofs of existence of elliptic lower-dimensional tori
were given by requesting only the first Mel'nikov conditions:
this allows treating degenerate frequencies.
The first result in this direction
is due to Bourgain \cita{B3}, where the ideas
introduced in Refs. \cita{CrW} and \cita{B1} to prove
existence of periodic and quasi-periodic solutions
in nearly integrable Hamiltonian partial differential equations
were adapted to construct lower-dimensional tori
in the finite-dimensional Hamiltonian systems \equ(A1.1)
corresponding to the case of constant normal frequencies.
New proofs, extending the results also to the case
of non-constant normal frequencies, are due
to Xu and You \cita{Y}, \cita{XY}.

An extension of the results of existence of periodic and quasi-periodic
solutions describing lower-dimensional invariant tori
for infinite-dimensional PDE systems has been provided in
a series of papers, which include Refs. \cita{Ku}, \cita{Wa}, \cita{CrW},
\cita{KP}, \cita{P2}, \cita{B1}, \cita{B2}, \cita{B4}, \cita{GM2}
and \cita{GMP}.

\*

On the other hand the problem \equ(1.2) has not been widely studied
in literature. It corresponds to a degenerate
case because in absence of perturbations the lower-dimensional tori
are neither elliptic nor hyperbolic: it is the perturbation
itself which determines if the tori, when continuing to exist,
become elliptic or hyperbolic (or of mixed type or parabolic).

\0(i) The case of hyperbolic tori is easier, and it was
the first to be studied, by Treshch\"ev \cita{T}.
Recently the problem was reconsidered in Ref. \cita{GG},
where the analyticity domain of the invariant tori
was studied in more detail. In the case of elliptic tori
the problem was considered in Refs. \cita{ChW} and \cita{WC},
where Treshch\"ev's approach to the study of the case
of hyperbolic tori, involving a preliminary change of coordinates,
is used to cast the Hamiltonian in a form
which is suitable for applying P\"oschel's results on elliptic tori:
in particular this imposes the same conditions as in Ref. \cita{P1}
on the normal frequencies which appear after the canonical change
of coordinates is performed.

\0(ii) The existence problem has been also considered in Ref. \cita{JLZ},
{\it where elliptic and hyperbolic tori were studied simultaneously},
again by imposing some non-degeneracy conditions on normal
frequencies. Ref. \cita{JLZ} does not investigate resummations of
Lindstedt's series; it is based on a rapid convergence method, close
in spirit to the original proofs of the KAM theorem: a concise
existence proof of lower-dimensional tori is achieved in both the
elliptic and hyperbolic cases.

\*

We stress that in all quoted papers,
except Ref. \cita{JLZ} and \cita{T}, the problem is considered
with $\e$ (i.e. the size of the perturbation) fixed
and the study deals with estimates of the measure of the
rotation vectors $\oo$ for which there exist invariant tori.
We suppose, instead, that $\oo$ is fixed,
hence  we study the dependence on $\e$ of the lower-dimensional
invariant tori and, in particular, the set of values of $\e$
for which the tori survive.
\*

Our techniques extend those in Refs. \cita{GG} and \cita{Ge}, and are
based on the method introduced in Refs. \cita{E2} and \cita{Ga}.
With respect to Ref. \cita{Ge}, where existence of quasi-periodic
solutions is proved for the {\it generalized Riccati equation}
considered in Ref. \cita{Ba}, the main difficulty is due to the
presence of several normal frequencies.  It is not surprising that this
generates extra technical difficulties: as already noted, it is well
known that the case $s=1$ is easier; see Refs. \cita{Mo} and \cita{C}.
An advantage of the present method is that it is fully constructive and
gives a very detailed knowledge of the solution.

\*\*
\appendix(A2,Excluded values of $\e$)

\0Define

\kern-5mm
$$\eqalign{
& \r_{\lis n_0-1}\defi
\Sqrt{\fra{\e}{a_{s}}} \min \Big\{
\min_{i>r} |\dpr_\e \Sqrt{\ul{\l}^{[0]}_{i}(\e)}|,
\min_{i \neq j \atop i,j>r}
|\dpr_\e \Sqrt{\ul{\l}^{[0]}_{j}(\e)} -
\dpr_\e \Sqrt{\ul{\l}^{[0]}_{i}(\e)} | \Big\} , \cr
& \r'_{\lis n_0-1} \defi \fra{1}{\Sqrt{\e a_{s}}} \max_{j} \Big\{
\Sqrt{\ul{\l}^{[0]}_{j}(\e)} \Big\} , \cr}
\Eqa(A2.1) $$
and note that $\r_{\lis n_0-1}$ is bounded from below
proportionally to $\r$, as defined in \equ(4.2), and $\r'_{\lis n_0-1}=1$.
Then \equ(3.3) excludes, for each $\nn$, an interval in $\e$ whose
measure is bounded (using $\Sqrt{a_s \e} \le C$; see \equ(3.2)) by
$$ 2^{-(\lis n_0-1)/2} C \,C_{0}\, K_0 |\nn|^{-\t_{1}} ,
\Eqa(A2.2)$$
where the constant $K_0$ can be estimated
by $K_0=s\,a_{s}^{-1} \r^{-1}_{\lis n_0-1}$.

The Diophantine condition on $\oo$ implies that if \equ(3.3) is invalid
then $|\nn|$ cannot be too small
$$ 2 \Sqrt{\e a_{s} \r'_{\lis n_0-1}} + 2^{-(\lis n_0-1)/2}
C_{0}|\nn|^{-\t_{1}} \ge |x| \ge C_{0}|\nn|^{-\t_{0}} .
\Eqa(A2.3)$$
Therefore $\Sqrt{\e a_{s} \r'_{\lis n_0-1}} \ge
\fra14 C_{0}|\nn|^{-\t_{0}}$ if $\lis n_0\ge3$, hence in this case 
we only have to consider
the values of $\nn$ with $|\nn|\ge (C_{0}/(4 \Sqrt{\e a_{s}
\r'_{\lis n_0-1}}))^{1/\t_{0}}$. Since $C/2<\Sqrt{\e a_s}\le C =
2^{-n_0}C_{0}$, we get the bound \equ(3.5) with
$\t_{1}=\t+r+1$ and $K= K_0\,C_{0}\,(4 C \Sqrt{\r'_{\lis n_0-1}}\,
C_{0}^{-1})^{(\t_{1}-r-1)/\t_{0}} \sum_{\nn\ne\V0}
\fra1{|\nn|^{r+1}}= 4 K_0 \Sqrt{{\r'}_{\lis n_0-1}}
\sum_{\nn\ne\V0} \fra1{|\nn|^{r+1}}$. Note that a condition
like $\t_{1}>\t+r$ is sufficient to obtain both summability over $\nn$
and a measure (of the excluded set) relatively small
with respect to that of $I_{C}$. If $\lis n_0<3$, hence $n_0<3$,
the same conclusion trivially holds possibly increasing
the value of $K$ by a factor $4$.

\*\*
\appendix(A3,Resummations: convergence and smoothness)

\0To prove Lemma 2, we first show that the series
defining $M^{[n]}(x;\e)$ for $0\le n \le \lis n_0$ converge and then
we check smoothness and the bounds.  This is done for completeness
as the argument is almost a word by word repetition of
the analysis in Ref. \cita{GG}, with a few slight changes of notations
necessary to adapt it to our present notations and scope.
To study convergence of the series defining $M^{[n]}(x,\e)$,
$n\le\lis n_0$, we remark that we have to consider only trees
in which all propagators have scales $[p]$ with $p\le \lis n_0$.
Therefore the propagators which do not vanish will be such that
their denominators satisfy $D(x) > 2^{-2(\lis n+1)} |x|^{2}$,
see \equ(4.4), so that they are effectively estimated from
below by $|x|^{2}$ times a constant. Note that the case $n=0$ is obvious
(and it is treated in Section \secc(3)).

\*

\asub(A3.1) {\it Convergence.}
We suppose that the eigenvalues of $\MM^{[\le p]}(x;\e)$,
$n=0,\ldots,n-1$, differ from the corresponding ones of
$\MM^{[\le 0]}(x;\e)\=M_{0}$ so that
$|\l^{[p]}_j(x,\e)-\l^{[0]}_j|<\g \e^2$ for some $\g>0$, and that $\e$
is small enough so that $\g\e^2<\fra12\e a_s 2^{-2\lis n-2}$ and,
therefore (see \equ(4.4)), if a line with frequency $x$ has scale
$[p]$, $p<n$, then $|x^2-\l^{[p]}_j(x,\e)|>2^{-2(\lis n+2)} x^2$.

We shall use that if a the propagator of a line is on a scale $[n]$
then one has $D(x) \le 2^{-2(n-2)}C^2_{0}$, even though {\it we could
allow} also a bound $D(x) \le 2^{-2(n-1)}C^2_{0}$.  The reason for this
is again for later use in bounds necessary to establish
the needed cancellations as commented in Section \secc(A3.2).

Consider a renormalized self-energy cluster $T \in \SS^{\RR}_{k,n-1}$, and
define $\L_{m}(T)=\{\ell\in\L(T)\,:\,n_{\ell}=m\}$, for $m \le n-1$, and
${\cal P}(T)$ the set of lines (path) connecting the external lines of $T$.

If $\nn$ is the momentum flowing in the line entering $T$ then the
momentum flowing in a line $\ell\in\L(T)$ of scale $[p]$, $p\le n-1$,
will be $\nn^{0}_{\ell} + \s_{\ell} \nn$, $\s_{\ell}=0,1$,
where $\nn_{\ell}^{0}$ is the momentum that would flow on $\ell$
if $\nn=\V0$. The corresponding frequency will be
$x'_{\ell}=x^{0}_{\ell} + \s_{\ell} x$, with obvious notations.

First of all we shall prove the bound
$$ \sum_{\vvvv \in V(T)} |\nn_{\vvvv}| \ge 2^{(n-\lis n-5)/\t_{0}} .
\Eqa(A3.1) $$
for $T\in \SS^{\RR}_{k,n-1}$. If there is a line $\ell\in\L_{n-1}(T)$
which does not belong to ${\cal P}(T)$ then $x_{\ell}=x^0_{\ell}$,
so that \equ(A3.1) follows from the Diophantine condition on $\oo$.
If all lines in $\L_{n-1}(T)$ belong to ${\cal P}(T)$ consider
the one among them, say $\ell$, which is closest to $\ell_{T}^{2}$,
\ie the entering line of $T$. Then call $T_{1}$
the connected set of nodes and lines between\footnote{${}^7$}{\nota
The lines between two lines $\ell_1$ and $\ell_2$ with $\ell_2<\ell_1$
are all the lines which precede $\ell_1$ but which do not
precede $\ell_2$ nor coincide with it.\vfil}
$\ell$ and $\ell_{T}^{2}$. If $T_{1}$ is a single node $\vvv$
then $\nn_{\vvvv}\neq \V0$, otherwise $\vvv$ would be
a trivial node; if $T_{1}$ is not a single node then by construction
all the lines of $T_{1}$ have scales strictly smaller than $n$,
hence $x_{\ell} \neq x$ otherwise $T_{1}$ would
be a self-energy cluster. In both cases one has
$|x_{\ell}-x|=|x_{\ell}^{0}| > C_{0} |\sum_{\vvvv\in V(T_{1})}
\nn_{\vvvv}|^{-\t_{0}}$. On the other hand both $D(x)$ and $D(x_{\ell})$
must be $\le (C_{0} 2^{-(n-2)+1})^2$ hence, by \equ(4.4)
$|x|,|x_\ell|\le C_0 2^{-n+\lis n+3}$, so that
$|x-x_{\ell}| \le C_{0} 2^{-n+\lis n + 4}$, and
\equ(A3.1) follows also in such a case.

The next task will be to show that the number $\NN_{m}(T)$ of lines
on scale $[m]$, with $m\le n-1$, contained in a cluster $T$ is bounded
by $\NN_{m}(T)\le \max \{ E_{m} \sum_{\vvvv\in V(T)}|\nn_{\vvvv}|-1,0\}$,
with $E_{m}= E\,2^{-m/\t_{0}}$ for a suitably chosen constant $E$; as
it will emerge from the proof one can take $E=2\,2^{(\lis n+4)/\t_{0}}$.

Before considering clusters we adapt to our context the
classical bound (Siegel-Bryuno-P\"oshel; see Ref. \cita{Ga}
and references quoted therein), stating that,
if $\NN_{m}(\th)$ denotes the number of
lines on scales $[m]$, then by induction
on the number of nodes of $\th$ one shows: $\NN_{m}(\th) \le
\max \{ E_{m} \sum_{\vvvv\in V(\th)}|\nn_{\vvvv}|-1,0\}$.
Indeed if $\th$ contains only one node $\vvv_0$ and the frequency
$x=\oo\cdot\nn_{\vvvv_{0}}$ of the root line has scale $[m]$ one has
$$ 2^{-m+1} C_0 \ge \Sqrt{D(x)} \ge 2^{-(\lis n+1)} |x|
\ge 2^{-(\lis n+1)} C_0 |\nn_{\vvvv_0}|^{-\t_0}\,\tto
|\nn_{\vvvv_{0}}|> 2^{(m-\lis n-2)/\t_0} ,
\Eqa(A3.2)$$
hence $E_{m}|\nn_{\vvvv_0}|-1\ge 2$ and the bound holds in this simple
case.

If $\th$ has $k$ nodes and the root line {\it does not} have scale
$[m]$ the inductive assumption, if it is assumed for the cases of
$k'<k$ nodes, gives the bound for $k$-nodes trees.

If the root line has scale $[m]$ then on each path of tree lines
leading to the root we select the line among the ones on scales $[m']$
with $m' \ge m$ {\it closest to the root} (if any is found on the
path) and we call the selected lines $\ell_{1},\ldots,\ell_{q}$. If
$q\neq 1$ either the bound follows just as in the case of $k=1$ (when
$q=0$) or from the inductive hypothesis (when $q\ge 2$).

The case $q=1$ and $[n_{\ell_{1}}]=[\io]$ (\ie $\nn_{\ell_{1}}=\V0$)
can be treated as the case $q=0$. If $q=1$ and
$\nn_{\ell_{1}}\neq \V0$, by construction all lines between
the root line $\ell$ and $\ell_{1}$, see footnote ${}^{7}$,
have scales $[m']$, with $m'<m$, so that such lines,
together with the nodes they connect, form a cluster $T$.
The frequencies $x_{\ell}$ and $x_{\ell_{1}}$ must be different
because the tree $\th$ contains no self-energy clusters.
On the other hand $\Sqrt{D(x_{\ell})},\Sqrt{D(x_{\ell_{1}})}
\le 2^{-m+1}C_{0}$, hence $|x_{\ell}|,|x_{\ell_{1}}| \le 
2^{-m + \lis n+3} C_{0}$ by \equ(4.4), and $C_0 |\nn_{\ell}-
\nn_{\ell_{1}}|^{-\t_0} \le |x_{\ell}-x_{\ell_{1}} |\le
2^{-m+\lis n+4} C_0$, so that we get $\sum_{\vvvv \in V(T)} |\nn_{\vvvv}|
\ge (2^{-m+\lis n+4})^{-1/\t_{0}}$, which gives
$\NN_m(\th)\le 1 + E_{m} \sum_{\vvvv\in V(\th)}|\nn_{\vvvv}|
-E_{m} \sum_{\vvvv\in V(T)}|\nn_{\vvvv}| -1 \le
E_{m} \sum_{\vvvv\in V(\th)}|\nn_{\vvvv}| -1$,
so that the bounds is completely proved.

\*

\0{\it Remark.} The above discussion  exploits the
property that the tree $\th$ that we consider {\it cannot}, by
definition of renormalized tree, contain self-energy clusters, and
follows Ref. \cita{GG} which was based on the possibility of bounding
the denominators proportionally to $x^2$ (in that case the
proportionality factor was $1$): a property also valid here for
$n\le \lis n_0$.

\*

For the bound on $\NN_{m}(T)$  we
consider a subset $G_{0}$ of the lines of a tree $\th$
between two lines $\ell_{\rm out}$ and $\ell_{\rm in}$.
Set $G=G_{0}\cup \ell_{\rm in}\cup\ell_{\rm out}$.
Let $[p_{\rm in}],[p_{\rm out}]$ be the scales of the lines
$\ell_{\rm out}$ and $\ell_{\rm in}$, respectively,
and suppose that $p_{\rm in},p_{\rm out}\ge m$, while all lines in
$G_{0}$ (if any) have scales $[p]$ with $p \le n-1$. Note that in general
$G_0$ is not even a cluster unless $p_{\rm in},p_{\rm out}\ge n$.
Then we can prove that $\NN_{m}(G_{0}) \le \max\{ E_{m}
\sum_{\vvvv\in V(G_0)}|\nn_{\vvvv}|-1,0\}$, where $V(G_{0})$ is the
set of nodes preceding $\ell_{\rm out}$ and following $\ell_{\rm in}$,
and $E_{m}$ is defined above.
If $G_{0}$ has zero lines then the harmonic $\nn_{0}$ of the (only)
node in $V(G_{0})$ is large, $|\nn_{0}|\ge 2^{(m-\lis n-2)/\t_{0}}$
(by the Diophantine property)
and the statement is true. Hence we proceed inductively
on the number of lines in $G_{0}$.

If no line of $G_{0}$ on the path ${\cal P}(G)$ connecting the
external lines of $G$ has scale $[m]$ then the lines in $G_{0}$
on scale $[m]$ (if any) belong to trees with root on
${\cal P}(G)$, and the statement follows from the bound on trees.

Suppose that $\ell\in {\cal P}(G)$ is a line on scale $[m]$,
then call $G_{1}$ and $G_{2}$ the disjoint subsets of $G$ such that
$G_{1}\cup G_{2} \cup \ell = G$. Then $G_1\cup\ell$ and $G_2\cup\ell$
have the same structure of $G$ itself but each has less lines: and
again the inductive assumption yields the result.

Therefore, as a particular case, by choosing $G_0=T$, with
$T\in\SS^{\RR}_{k,n-1}$, the bound for $\NN_{m}(G)$ implies the bound
on $\NN_{m}(T)$ we are looking for. 

The above analysis is taken from Ref. \cita{Ge} and differs
from Ref. \cita{GG} because here the scales depend on $\e$ and it is
not clear how to define a ``strong Diophantine condition'', which
would
allow a one-to-one correspondence between line scales and line momenta.

The bound on the contribution of a single self-energy cluster
$T\in \SS^{\RR}_{k,n-1}$ is then
$$ \eqalignno{
& \fra{\e^{k}}{k!} C_{0}^{-2k}
F^{2k} e^{-\fra12\k_0\sum_{\vvvv}|\nn_{\vvvv}|}
\Big(\prod_{m=0}^{m_0} 2^{(m+3)2\NN_{m}(T)} \Big)\cdot
& \eqa(A3.3) \cr
&\cdot\Big( e^{- \fra{1}{2}\k_{0} \sum_{\vvvvv\in V(T)} |\nn_{\vvvv}|}
\prod_{m=m_{0}+1}^{\io} 2^{2(m+3) 2^{-m/\t_{0}} E
\sum_{\vvvvv\in V(T)} |\nn_{\vvvv}|} \Big)\ \le\ \fra{\e^{\k}
G^k}{k!}  e^{-\fra12
\k_0\sum_{\vvvvv\in V(T)} |\nn_{\vvvv}|} , \cr}$$
with $F$ an upper bound on the constants $F_0,F_1$ bounding the
Fourier transform of the perturbation (see \equ(1.4)), while $m_0$ is
defined so that $\log 2 \sum_{m>m_0} 2(m+3)2^{-m/\t_{0}} E
\le \fra12\k_{0}$ and $G$ is a suitable constant.

The number of trees can be bounded by $4^k k!$, and the sum over
the scale labels involves at most $2$ possible values per line because
of the upper and lower cut-offs present in the propagators definition.
The sum over the harmonics can be estimated by making use of {\it
part} of the exponential factor in \equ(A3.3) (say $\fra14\k_0$) while
the other $\fra14\k_0$ will be used as a factor bounded by
$e^{-\fra14 \k_0 2^{(n-\lis n-5)/\t_{0}}}$, by \equ(A3.1).

Hence we get convergence at exponential rate $2^{-1}$ for $\e<\e_{1}$
(and $\e_{1}$ is an explicitly computable constant) and the matrix
$M^{[n]}(x;\e)$ is defined by a convergent series and it is bounded by
$$ \| M^{[n]}(x;\e) \| < \lis B \e^{2}
e^{-\fra14 \k_0 2^{(n-\lis n-4)/\t_{0}}} ,
\Eqa(A3.4) $$
for a suitable $\lis B$ which can be read from \equ(A3.3),
\ie we get the first of the first line in
\equ(5.13) with the constant $B$ replaced by $\lis B$, $\t=\t_{0}$,
and $\k_1=\fra14\k_0e^{-(\lis n +4)/\t_{0}}$.
The $\e^2$ factor is due to the parallel
remark that, in any self-energy cluster whose value
contributes to $\MM^{[n]}(x;\e)$, $k$ is certainly $\ge
2$ (see Remark to Definition 4 in Section \secc(5)).

Therefore if $\e$ is small enough (that is smaller than a
constant independent of $n\le \lis n_0$)
$$ \|\MM^{[\le n]}(x;\e)- M_0\|\le \lis B \e^{2} \sum_{n=1}^\io
e^{-\fra14 \k_0 2^{(n-\lis n-4)/\t_{0}}} \defi B'\e^2 ,
\Eqa(A3.5)$$
so that the eigenvalues of $\MM^{[\le n]}(x;\e)$ will be shifted
with respect to the corresponding eigenvalues of $M_0$ by $\g \e^2$
at most, with $\g\defi B' C$, see (I) in Appendix \secc(A4).

Hence if we define $\g$ as $B'C$ and $\e$ is chosen small enough, say
$\e<\e_2$, so that $\g\e^2<\fra12\e a_s 2^{-2\lis n-2}$ (as it must be
in order that the above argument be consistent, see the beginning of
the current Section) we obtain the validity of the assumed inductive
hypothesis for all $n\le\lis n_0$ and of the first inequality in the
first line of \equ(5.13) where $B$ can be chosen equal to $\lis B$ above.

\*

The symmetries in items (i) and (ii) are an algebraic consequence
of the form of the Lindstedt series: hence they are a necessary
consequence of the proved convergence, see Ref. \cita{GG}.

\*

\asub(A3.2) {\it Smoothness.}
The function $M^{[n]}(x;\e)$ which we have just shown to be well defined
for all $\e$ small enough will be smooth in $\e,x$.
We assume inductively that this is the case for $M^{[p]}(x;\e)$,
$0\le p<n-1$, and that the bounds in the first line of \equ(5.13) hold for
such $p$'s (the case $p=0$ is obvious as $\MM^{[0]}(x;\e)\=M_0$).

Each derivative with respect to $x$ or, respectively, to $\e$ will
replace the value of a self-energy cluster with $k$ nodes by a sum
of $k$ terms which can be bounded by a bound like \equ(A3.3).

In fact, given a self-energy cluster $T$, the right derivative
$\dpr^+_x$ may fall on a denominator of one of the $k-1$ cluster lines.
If its frequency is $x+x_0$ with scale label $[m]$,
derivation yields, up to a sign, a product of two matrices 
$((x_0+x)^2-\MM^{[\le m]}(x_{0}+x;\e))^{-1}$
times $2\,(x_0+x)-\dpr_{x}^{\pm}\MM^{[\le m]}(x_{0}+x;\e)$ with an
appropriate order of multiplication.
The  term $2\,(x+x_{0})\,((x_0+x)^2-\MM^{[\le m]}(x_{0}+x;\e))^{-2}$
can be bounded proportionally to $(C_0^{-2}2^{2(m-1)})^{3/2}
\le (C_0^{-2}2^{2(m-1)})^{2}$, while the remaining term can be
studied by making use of the inductive assumption
$\|\dpr_x\MM^{[\le m]}(x_{0}+x;\e) \|\le B\e^2 a_s^{-1/2}$
and it leads to the same bound found for the first term,
\ie $(C_0^{-2}2^{2(m-1)})^{2}$, {\it multiplied}
by $B\e^2$.\footnote{${}^8$}{\nota Since the matrix $M^{[m]}(x_{0}+x;\e)$
is generated by self-energy clusters clusters of degree at least $2$.}

If the derivative falls on either a $\psi_{p}$ or
a $\chi_{p}$ function, we can use that such derivative can be bounded
proportionally to $C_0^{-1}2^{p}$ and $\sum_{p=0}^{m-1}2^{p} = 2^{m}$,
to obtain again the same bound as the first case. 

Hence the final bound has the form $B_1+\e^2 B b$ with $B_1,b$
suitable constants, provided $\e$ is small enough, say $\e<\e_3$.
The value of the constants $B_1,b$ do not depend on the
inductively assumed value for $B$: in particular
$B_{1}$ can be obtained (see Remark (2) below for a smarter
bound) by replacing $2^{(m+3)}$ in the two factors in the \lhs
of \equ(A3.3) by $2^{2(m+3)}$ and by inserting a factor $k$
times a constant (to keep track of all the constant factors
arising from differentiation). Therefore if $B=2 B_1$
the estimate on $\dpr^{+}_x \MM^{[\le n]}(x;\e)$ follows if $\e$
is small enough, say $\e<\e_4$.
The same can be said about the left derivative $\dpr^-_x$.

The right and left differentiability of $\MM^{[n]}(x;\e)$ with respect
to $x$ is due to the dependence of $\MM^{[n]}(x;\e)$ on the function
$D(x)$: the latter has a discontinuous derivative at a finite number
of points (roughly at midpoints between the eigenvalues $\l^{[0]}_{j}$
of $M_0$).\footnote{${}^{9}$}{\nota One could avoid having only left and
right differentiability by using a regularized version of the function
$D(x)$ as discussed in Remark (2) after Lemma 2 in Section
\secc(5).\vfil}
Note that the denominators in the self-energy values defining
$M^{[n]}(x;\e)$ cannot vanish, and actually stay well away from $0$,
permitting the above bounds, because of the lower cut-off
$\psi_0(D(x))$ appearing in the definition of the propagators
$g^{[0]}(x;\e)(x;\e)$; see \equ(5.6) and \equ(5.7).

The same argument holds for $\dpr_\e$: however the bound will be only
$B\e$ instead of $B\e^2$ because the derivative with respect to $\e$
might decrease by one unit the degree of the self-energy values
involved. Thus the first line of \equ(5.13) is completely proved.
Of course for each of the three terms we get a different constant $B$,
but for simplicity we use for them all the largest, still calling it $B$.

\*

\0{\it Remark.} (1) We could also prove existence of higher
$x,\e$-derivatives of $\MM^{[n]}(x;\e)$ and of its eigenvalues
$\l^{[n]}_j(x,\e)$ for $j>r$ via the above argument.
\\
(2) {\it The more derivatives we try to estimate with the
above method the smaller would become the set of allowed values of
$\e$}. This however is avoidable. Instead of imagining to include the
bound $C_0^{-2} 2^{2m}$ arising above as a consequence of the ``extra''
$D(x+x_0)$ or of the other derivatives into the factors
$2^{m+3}$ associated with the divisors in \equ(A3.3) one could simply
further bound this by $C_0^{-2} 2^{2n}$ and use part of the factor
$e^{-\fra14 \k_0 2^{(n-\lis n-4)/\t_{0}}}$ (say replacing $\fra14\k_0$
with $\fra18\k_0$): this eventually leads to a bound on the $s$-th
right-derivative with respect to $x$ of a value of a self-energy
cluster proportional to $2^{ns}e^{-\fra18 \k_0 2^{(n-\lis
n-4)/\t_{0}}}$ but with an $s$-independent estimate of the radius of
convergence (as the constant $G$ in \equ(A3.3) remains the same). This
is sufficient to get the existence of the $s$-th derivatives without
any further restriction on $\e$: and a similar argument holds for the
$\e$-derivatives.

\*

\asub(A3.3) {\it Cancellations.}
Only the bound in the fourth line of \equ(5.13) follows from
those in the first line. The bounds in the second and third lines
express remarkable properties of Lindstedt series and are essentially
algebraic properties: they are the ``same'' cancellations which occur
in KAM theory, see Refs. \cita{Ga}, \cita{GM1}, and are based on the
remark that if $T$ is a self-energy cluster the entering and exiting
lines have the same momentum $\nn$: hence the sum of the harmonics of
the nodes of $T$ vanishes $\sum_{\vvvv\in V(T)}\nn_{\vvvv}=\V0$.

We start by dealing with the trivial cases.
Consider first self-energy clusters $T$ such that 
$$ \sum_{\vvvv\in V(T)}|\nn_{\vvvv}|\ge
(C_{0}/2^{6}|x|)^{1/\t_{0}} .
\Eqa(A3.7)$$
For such a self-energy cluster $T$ one can use part
(say $1/8$) of the exponential decay of the node factors to obtain
a bound $e^{-\fra{\k_{0}}{8}\sum_{\vvvvv\in V(T)}|\nn_{\vvvvv}|}
\le e^{-b_{1} |x|^{-1/2\t_{1}}} \le b_{2} x^{2}$,
with $b_{1}$ and $b_{2}$ two suitable positive constants,
while a factor $\e^{2}$ simply follows from the fact that
any self-energy cluster has at least two nodes. 

So we can assume that \equ(A3.7) {\it does not hold}.
If $\nn$ is the momentum flowing in the entering line then the
momentum flowing in a line $\ell\in \L(T)$ of scale $[p]$, $p\le n$, if the
scale of the cluster is $[n]$, will be $\nn^{0}_{\ell} + \s_{\ell} \nn$,
$\s_{\ell}=0,1$, where $\nn_{\ell}^{0}$ is the momentum that would
flow on $\ell$ if $\nn=\V0$. The corresponding frequency will be
$x_{\ell}=x^{0}_{\ell} + \s_{\ell} x$, with obvious notations.

If the entering and exiting lines
are imagined attached to the internal nodes of $T$ in all possible
ways (\ie in $k^2$ ways if $T$ contains $k$ nodes) {\it keeping all
their labels unaltered} then one obtains a
family $\FF_T$ of self-energy clusters. The
contribution of each self-energy cluster of $\FF_T$ to each of the
entries of the matrix $\MM^{[n]}(x;\e)$ with labels $i,j\le r$ (the $\a\a$
entries in the notations of Lemma 2) and with labels $i\le r, j>r$ (the
$\a\b$ entries) has the form $
M_{i,j;\vvvv,\wwww}(x,T)\, \n_{\vvvv,i}\,\n_{\wwww,j}$ or, respectively,
$ M'_{i,j;\vvvv}(x,T)\,\n_{\vvvv,i}$, with
$$\eqalign{
&M_{i,j; \vvvv,\wwww}(x,T)=M_{i,j}(T)+ x
M^{(1)}_{i,j}(T)+x^2M^{(2)}_{i,j,\vvvv,\wwww}(x,T),\qquad i,j\le r ,
\cr
&M'_{i,j; \vvvv}(x,T)=M'_{i,j}(T)+ x
M^{\prime(1)}_{i,j,\vvvv}(x,T),\qquad \kern2.1cm i\le r,\, r<j ,\cr}
\Eqa(A3.6)$$
so that after performing the sum over the self-energy clusters
of $\FF_T$, \ie after performing the sums $\sum_{\vvvv,\wwww\in V(T)}$
or, respectively, $\sum_{\vvvv\in V(T)}$, the first two terms
in the first line and the first term in the second
line do not contribute because $\sum_{\vvvv}\nn_{\vvvv}=\V0$.  However
one has to show that the matrices $M$ and $M'$ in the {\it r.h.s.} of
\equ(A3.6) satisfy appropriate bounds once the factors $x$ determining
the order of zero at $x=0$ are extracted. From the convergence one
expects that the bounds should still be proportional to $\e^2 $ while
the derivatives $\dpr_x^{\pm}$ or $\dpr_\e$ should satisfy bounds
proportional to $\e^2$ or to $\e$ respectively.

The \equ(A3.6) are proved by means of interpolations, see \cita{GM1},
between the contributions of the self-energy clusters in the family
$\FF_T$. When we collect together the values of the self-energy
clusters in $\FF_T$ then the arguments of some of the propagators can
fall outside the supports of the respective cut-off function (because the
lines are shifted but their scale labels are kept fixed so that scales
of the propagators of the self-energy clusters $T'\in\FF_{T}$ are the
ones inherited by $T$ while the momentum flowing in them may change).

This generates trees and clusters for which we made no estimates
(because they are just $0$). However when interpolating we may end up
computing values of trees, with scale assignments which would give a value
$0$, at {\it intermediate} frequencies where the values no longer
vanish. In estimating such interpolated values we can proceed as in
the cases already treated, but it will not be necessarily true
that a line of frequency $x$ and scale $[n]$
will satisfy $2^{-2(n-1)} C^2_0<D(x)\le 2^{-2n} C^2_0$.
Nevertheless a slightly weaker version of this inequality has to
hold in which the the {\it l.h.s.} is divided by $4$ and the {\it
r.h.s.}  is multiplied by $4$  (cf. also Remark (3) after the inductive
assumption in Section \secc(6)), and the estimates will not only be
possible but they can be regarded as already obtained because, 
as the reader can check, we have been careful in discussing
the bounds obtained so far under such weaker condition.
This also clarifies why we have defined $\lis n$ in \equ(4.2)
one unit larger than what appeared there as necessary so that the
estimate \equ(4.4) is apparently worse than it should.

In some cases, however, a serious problem seems to arise when actually
attempting to derive bounds: namely the bounds on the matrices which
appear as coefficients in \equ(A3.6) can really be checked as just
outlined by the above hints, and without affecting the values of $\e$
for which one has convergence, only if $x$ verifies the condition of
being so small that the variations of the momenta flowing in the inner
lines of $T$, when the entering or exiting lines are moved and
re-attached to all nodes of $T$, remain so small that the quantities
$D(x_{\ell})$ corresponding to the lines $\ell$ in the cluster $T$
stay essentially unchanged.

In certain cases shifting the entering or exiting lines to the nodes
of the self-energy cluster $T$ may considerably change the scales of
the lines $\ell$ in $T$, but this is the case in which
\equ(A3.7) holds. And precisely in such a case the
cancellations are not needed to prove the bound, because we have
checked that the value of {\it each} self-energy cluster
contributing to $M^{[n]}$ individually already
verifies that bound that we want to prove.

If \equ(A3.7) does not hold, then two  cases are possible:
either $|x|$ is close to $\l_j^{[p]}$
for some $j>r$ or larger, and no cancellation occurs,
or $|x|$ is $< C_0 2^{-n} $.
In the latter case the inequality opposite to \equ(A3.7)
implies that for $\ell\in\L(T)$ one has
$|x^0_\ell|\ge 4|x|$, hence $2|x^0_\ell|\ge |x_\ell|\ge \fra12
|x^0_\ell|$, so that the scales can change by at most one unit by 
shifting the external lines of $T$. Then the quantities $D(x_{\ell})$
do not change much for all lines $\ell \in \L(T)$, and we shall have the
cancellation through the mentioned mechanism. Therefore the contribution
of $\MM^{[p]}(x;\e)$ to $\MM^{[\le n]}(x;\e)$ can be bounded in both cases
proportionally to $e^{-\fra14 \k_0 2^{(n-\lis n-4 )/\t_0}}$ times
$\min\{\e^2,\e\,|x|^2\}$ for the entries $\a\a$ or times
$\min\{\e^2,\e^{\fra32} |x|\}$ for the $\a\b$ entries:
either by the cancellation (second case) or by the general
bound $O(\e^2)$ on matrix elements (first case),
because $x^2$ is of order $O(\e)$.

Finally we note that in the estimates of the $M$'s in \equ(A3.6) we
have to sum over the scale labels and this gives a factor per line
larger than the one arising in the bound \equ(A3.3) (which was $2$);
in fact we have to consider also trees with vanishing value: but the
scales of the divisors associated with their lines can change at most
by one unit with respect to the scale, hence we can have at most $4$ scale
labels per line.

\*

\0{\it Remark.}
We stress once more that the above analysis holds if $\e$ is small
enough, say $\e<\lis \e_{1}$ with $\lis\e_{1}$ determined by
collecting all the (three)
restrictions imposed by requiring $\e$ to be ``small enough'', derived
above and $\lis\e_1$ is {\it independent} of $n_0$ (otherwise it would
be uninteresting). The reason is that as long as we do not deal with
$x$'s which are too close to the eigenvalues of $M_0$, so that the key
inequality \equ(4.4) holds, we do not really see the difference
between the hyperbolic and the elliptic cases: and in the hyperbolic
cases there is no need for a lower cut-off at scale $\sim n_0$ where
resonances between the proper frequencies (which are 
of order $\e$) and the elliptic normal
frequencies become possible (as $\e \simeq C^2_{0} 2^{-2n_0}$).

\*

\asub(A3.4) {\it Resonant resummations.} Concerning the proof of Lemma 6 we
only need to add a few comments.
The bounds on $\NN_{m}(\th)$ and $\NN_{m}(T)$ can be discussed
exactly as for the scales $[n]$ with $n \le \lis n_{0}$,
with the only difference that now one has to use also
the second part of the Diophantine conditions \equ(6.7),
as already done in the argument leading to \equ(6.12);
in particular the role of the exponent $\t_{0}$ is now
played by $2\t_{1}$ (because of the Diophantine conditions
in \equ(6.7) which replaces \equ(1.3) in the discussion),
while in the analogues of \equ(A3.1)
and the following bounds no $\lis n$ appear, as the
propagator divisors are bounded directly in terms of the
corresponding scales, and not in terms of the frequencies.

Also the argument given above about the cancellations extends easily
to the scales $[n]$, with $n \ge \lis n_{0}$. The only difference
is that in \equ(A3.7) the exponent $1/\t_{0}$ has to
be replaced with $1/(2\t_{1})$, in such a way that
for any line $\ell\in \L(T)$ one has $||x_{\ell}^{0}|-
\Sqrt{|\l^{[n_{\ell}-1]}_{j}(\e)|}| \ge 4|x|$, hence
the chain of inequalities
$$ 2 \left| |x_{\ell}^{0}|-\Sqrt{|\l^{[n_{\ell}-1]}_{j}(\e)|} \right|
\ge \left| |x_{\ell}|-\Sqrt{|\l^{[n_{\ell}-1]}_{j}(\e)|} \right| \ge
\fra{1}{2} \left| |x_{\ell}^{0}|-
\Sqrt{|\l^{[n_{\ell}-1]}_{j}(\e)|} \right| ,
\Eqa(A3.8) $$
follows, and again by shifting the external lines of $T$
the scales of the internal lines can change at most by one unit,
when \equ(A3.7) is not satisfied

\*\*
\appendix(A4,Matrix properties)

\0{\bf (I)}
{\it Let $M_0$ be a $d\times d$ Hermitian matrix with eigenvalues
$\l_1<\ldots<\l_p$ with multiplicities $n_1,\ldots,n_p$ and
eigenspaces $\P_1,\ldots,\P_p$ on which we fix orthonormal bases $\ul
e_{j,k}, j=1,\ldots,p, k=1,\ldots,n_j$. Let $M_1$ be Hermitian and
consider the matrix $M=M_0+\e M_1$. There exists a constant $C$ such
that, for $\e$ small enough,
there will be $n_j$ eigenvalues of $M$ (not necessarily all distinct)
which are analytic in $\e$ and one has
$|\l_{j,k}(\e)-\l_{j,k'}(\e)|\le C\e$
for $k,k'=1,\ldots,n_j$.}

\*

\0{\it Hint.} If $n_{j}=1$ this follows immediately form the formula
$\l_j(\e)={\rm Tr}\,\Big( \fra1{2\p i}\oint_{\g_j} \fra{z
dz}{z-M}\Big)$, where $\g_j$ is a circle around $\l_j(0)$ of
$\e$-independent radius smaller than half the minimum separation $\d$
between the $\l_j$ for $\e$ small enough (so that $C_1 \e^{\fra1d}<\d
$ for a suitable $C_1$)\footnote{${}^{10}$}{\nota Because the
characteristic polynomials $P(\l),P_0(\l)$ are
related by $P(\l)=P_0(\l)+\e Q(\l,\e)$ with $Q$ of lower degree.
Therefore there is $L$ such that if $|\l|>L$ then for
all $|\e|<1$ (say) it is $P(\l)\ne0$. Furthermore if all roots of $P$
differ by at least $y$ from those of $P_0$ one has $|P(\l)|\ge y^d-\e
C^d$ where $C^d=\max_{|\L|\le L, |\e|\le 1}|Q(\l,\e)$. Hence $y\le
C\e^{d^{-1}}$.\vfil}.

Otherwise it follows form similar formulae for the projection operator
$E_j$ on $\P_j$ and for $E_jME_j$:
$$ E_j=\fra1{2\p i}\oint_{\g_j} \fra{dz}{z-M},\qquad
E_jM E_j=\fra1{2\p i}\oint_{\g_j} E_j\fra{z\, dz}{z-M}E_j ,
\Eqa(A4.1)$$
which, for $\e$ small, can be expanded into a convergent power series
in $\e$ (as done explicitly in a similar context in \equ(A4.3) below)
because of the $\e$-independence of the radii of $\g_j$.
One can also construct an orthonormal basis on $\P_j$
with vectors of the form ${\bf v}_{j,k}={\bf e}_{j,k}+\sum_{q\ge 1}^\io
\e^q {\bf e}^{(q)}_{j,k}$ (simply applying the Hilbert-Schmidt
orthonormalization to the vectors $E_j {\bf e}_{j,k}, k=1,\ldots,n_j$).
One then remarks that the matrix $E_jME_j$ has $n_j$ eigenvalues and
that it has the form $\l_j+\e \widetilde M(\e)$.

So the problem is reduced to the case in which $M_0$ is the identity
perturbed by an analytic matrix. Either $\widetilde M(\e)$ is
proportional to the identity and there is nothing more to do, or it is
not: hence there will be an order in $\e$ at which the degeneracy is
removed and repeating the argument we reduce the problem to a similar
one for matrices of dimension lower than $n_j$: and so on until we
find a matrix (possibly one dimensional)
proportional to the identity to all orders. \qed

\*

In our analysis we need the following corollary.

\*

\0{\bf (II)} {\it Let $M_0$ be Hermitian with $r$ degenerate
eigenvalues equal to $0$ and $s=d-r$ simple eigenvalues $\e a_j$,
$j=1,\ldots,s$.
\\
(i) The matrix $M_0+\e^2 M_1$ with $M_1$ Hermitian and differentiable in
$\e$ with bounded derivative will have $s$ non-degenerate eigenvalues
$\e a_j+O(\e^2)$, $j=1,\ldots,s$, and $r$ eigenvalues $\l_{1}(\e),
\ldots,\l_{r}(\e)$, all analytic in $\e$, with the property
that for all $k=1,\ldots,r$ one has $|\l_k(\e)|< C\, \e^2$,
if $\e$ is small enough and $C$ is a suitable constant.
\\
(ii) If $M_1$ depends on a parameter $x$ and is differentiable
also in $x$ with bounded derivative then
$$\eqalign{
& |\dpr_x \l_j(x;\e)|\le C\e^2, \quad|\dpr_\e\l_j(x;\e)|\le C ,
\qquad j>r ,\cr
&|\l_j(x;\e)-\l_{j}(x';\e)|\le C \e^2 |x-x'|^{1/r} ,
\qquad j\le r , \cr}
\Eqa(A4.2)$$
if $\e$ is small enough and $C$ is a suitable constant.}

\*
The second relation in \equ(A4.2) is not used in this paper and is
given only for completeness.
\*

\0{\it Hint.} We apply the previous lemma to the matrices $\e^{-1}M_0$
and $\e M_1$ and we get (i). To get \equ(A4.2) we note that the
$x$-derivative of $M_0+\e^2 M_1$ is $\e^2\dpr_x M_1$ and the first of
\equ(A4.2) follows. To obtain the second of \equ(A4.2) we
have to compare the eigenvalues of $M_0+\e^2 M_1(x;\e)$ with those
of $M_0+\e^2 M_1(x;\e)+ \e^2 O(|x-x'|)$. By the above expression for
the projection on the plane of the first $r$ eigenvalues this is
reduced to the problem of comparing two $r\times r$ matrices of order
$\e^2$ and differing by $O(|x-x'|)$.
The power $1/r$ arises from the estimate that the considered 
projection of the matrix $M_1$
(which is only differentiable in $x$)
has $r$ eigenvalues close to $0$ within $C_1 |x-x'|^{1/r}$,
for some $C_1>0$ (by (I) above), and $\e$ is small enough.
Hence we get the second of \equ(A4.2). \qed

\*

A third property that we need is the following one.

\*

\0{\bf (III)} {\it If $M_0$ is as in (II) and $M_1$ is Hermitian and has
the form $\pmatrix{\e^2 x^2 N& \e^2 x P\cr \e^2 x P^*& \e^2 Q\cr}$,
with $N,Q$ two $r\times r$ and $r\times s$ matrices and $P$ a
$r\times s$ matrix then the first $r$ eigenvalues of $M_0+M_1$ are
bounded by $|\l_j(x,\e)|< C\e^2 x^2$, for $j=1,\ldots,r$.}

\*

\0{\it Hint.} This is obtained by using \equ(A4.1) which
gives the projection over the plane of the $r$ eigenvalues within
$O(\e^2)$ of $0$ as integral over a circle of radius $\fra12 a_1 \e$
$$ E=\fra1{2\p i}\oint_\g \fra{dz}{z-M_0}\sum_{k=0}^\io
(M_1\fra1{z-M_0})^k ,
\Eqa(A4.3) $$
and one sees that $(M_1\fra1{z-M_0})^k$ has for all $k\ge1$ the same
form of $M_1$, with $\e^{2}$ replaced by $\e^{2k}$,
so that the sum of the series is the matrix
$\pmatrix{1&0\cr0&0\cr}$ corresponding to the $k=0$ term (it is a
$d\times d$ block matrix with the first
$r\times r$ block $1$ and the other blocks $0$) plus a matrix
of the same form of $M_1$. Likewise the basis ${\bf v}_h=E{\bf e}_{h}$,
$h=1,\ldots,r$ consists of vectors of the form ${\bf e}_{h}+
\pmatrix{\e^2 x^2 {\bf u}_h\cr \e^2 x {\bf u}'_h}$, so that
one checks that the matrix $({\bf v}_h,(M_0+M_1) {\bf v}_{h'})$
is a $r\times r$ matrix which is proportional to $\e^2 x^2$
(\ie it has the form $\e^2 x^2 M_2(x,\e)$, with $M_2$
bounded for $\e$ small and for $|x|<1$) and
which, by construction, has the same eigenvalues as the first $r$
eigenvalues of the matrix $M_0+M_1$.

\*

For the above properties see also \cita{RS} and \cita{Ka}.

\*\*
\appendix(A5, Algebraic identities for the renormalized expansion)

\0We show that the function $\hh$ defined through the renormalized
expansion solves the equations of motion \equ(1.5) for all
$\e\in\EE$. This is essentially a repetition of Ref. \cita{Ge}.
We shall check that $\hhh =\e g \dpr_{\ff} f(\pps +\aaa, \bb_{0}+\bbb)$,
where $\ff=(\aa,\bb)$ and $g$ is the pseudo-differential operator
with kernel $g(\oo\cdot\nn)=(\oo\cdot\nn)^{-2}$.
We can write $\hhh = \sum_{\nn\in \zzz^{r}}
e^{i\nn\cdot\pps} \hhh_{\nn}$, $\hhh_{\nn} =
\sum_{n=0}^{\io} \hhh_{n,\nn}$ (only two terms in this series are
different from $0$ for each $\nn$), with $\hhh_{n,\nn} =
\sum_{k=1}^{\io} \sum_{\th \in \Th^{\RR}_{k,\nn}(n)} \Val(\th)$,
where $\Th^{\RR}_{k,\nn}(n)$ is the set of trees
in $\Th^{\RR}_{k,\nn}$ such that the root line has scale $n$.
With respect to the previous sections we have dropped the
component label $\g\in\{1,\ldots,d\}$ in the definition of the
set of trees, for notational convenience.

Note that, for all $x\neq 0$ and for all $p\ge 0$ one has
$$ 1 = \sum_{n=p}^{\io}
\psi_n(\D^{[n]}(x,\e)) \prod_{q=p}^{n-1} \chi_{q}(\D^{[q]}(x,\e)) ,
\Eqa(A5.1) $$
where the term with $n=p$ has to be interpreted
as $\psi_{p}(\D^{[p]}(x;\e))$. 

Set $\Psi_{n}(x;\e)=
\psi_{n}(\D^{[n]}(x;\e))\prod_{p=0}^{n-1}\chi_{p}( \D^{[p]}(x;\e))$
for $n\ge 1$, $\Psi_{0}(x;\e)=\psi_{0}(\D^{[0]}(x;\e))$: 
by using \equ(A5.1) one can write, in Fourier space and evaluating
the functions of $\ff$ at $\ff=(\pps+\aa,\bb_0+\bb)$,
$$ \eqalign{
& g(\oo\cdot\nn) \left[ \e \dpr_{\ff} f(\ff) \right]_{\nn}
= g(\oo\cdot\nn) \sum_{n=0}^{\io} \Psi_{n}(\oo\cdot\nn;\e)
\left[ \e \dpr_{\ff} f(\ff) \right]_{\nn} \cr
& \qquad \qquad
= g(\oo\cdot\nn) \sum_{n=0}^{\io} \Psi_{n}(\oo\cdot\nn;\e)
(g^{[n]}(\oo\cdot\nn;\e))^{-1}
g^{[n]}(\oo\cdot\nn;\e) \left[ \e \dpr_{\ff} f(\ff) \right]_{\nn} \cr
& \qquad \qquad
= g(\oo\cdot\nn) \sum_{n=0}^{\io}
\left( (\oo\cdot\nn)^{2} -\MM^{[\le n]}(\oo\cdot\nn;\e) \right)
g^{[n]}(\oo\cdot\nn;\e) \left[ \e \dpr_{\ff} f(\ff) \right]_{\nn} \cr
& \qquad \qquad
= g(\oo\cdot\nn) \sum_{n=0}^{\io}
\left( (\oo\cdot\nn)^{2} -\MM^{[\le n]}(\oo\cdot\nn;\e) \right)
\sum_{k=1}^{\io}
\sum_{\th \in \lis \Th^{\RR}_{k,\nn}(n)} \Val(\th) , \cr}
\Eqa(A5.2) $$
where $\lis\Th^{\RR}_{k,\nn}(n)$ differs from $\Th^{\RR}_{k,\nn}(n)$
{\it as it contains also trees which can have one
renormalized self-energy cluster $T$ with exiting line $\ell_{0}$},
if $\ell_{0}$ denotes the root line of $\th$;
for such trees the line entering $T$ will be on a scale $p\ge 0$,
while the renormalized self-energy cluster $T$ will have
a scale $n_{T}=q$, with $q+1 \le \min\{n,p\}$.

\*

\0{\it Remark.} Note that in both \equ(A5.1) and \equ(A5.2) only
a finite number of addends is different from zero,
as the analysis of Section \secc(6) shows,
so that the two series are well defined.
The same observation applies to the following formulae,
where appear series which, in fact, are finite sums.

\*

By explicitly separating in \equ(A5.2)
the trees containing such self-energy clusters from the others,
$$ \eqalignno{
& g(\oo\cdot\nn) \left[ \e \dpr_{\ff} f(\ff) \right]_{\nn}
= g(\oo\cdot\nn) \sum_{n=0}^{\io}
\left( (\oo\cdot\nn)^{2} -\MM^{[\le n]}(\oo\cdot\nn;\e) \right)
\sum_{k=1}^{\io} \sum_{\th \in \Th^{\RR}_{k,\nn}(n)} \Val(\th) \cr
& \qquad \qquad +
g(\oo\cdot\nn) \sum_{n=1}^{\io}
\left( (\oo\cdot\nn)^{2} -\MM^{[\le n]}(\oo\cdot\nn;\e) \right)
g^{[n]}(\oo\cdot\nn;\e) & \eqa(A5.3) \cr
& \qquad\qquad\qquad
\sum_{p=n}^{\io} \sum_{q=0}^{n-1}
M^{[q]}(\oo\cdot\nn;\e)
\sum_{k=1}^{\io}
\sum_{\th \in \Th^{\RR}_{k,\nn}(p)} \Val(\th) \cr
& \qquad \qquad +
g(\oo\cdot\nn) \sum_{n=2}^{\io}
\left( (\oo\cdot\nn)^{2} -\MM^{[\le n]}(\oo\cdot\nn;\e) \right)
g^{[n]}(\oo\cdot\nn;\e) \cr
& \qquad\qquad\qquad
\sum_{p=1}^{n-1} \sum_{q=0}^{p-1}
M^{[q]}(\oo\cdot\nn;\e)
\sum_{k=1}^{\io}
\sum_{\th \in \Th^{\RR}_{k,\nn}(p)} \Val(\th) , \cr} $$
which, by the definitions of $\hh$, can be written as
$$ \eqalignno{
& g(\oo\cdot\nn) \left[ \e \dpr_{\ff} f(\ff) \right]_{\nn}
= g(\oo\cdot\nn) \Big[ \sum_{n=0}^{\io}
\left( (\oo\cdot\nn)^{2} -\MM^{[\le n]}(\oo\cdot\nn;\e) \right)
\hhh_{n,\nn}
& \eqa(A5.4) \cr
& \qquad
+ \sum_{n=1}^{\io} \Psi_{n}(\oo\cdot\nn;\e)
\sum_{p=n}^{\io} \sum_{q=0}^{n-1}
M^{[q]}(\oo\cdot\nn;\e) \hhh_{p,\nn} +
\sum_{n=2}^{\io} \Psi_{n}(\oo\cdot\nn;\e)
\sum_{p=1   }^{n-1} \sum_{q=0}^{p-1}
M^{[q]}(\oo\cdot\nn;\e) \hhh_{p,\nn}
\Big] . \cr} $$
The terms in the second line of \equ(A5.4) can be written as
$$ \eqalign{
& \sum_{p=1}^{\io} \Big(
\sum_{q=0}^{p-1} \sum_{n=q+1}^{p}
M^{[q]}(\oo\cdot\nn;\e) \Psi_{n}(\oo\cdot\nn;\e) +
\sum_{q=0}^{p-1} \sum_{n=p+1}^{\io}
M^{[q]}(\oo\cdot\nn;\e) \Psi_{n}(\oo\cdot\nn) \Big) \hhh_{p,\nn} \cr
& \qquad = \sum_{p=1}^{\io} 
\sum_{q=0}^{p-1} M^{[q]}(\oo\cdot\nn;\e)
\sum_{n=q+1}^{\io} \Psi_{n} (\oo\cdot\nn;\e) \, \hhh_{p,\nn} \cr}
\Eqa(A5.5) $$
and, by changing $p\to n$ and $n\to s$, we obtain
$$ \eqalign{
& \sum_{n=1}^{\io} \Big( 
\sum_{q=0}^{n-1} M^{[q]}(\oo\cdot\nn;\e)
\chi_{0}(\D^{[0]}(\oo\cdot\nn;\e)) \ldots
\chi_{q}(\D^{[q]}(\oo\cdot\nn;\e)) \, \cdot \cr
& \qquad \qquad \qquad
\cdot \sum_{s=q+1}^{\io} \chi_{q+1}(\D^{[q+1]}(\oo\cdot\nn;\e)) \ldots
\psi_{s}(\D^{[s]}(\oo\cdot\nn;\e)) \Big) \hhh_{n,\nn} \cr
& \qquad = \sum_{n=1}^{\io} 
\sum_{q=0}^{n-1} M^{[q]}(\oo\cdot\nn;\e)
\chi_{0}(\D^{[0]}(\oo\cdot\nn;\e)) \ldots \chi_{q}
(\D^{[q]}(\oo\cdot\nn;\e)) \, \hhh_{n,\nn} , \cr}
\Eqa(A5.6) $$
where the identity \equ(A5.1) has been used in the last line
(with the correct interpretation of the term with $s=j+1$
explained after \equ(A5.1)).
By the definition of the matrices $\MM^{[\le n]}(x;\e)$ one has
$$ \sum_{q=0}^{n-1}
M^{[q]}(\oo\cdot\nn;\e) \chi_{0}(\D^{[0]}(x;\e)) \ldots
\chi_{q}(\D^{[q]}(x;\e)) = \MM^{[\le n]}(x;\e) ,
\Eqa(A5.7) $$
so that, by inserting \equ(A5.6) in \equ(A5.3),
after having used \equ(A5.7), we obtain
$$ \eqalignno{
g(\oo\cdot\nn) \left[ \e \dpr_{\ff} f(\ff) \right]_{\nn}
& = g(\oo\cdot\nn) \sum_{n=0}^{\io} \Big[
\left( (\oo\cdot\nn)^{2} - \MM^{[\le n]}(\oo\cdot\nn;\e) \right) +
\MM^{[\le n]}(\oo\cdot\nn;\e) \Big] \hhh_{n,\nn} \cr
& = g(\oo\cdot\nn) \sum_{n=0}^{\io} (\oo\cdot\nn)^{2} \hhh_{n,\nn}=
\sum_{n=0}^{\io} \hhh_{n,\nn} = \hhh_{\nn} ,
& \eqa(A5.8) \cr} $$
so that the assertion is proved.

\*

\0{\it Remark.}
Note that at each step only absolutely converging series
have been dealt with, so that the above analysis is
rigorous and not only formal.

\*\*
\0{\bf Acknowledgments.} We are indebted to V. Mastropietro for
many discussions and, in particular, to A. Giuliani for
critical reading and several suggestions.

\*\*

{\nota

\rife{B1}{1}{
J. Bourgain,
{\it Construction of quasi-periodic solutions for Hamiltonian
perturbations of linear equations and applications to nonlinear PDE},
Internatational Mathematics Research Notices
{\bf 11} (1994) 475--497. }
\rife{B2}{2}{
J. Bourgain,
{\it Construction of periodic solutions of nonlinear
wave equations in higher dimension},
Geometric and Functional Analysis
{\bf 5} (1995), 629--639. }
\rife{B3}{3}{
J. Bourgain,
{\it On Melnikov's persistency problem},
Mathematical Research Letters
{\bf 4} (1997), 445--458. } 
\rife{B4}{4}{
J. Bourgain,
{\it Quasi-periodic solutions of Hamiltonian perturbations of 2D
linear Schr\"odinger equations},
Annals of Mathematics
{\bf 148} (1998), no. 2, 363--439. }
\rife{Ba}{5}{
J. C. A. Barata,
{\it On formal quasi-periodic solutions of the Schr\"odinger
equation for a two-level system with a Hamiltonian
depending quasi-periodically on time},
Reviews in Mathematical Physics
{\bf 12} (2000), no. 1, 25--64. }
\rife{BaG}{6}{M.V. Bartuccelli, G. Gentile,
{\it Lindstedt series for perturbations of isochronous systems.
A review of the general theory},
Reviews in Mathematical Physics
{\bf 14} (2002), no. 2, 121--171.}
\rife{BGGM}{7}{\ifnum\tipobib=2\hskip3mm\fi%
F. Bonetto, G. Gallavotti, G. Gentile, Mastropietro,
{\it Lindstedt series, ultraviolet divergences and Moser's theorem},
Annali della Scuola Normale Superiore di Pisa Classe di Scienze
{\bf 26} (1998), no. 3, 545--593. }
\rife{C}{8}{
C.-Q. Cheng,
{\it Lower-dimensional invariant tori in the regions of instability
for nearly integrable Hamiltonian systems},
Communications in Mathematical Physics
{\bf 203} (1999), no. 2, 385--419. }
\rife{ChW}{9}{
C.-Q. Cheng, S. Wang,
{\it The surviving of lower-dimensional tori from a resonant
torus of Hamiltonian systems},
Journal of Differential Equations
{\bf 155} (1999), no. 2, 311--326. }
\rife{CrW}{10}{
W. Craig and C.E. Wayne,
{\it Newton's method and periodic solutions of nonlinear
wave equations},
Communications in Pure and Applied Mathematics
{\bf 46} (1993), no. 11, 1409--1501. }
\rife{E1}{11}{
L.H. Eliasson,
{\it Perturbations of stable invariant tori for Hamiltonian systems},
Annali della Scuola Normale Superiore di Pisa Classe di Scienze
{\bf 15} (1988), no. 1, 115--147. }
\rife{E2}{12}{
L.H. Eliasson,
{\it Absolutely convergent series expansions for quasi-periodic motions},
Mathematical Physics Electronic Journal
{\bf 2} (1996), paper 4, 1--33 (Preprint 1988). }
\rife{Fe}{13}{
C. Fefferman,
{\it Pointwise convergence of Fourier series},
Annals of Mathematics
{\bf 98}, 551--571, 1973.}
\rife{G3}{14}{
G. Gallavotti, {\it Twistless KAM tori, 
quasi flat homoclinic intersections, and
other cancellations in the perturbation series of certain completely
integrable hamiltonian systems. A review},
Reviews in Mathematical Physics
{\bf 6}, 343--411, 1994.}
\rife{Ga}{15}{
G. Gallavotti,
{\it Twistless KAM tori},
Communications in Mathematical Physics
{\bf 164} (1994), no. 1, 145--156. }
\rife{Ga01}{16}{
G. Gallavotti: 
{\it Renormalization group in Statistical Mechanics and Mechanics:
gauge symmetries and vanishing beta functions},
Physics Reports
{\bf352}, 251--272, 2001.}
\rife{Ga02}{17}{
G. Gallavotti: {\it Exact Renormalization Group}, 
Paris IHP, 12 october 2002, Seminaire Poincar\'e, Editors
B. Duplantier, V. Rivasseau, Institut H. Poincar\'e-CNRS-CEA.}
\rife{GG}{18}{
G. Gallavotti, G. Gentile,
{\it Hyperbolic low-dimensional tori and summations of divergent series},
Communications in Mathematical Physics
{\bf 227} (2002), no. 3, 421--460. }
\rife{GBG}{19}{
G. Gallavotti, F. Bonetto, G. Gentile,
{\it Aspects of the ergodic, qualitative and statistical properties of motion}
Springer--Verlag, Berlin, 2004.}
\rife{Ge}{20}{
G. Gentile,
{\it Quasi-periodic solutions for two-level systems},
Communications in Mathematical Physics
{\bf 242} (2003), no. 1--2, 221--250. }
\rife{GM1}{21}{
G. Gentile, V. Mastropietro,
{\it Methods for the analysis of the Lindstedt series for KAM tori
and renormalizability in classical mechanics.
A review with some applications},
Reviews in Mathematical Physics
{\bf 8} (1996), no. 3, 393--444. }
\rife{GM2}{22}{
G. Gentile, V. Mastropietro,
{\it  Construction of periodic solutions of nonlinear wave equations
with Dirichlet boundary conditions by the Lindstedt series method},
to appear on
Journal de Math\'ematiques Pures et Appliqu\'ees. }
%
\rife{GMP}{23}{
G. Gentile, V. Mastropietro, M. Procesi,
{\it Periodic solutions for completely resonant wave equations},
Preprint, 2004. }
\rife{JLZ}{24}{
A. Jorba, R. de la Llave, M. Zou,
{\it Lindstedt series for lower-dimensional tori}, in
{\it Hamiltonian systems with three or more degrees of
freedom} (S'Agar\'o, 1995), 151--167,
NATO Adv. Sci. Inst. Ser. C Math. Phys. Sci., 533, Ed. C. Sim\'o,
Kluwer Acad. Publ., Dordrecht, 1999. }
\rife{Ka}{25}{
T. Kato,
{\it Perturbation theory for linear operators},
Grundlehren der Mathematischen Wissenschaften, Band 132,
Springer-Verlag, Berlin-New York, 1976. }
\rife{Ku}{26}{
S.B. Kuksin,
{\it Nearly integrable infinite-dimensional Hamiltonian systems},
Lecture Notes in Mathematics 1556,
Springer-Verlag, Berlin, 1993.}
\rife{KP}{27}{
S.B. Kuksin, J. P\"{o}schel,
{\it Invariant Cantor manifolds of quasi-periodic oscillations
for a nonlinear Schr\"odinger equation},
Annals of Mathematics
{\bf 143} (1996), no. 1, 149--179. }
\rife{LW}{28}{
R. de la Llave, C.E. Wayne,
{\it Whiskered and low dimensional tori in nearli integrable 
Hamiltonian systems},
Mathematical Physics Electronic Journal, 2004. }
%
\rife{Me1}{29}{
V.K. Mel'nikov,
{\it On certain cases of conservation of conditionally periodic motions
under a small change of the Hamiltonian function},
Doklady Akademii Nauk SSSR 
{\bf 165} (1965), 1245--1248;
english translation in
Soviet Mathematics Doklady 
{\bf 6} (1965), 1592--1596. }
\rife{Me2}{30}{
V.K. Mel'nikov,
{\it A certain family of conditionally periodic solutions of a
Hamiltonian systems},
Doklady Aka\-de\-mii Nauk SSSR 
{\bf 181} (1968), 546--549;
english translation in
Soviet Mathematics Doklady
{\bf 9} (1968), 882-886. }
\rife{Mo}{31}{
J. Moser,
{\it Convergent series expansions for quasi-periodic motions},
Mathematische Annalen
{\bf 169} (1967), 136--176. }
\rife{P1}{32}{
J. P\"{o}schel,
{\it On elliptic lower-dimensional tori in Hamiltonian systems},
Mathematische Zeitschrift 
{\bf 202} (1989), no. 4, 559--608. }
\rife{P2}{33}{
J. P\"oschel,
{\it Quasi-periodic solutions for a nonlinear wave equation},
Commentarii Mathematici Helvetici
{\bf 71} (1996), no. 2, 269--296. }
\rife{RS}{34}{
M. Reed, B. Simon,
{\it Methods of modern mathematical physics. IV. Analysis of operators},
Academic Press, New York-London, 1978. }
\rife{R}{35}{
H. R\"ussmann,
{\it Invariant tori in non-degenerate nearly
integrable Hamiltonian systems}, 
Regular and Chaotic Dynamics
{\bf 6} (2001), 119--204. }
\rife{T}{36}{
D.V. Treshch\"ev,
{\it A mechanism for the destruction of resonance tori
in Hamiltonian systems},
Rossi\u\i skaya Akademiya Nauk. Matematicheski\u\i Sbornik 
{\bf 180} (1989), no. 10, 1325--1346;
english translation in
Mathematics of the USSR-Sbornik
{\bf 68} (1991), no. 1, 181--203. }
\rife{WC}{37}{
S. Wang, C.-Q. Cheng,
{\it Lower-dimensional tori for generic Hamiltonian systems},
Chinese Science Bulletin
{\bf 44} (1999), no. 13, 1187--1191. }
\rife{Wa}{38}{
C.E. Wayne,
{\it Periodic and quasi-periodic solutions of nonlinear wave equations
via KAM theory},
Communications in Mathematical Physics
{\bf 127} (1990), no. 3, 479--528. }
\rife{XY}{39} {
J. Xu, J. You:
{\it Persistence of lower dimensional 
tori under the first Melnikov's non-resonance condition},
Journal de Math\'ematiques Pures et Appliqu\'ees
{\bf 80} (2001), no. 10, 1045--1067. }
\rife{Y}{40}{
J. You:
{\it Perturbations of lower-dimensional tori for Hamiltonian systems},
Journal of Differential Equations
{\bf 152} (1999), no. 1, 1--29. }

\biblio
}

\end